\documentclass[siam11,leqno,onefignum,onetabnum]{siamltex}
\usepackage{amsmath}
\usepackage{amsfonts}
\usepackage{graphicx}
\graphicspath{{figures/}}
\def\d{\delta}
\def\l{\lambda}

\def\T{\mathbb T}
\def\Z{\mathbb Z}
\def\R{\mathbb R}
\def\cH{\mathcal H}
\def\cF{\mathcal F}

\title{Ergodic problems\\ for Hamilton-Jacobi equations:\\yet another but efficient numerical method}

\author{Simone Cacace and Fabio Camilli}

\begin{document}
\maketitle
\slugger{sisc}{xxxx}{xx}{x}{x--x}

\begin{abstract}
We propose a new approach to the numerical solution of ergodic problems arising in the homogenization of Hamilton-Jacobi (HJ) equations.
It is based on a Newton-like method for solving inconsistent systems of nonlinear equations,
coming from the discretization of the corresponding ergodic HJ equations.  We show that our method is
able to solve efficiently cell problems in very general contexts, e.g., for first and second order scalar convex and nonconvex Hamiltonians,
weakly coupled systems, dislocation dynamics and mean field games, also in the case of more competing populations. 
A large collection of numerical tests in dimension one and two shows the performance of the proposed method, both in terms of accuracy and computational time.
\end{abstract}

\begin{keywords}
Hamilton-Jacobi equations, homogenization, effective Hamiltonian, Newton-like methods, inconsistent nonlinear systems, dislocation dynamics, mean field games.
\end{keywords}

\begin{AMS}35B27, 35F21, 49M15.\end{AMS}

\pagestyle{myheadings}
\thispagestyle{plain}

\section{Introduction}\label{INTRO}
In many problems, such as  homogenization and long time behavior of first and second order Hamilton-Jacobi equations, weak KAM
theory, ergodic mean field games and dislocation dynamics, an essential step for the qualitative
analysis of the problem is the computation of the effective Hamiltonian. This function plays the
role of an eigenvalue and in general it is unknown except in some very special cases. Hence
the importance of  designing  efficient algorithms for its computation, taking also into account that the evaluation
at each single point of this function requires the solution of a nonlinear partial differential equation.
Moreover, the problem characterizing the effective Hamiltonian is in many cases ill-posed. Consider for example
the cell problem for a first order Hamilton-Jacobi equation
\begin{equation}\label{HJ}
 H(x, Du+p)=\l\,,\qquad x\in\T^n,
\end{equation}
where $\l\in\R$, $p\in\R^n$ and $\T^n$ is the unit $n$-dimensional torus.
It involves, for any given $p$, two unknowns $(u,\l)$ in a single
equation. Moreover, despite the effective Hamiltonian $\bar H(p):=\l$ is uniquely identified by \eqref{HJ},
the corresponding viscosity solution $u$ is in general not unique, not even for addition of constants.\par
In the recent years, several numerical schemes for the approximation of the effective Hamiltonian have been proposed 
(see \cite{ACC},\cite{Ba},\cite{E},\cite{FR},\cite{GO},\cite{Q},\cite{R06}), and they are mainly based on two different approaches. 

The first approach consists in the regularization of the cell problem \eqref{HJ} via well-posed problems, such as the stationary problem
$$
  \d u_\d +H(x,Du_\d+p)=0\,,\qquad  x\in \T^n,
$$
for $\d>0$, or the evolutive one
$$
  u_t+H(x,Du+p)=0\,,\qquad  x\in \T^n,\quad t\in (0,\infty).
$$
Indeed, it can be proved that both  $-\d u_\d$  and $-u(x,t)/t$ converge to $\bar H(p)$, respectively for $\d\to 0$ and $t\to +\infty$
(see \cite{LPV}). In \cite{Q}, these regularized problems are discretized by finite-difference schemes, obtaining respectively
the so called small-$\d$ and large-$t$ methods. In the limit for $\d\to 0$ or $t\to +\infty$, and simultaneously for the discretization step $h\to 0$,
one gets an approximation of $\bar H(p)$ (the convergence of this method is proved in \cite{ACC}).
Hence the computation of an approximation of  $\bar H(p)$, for any fixed $p$, requires the solution
of a sequence of nonlinear finite-difference systems,  which become more and more  ill-conditioned when $\d$ is small
 or $t$ is large. It is worth noting that the idea of approaching   ergodic problems via  small-$\d$ or  large-$t$ methods has been applied 
 to several other contexts, for example in mean field games theory \cite{AC}, $u$-periodic homogenization problems  \cite{AP} and dislocation dynamics \cite{CCM}.

The second approach for computing the effective Hamiltonian is based on the following $\inf$-$\sup$ formula:
\[\bar H(p)=\inf_{u\in C^\infty(\T^n)}\sup_{x\in\T^n} H(x, Du+p)\,.\]
In \cite{GO}, this formula is discretized on a simplicial grid, by taking the infimum on the subset of 
piecewise affine functions and the supremum on the barycenters of the grid elements. The resulting discrete problem is then solved via standard minimax methods.
An alternative method  is suggested in \cite{E}, where it is shown that the solution of the Euler-Lagrange equation
$$\mbox{div}\left(e^{k\,H(x,Du+p)} H_p(x,Du+p)\right)=0,\qquad  x\in \T^n$$
approximates, for $k\to+\infty$, the infimum in the $\inf$-$\sup$ formula. A finite difference implementation of this method is presented in \cite{FR}, 
for the special class of eikonal Hamiltonians.\\
\par
 In this paper we propose a new approach which allows  to compute   solutions of   ergodic problems {\em directly}, i.e.,
 avoiding  small-$\d$, large-$t$ or $\inf$-$\sup$ approximations. All these problems involve  a couple of unknowns $(u,\l)$, possibly depending  on some parameter,
 where $u$ is either a scalar or vector function and $\l$ is a constant.
 After performing a discretization of the ergodic problem (e.g., using finite-difference schemes), we   collect all the unknowns $(U,\Lambda)$
 of the discrete problem in a single
 vector $X$ of length $N$ and  we  recast  the $M$ equations   of the discrete system as functions  of $X$, for some $N,M\in\mathbb{N}$.
 We get a  nonlinear map $ F:\R^{N}\to\mathbb R^M$ and the discrete problem is equivalent to find $X^\star\in\R^N$ such that
\begin{equation}\label{I3}
 F(X^\star)=0\in\R^M\,,
\end{equation}
where the system \eqref{I3} can be inconsistent, e.g., {\em underdetermined}  ($M<N$) as for the cell problem \eqref{HJ},  or {\em overdetermined} ($M>N$)
as for stationary mean field games (see Section \ref{meanfieldgames}). Note that this terminology is properly employed for linear systems, but it is commonly adopted also for nonlinear systems.
In each case, \eqref{I3} can be solved  by a  generalized Newton's method involving the Moore-Penrose pseudoinverse of the Jacobian of 
$F$ (see \cite{BI}), and efficiently implemented
via suitable $QR$ factorizations.
\par
This approach has been experimented by the authors in the context of stationary mean field games on networks (see \cite{CCM16}) and,
to our knowledge, this is the first time a cell problem in homogenization of Hamilton-Jacobi
equations is solved directly, by interpreting the effective  Hamiltonian as an unknown (as it is!).
We realized that, once a consistent discretization of
the Hamiltonian is employed to correctly approximate viscosity solutions, all the job is reduced to computing zeros of nonlinear maps.
Moreover, despite the cell problem does not admit in general a unique viscosity solution, the ergodic constant defining the
effective Hamiltonian is often  unique. This ``weak'' well-posedness of the problems, and also the fact that the effective Hamiltonian is usually the main object of
interest more than the viscosity solution itself, encouraged the development of the proposed method.\par

The paper is organized as follows. In Section \ref{newtonlike} we introduce our approach for solving ergodic problems, 
and we present a Newton-like method for inconsistent nonlinear systems, discussing some basic features and implementation issues.
In the remaining sections, we apply the new method to more and more complex ergodic problems, arising in very different contexts.
More precisely, Section \ref{Eikonal} is devoted to the eikonal Hamiltonian, which is the benchmark for our algorithm, due to the availability of 
an explicit formula for the effective Hamiltonian. Section \ref{qnorms} concerns more general convex Hamiltonians, while Section \ref{nonconvex}
is devoted to a nonconvex case and Section \ref{2ndorder} to second order Hamiltonians. In Section \ref{WCS} we solve some vector problems for weakly coupled
systems, whereas Section \ref{dislocations} is devoted to a nonlocal problem arising in dislocation dynamics. Finally, in Section \ref{meanfieldgames}
we solve some stationary mean field games and in Section \ref{multipopMFG} an extension to the vector case with more competing populations.
\section{A Newton-like method for ergodic problems}\label{newtonlike}
In this section we introduce a new numerical approach for solving ergodic problems arising in the homogenization of Hamilton-Jacobi equations. Then we present
a Newton-like method for inconsistent nonlinear systems and we discuss some of its features, also from an implementation point of view.

We assume that a generic continuous ergodic problem is defined on the torus $\T^n$  and we denote by  $\T^n_h$ a numerical grid on $\T^n$.
We also assume that the  discretization  scheme for the continuous problem  results in  a system of nonlinear equations of the form
\begin{equation}\label{scheme}
S(x,h,U)=\Lambda
\end{equation}
where
\begin{itemize}
\item $h>0$ is the discretization parameter ($h$ is meant to  tend to $0$);
\item $x\in \T^n_h$ is the point where the continuous problem  is approximated;
\item $U$ is a real valued mesh function on $\T^n_h$ and $\Lambda$ is a real number,
meant to approximate respectively the continuous solution $u$ of the (possibly vector) ergodic problem and the corresponding ergodic constant $\l$;
\item $S$ represents a generic numerical scheme;
\item $N$ and $M$ are respectively the length of the vector $(U,\Lambda)$ and the number of equations in \eqref{scheme}.
\end{itemize}
\medskip
\noindent We remark that our aim is to efficiently solve the nonlinear system \eqref{scheme}.
Therefore, we do not specify the type of grids or schemes employed in the discretization. In particular, 
we do not consider properties of the scheme itself, such as consistency, stability, monotonicity and, more important,
the ability to correctly select approximations of viscosity solutions, assuming that they are included in the form of the operator $S$.
We just point out that, in our tests, we always perform finite difference discretizations on uniform grids. Moreover, if not differently specified,
we mainly employ the well known Engquist-Osher numerical approximation for the first order terms in the equations, due to its simple implementation.
For instance, in dimension one, we use the following upwind approximation of the gradient:
$$|Du|(x)\sim \sqrt{{\min}^2\left\{\frac{U(x+h)-U(x)}{h},0\right\} +{\max}^2\left\{\frac{U(x)-U(x-h)}{h},0\right\}}\,.$$

Here, the only assumption on the discrete ergodic problem is the following: 
$$
    \begin{array}{c}
    \text{for each fixed $h$, there exists a unique $\Lambda$ for which \eqref{scheme}}\\
    \text{ admits a solution $U$ (in general not unique). }
    \end{array}
$$
To compute a solution of \eqref{scheme}, we   collect   the unknowns $(U, \Lambda)$   in a single vector $X$ of length $N$  and we recast
the $M$ equations as functions of $X$. Hence  we get the nonlinear map
$ F:\mathbb{R}^{N}\to\mathbb{R}^{M}$ defined by $F(X)=S(x ,h,U)-\Lambda$, 
and \eqref{scheme} is equivalent to the nonlinear system
\begin{equation}\label{s:3}
 F(X)=0.
\end{equation}
The system \eqref{s:3} is said {\em underdetermined}  if $M<N$ and {\em overdetermined} if $M > N$.
As already remarked, this terminology applies to linear systems, nevertheless it is commonly adopted, with a slight abuse, also in the nonlinear case.

Assuming that $F$ is Fr\'echet differentiable, we consider the following generalized Newton-like method \cite{BI}: given $X^0\in\R^N$, iterate up to convergence
\begin{equation}\label{NL1}
    X^{k+1}=X^k-J_F(X^k)^\dag F(X^k),\qquad k\ge 0\,,
\end{equation}
where $J_F(\cdot)=\frac{\partial  F_i}{\partial X_j}(\cdot)$ is the Jacobian of $F$ and  $J_F(\cdot)^\dag$ denotes the Moore-Penrose {\em pseudoinverse} of $J_F(\cdot)$. 
As in the case of square systems, we can rewrite \eqref{NL1} in a form suitable for computations, i.e., for $k\ge 0$
\begin{align}
    J_F(X^k)\delta=-F(X^k)\label{NL2a}\\
     X^{k+1}=X^k+\delta\,,\label{NL2b}
\end{align}
where the solution of the system \eqref{NL2a} is meant, for arbitrary $M$ and $N$, in the following generalized sense.\\

\begin{proposition} (\cite{GWP})
    The vector
    \begin{equation}\label{NL3}
   \d^\star=-J_F(X^k)^\dag F(X^k)
\end{equation}
is the unique vector of smallest Euclidean norm which minimizes the Euclidean norm of the residual $J_F(X^k)\delta+F(X^k)$.
\end{proposition}\\

\noindent It is easy to see that the generalized solution \eqref{NL3} is given  \\
\begin{itemize}
    \item  for square systems ($M=N$)  by
     $$\delta^\star=-J_F^{-1}(X^k) F(X^k)\,,$$
     provided that the Jacobian is invertible.\\
    \item  for overdetermined systems ($M>N$)  by  the   {\em least-squares} solution
      \begin{equation}\label{NL4}
      \delta^\star=\displaystyle\arg\min_{d\in\R^N}\|J_F(X^k) d+F(X^k)\|_2^2,
      \end{equation}
      provided that the Jacobian has full column rank $N$.\\
    \item  \noindent for  underdetermined systems ($M<N$) by the {\em min Euclidean norm least-squares} solution
      \begin{equation}\label{NL5}
      \delta^\star=\displaystyle\arg\min_{d\in\R^N}\|d\|_2^2 \qquad\mbox{subject to } J_F(X^k)d+F(X^k)=0,
      \end{equation}
      provided that the Jacobian has full row rank $M$.\\
\end{itemize}
In each case, the generalized solution $\d^\star$ can be efficiently obtained avoiding the computation of the Moore-Penrose pseudoinverse.
Indeed, it suffices to perform a $QR$ factorization of the Jacobian
$J_F(X^k)$ in the overdetermined case (or its transpose in the underdetermined case), i.e., a factorization of the form $J_F(X^k)=QR$ in which
\begin{align*}
Q=\left[\begin{array}{cc}Q_1 & Q_2\end{array}\right]\,  &\text{is an $M\times M$ orthogonal matrix with $Q_1$ of size $M\times N$,}\\
& \text{$Q_2$ of size $M\times(M-N)$;}\\
 R=\left[\begin{array}{c}R_1 \\ 0\end{array}\right] \, &\text{is an $M\times N$
matrix with $R_1$ upper triangular of size $N\times N$}\\
& \text{and the null matrix of size $(M-N)\times N$.}
\end{align*}
More precisely:\\
\begin{itemize}
 \item  in the square case, factoring $J_F(X^k)=QR$, we get $J_F(X^k)^\dag=R^{-1}Q^T$ and therefore
$$\delta^\star=-R^{-1}Q^T F(X^k) \, \Longleftrightarrow\, R\delta^\star=-Q^T F(X^k) \,,$$
where the last step is readily computed by back substitution;\\
\item in the overdetermined case, factoring $J_F(X^k)=QR$,  we get by \eqref{NL4}
\begin{align*}
    \delta^\star&=\displaystyle\arg\min_{d\in\R^N}\left\| \left[\begin{array}{c}
                   R_1 \\ 0
                  \end{array}
    \right] d+\left[\begin{array}{c}
                   Q_1^T \\ Q_2^T
                  \end{array}
    \right]F(X^k)\right\|_2^2\\[4pt]
    &= \|Q_2^T F(X^k)\|_2^2+\displaystyle\arg\min_{d\in\R^N}\|R_1 d+Q_1^T F(X^k)\|_2^2
\end{align*}
so that, minimizing the second term,  we get
 \begin{equation}\label{NL6}
 R_1 \delta^\star=-Q_1^T F(X^k)\,,
 \end{equation}
 which is computed again by back substitution.\\
\item in the underdetermined case, factoring $J_F(X^k)^T=QR$ (note that now $M$ and $N$ are exchanged), we have $J_F(X^k)=R^T Q^T=R_1^T Q_1^T$.
Moreover, setting $$d=Qz=Q_1 z_1+Q_2 z_2\,,$$ we get\vspace{-3pt}
\begin{align*}
0&=J_F(X^k)d+F(X^k)=R_1^TQ_1^T(Q_1 z_1+Q_2z_2)+F(X^k)\\&=R_1^Tz_1+R_1^TQ_1^TQ_2z_2+F(X^k)\,.
\end{align*}
Since $Q_1^TQ_2=0$ by the orthogonality of $Q$, we obtain the constraint
\begin{equation}\label{aux}
R_1z_1+F(X^k)=0\,.
\end{equation}
It follows that we can minimize (see \eqref{NL5})
$$\|d\|_2^2=\|Q^Td\|_2^2=\|z_1\|_2^2+\|z_2\|_2^2\,,$$
just taking $z_2=0$, and we conclude that
$$\delta^\star=Q_1z_1=-Q_1R_1^{-T}F(X^k)\,,$$
where $z_1$ is computed by \eqref{aux} again via back substitution.\\
\end{itemize}

From now on we will refer to the generalized solution \eqref{NL3} for arbitrary $M$ and $N$  as to the {\em least-squares} solution.\\

\noindent Summarizing,  we consider the following   algorithm for the solution of \eqref{s:3}: \\\\
\textsc{
Given an initial guess $X$ and a tolerance $\varepsilon>0$,\\
repeat
\begin{itemize}
 \item[$(1)$] Assemble $F(X)$ and $J_F(X)$
 \item[$(2)$] Solve the linear system $J_F(X)\delta=-F(X)$ in the least-squares sense, using the $QR$ factorization of $J_F(X)$ (or of $J_F(X)^T$)
 \item[$(3)$] Update $X\,\leftarrow\,X+\delta$
\end{itemize}
until $\|\delta\|_2^2<\varepsilon$ and/or $\|F(X)\|_2^2<\varepsilon$}\\\par
In the actual code implementation of the algorithm above,
we employ several well known variants and modifications of the classical Newton method, as discussed in the following remarks.\\\\
$\bullet$
The convergence of Newton-like methods is in general local. Nevertheless in some cases,
 as in  \eqref{HJ}, if $H$ is convex and a proper discretization preserving this property is performed, the map $F$ in \eqref{s:3} is also convex
and therefore the convergence is global. Moreover, in every example we consider (except for multi-population mean field games, see Section \ref{multipopMFG}) the ergodic constant is unique.\\\\
$\bullet$ Sometimes Newton-like  methods do not converge, due to oscillations around a minimum of the residual function $\|F(X)\|_2^2$. This situation can be successfully
 overcome by introducing a {\em dumping parameter} in the update step, i.e., by replacing $X$ in the step $(3)$ of the algorithm with $X+\mu\delta$ for some $0<\mu<1$.
 Usually a fixed value of $\mu$ works fine, possibly affecting the number of iterations to obtain convergence.
 A more efficient (but costly) selection of the dumping parameter can be implemented using {\em line search} methods, such as the {\em inexact backtracking} line search with
 {\em Armijo-Goldstein condition}. Especially when dealing with {\em nonconvex} residual functions, the Newton step can be trapped into a local minimum.
 In this case a re-initialization of the dumping parameter to a bigger value can resolve the situation, acting as a thermalization in {\em simulated annealing} methods. \\\\
$\bullet$ It may happen (usually if the initial guess is $X=0$) that $J_F(X)$ is nearly singular or rank deficient, so that the least-squares solution cannot be
 computed. In this case, in the spirit of the {\em Levenberg-Marquardt} method or more in general quasi-Newton methods, we can add a regularization term
 on the principal diagonal of the Jacobian, by replacing $J_F(X)$ with $\tau I+J_F(X)$, where $\tau>0$ is a tunable parameter and $I$ denotes the identity (not necessarily square) matrix.
 This correction does not affect the solution, but it may slow down the convergence if $\tau$ is not chosen properly.
 This depends on the fact that the method reduces to a simple gradient descent if the term
 $\tau I$ dominates $J_F(X)$. In our implementation, we switch-on the correction (with a fixed $\tau$) only at the points $X$ where $J_F(X)$ is nearly singular or
 rank deficient. In this way we can easily handle, for instance, second order problems with very small diffusion coefficients.\\\\
$\bullet$ Newton-like methods classically require the residual function to be Fr\'echet differentiable. Nevertheless, this assumption can be weakened to include
 important cases, such as cell problems in which the Hamiltonian is of the form $H(x,p)=\frac{1}{q}|p|^q-V(x)$ with $q\ge 1$.
 Note that the derivative in $p$ is given by $H_p(x,p)=|p|^{q-2}p$, so that the Jacobian in the corresponding Newton step is not differentiable at the
 origin for $1\le q<2$.
 In this situation, in the spirit of nonsmooth-Newton methods, we can replace the usual gradient with any element of the sub-gradient. Typically we choose
 $H_p(x,p)=0$ for $p=0$.\\\\
 $\bullet$
It is interesting to observe that in the overdetermined case, the iterative method \eqref{NL2a}-\eqref{NL2b} coincides with
the {\em Gauss-Newton} method for the optimization problem
$\min_{X}\frac 1 2\|  F(X)\|_2^2$.
Indeed, defining $\cF(X)=\frac12\| F(X)\|_2^2$, the classical  Newton method for the critical points  of $\cF$ is given by
$$
H_\cF(X^k)(X^{k+1}-X^k)=-\nabla \cF(X^k) \qquad k\ge 0\,.
$$
Computing the  gradient $\nabla \cF$ and the Hessian $H_\cF$ of $\cF$ we have
$$\nabla \cF(X)=J_F(X)^T   F(X)\,,\qquad H_\cF(X)=J_F(X)^T J_F(X)+\sum_{i=1}^{m}\frac{\partial^2 F_i}{\partial X^2}(X)
  F_i(X)\,,$$
where the second order term is given by
$$
\left(\frac{\partial^2  F_i}{\partial^2 X}(X)\right)_{k,\ell}=\frac{\partial^2 F_i}{\partial X_k \partial X_\ell}(X)\,.$$
Since  the minimum of $\cF(X)$  is zero, we expect that  $F(X^k)$ is small for $X^k$ close enough to a solution $X^\star$.
Hence we  approximate $\cH_\cF(X^k)\simeq J_F(X^k)^T J_F(X^k)$, obtaining
the Gauss-Newton method, involving only first order terms:
$$
J_F(X^k)^T J_F(X^k)\left(X^{k+1}-X^k\right)=-J_F(X^k)^T   F(X^k)\qquad k\ge 0\,.
$$
Applying again $QR$ decomposition to $J_F(X^k)$ we finally get the iterative method \eqref{NL2a}-\eqref{NL2b}, with $\d$ given
by the least-squares  solution \eqref{NL6}.
This is the approach we followed for solving stationary mean field games on networks in \cite{CCM16}.\\

Throughout the next sections we present several ergodic problems, that can be set in our framework and solved by
the proposed Newton-like method for inconsistent systems.
Each section contains numerical tests in dimension one and/or two, including some experimental convergence analysis
and also showing the performance of the proposed method, both in terms of accuracy and computational time.
All tests were performed on a Lenovo Ultrabook X1 Carbon, using 1 CPU Intel Quad-Core i5-4300U 1.90Ghz with 8 Gb Ram,
running under the Linux Slackware 14.1 operating system. The algorithm is implemented in \textsc{C} and  employs the library {\em SuiteSparseQR}
\cite{SPQR}, which is designed to efficiently compute in parallel the $QR$ factorization and the least-squares solution to very large and sparse linear systems.

\section{The eikonal Hamiltonian}\label{Eikonal}
We start considering the simple case of the eikonal equation in dimension one, namely the cell problem
$$
\frac12|u^\prime+p|^2-V(x)=\lambda\qquad \mbox{in }\mathbb{T}^1\,,
$$
where $p\in\R$, $\l\in\R$ and $V$ is a $1\hbox{-}$periodic potential. This is a good benchmark for the proposed method, since
a formula for the effective Hamiltonian is available (see \cite{LPV}):
$$
\bar{H}(p)=\left\{
\begin{array}{ll}
 -\min V & \mbox{if } |p|\le p_c\\
 \lambda & \mbox{if } |p|> p_c\quad\mbox{s.t. }|p|=\int_0^1\sqrt{2(V(s)+\lambda)}ds
\end{array}
\right.
$$
where $p_c=\int_0^1\sqrt{2(V(s)-\min V))}ds$. Note that the effective Hamiltonian $\bar{H}$ has a plateau in the whole interval
$[-p_c,p_c]$. With a slight abuse of notation, in what follows we will refer to this interval as to {\em the} plateau.

\begin{figure}[h!]
 \begin{center}
 \begin{tabular}{cc}
\includegraphics[width=.32\textwidth]{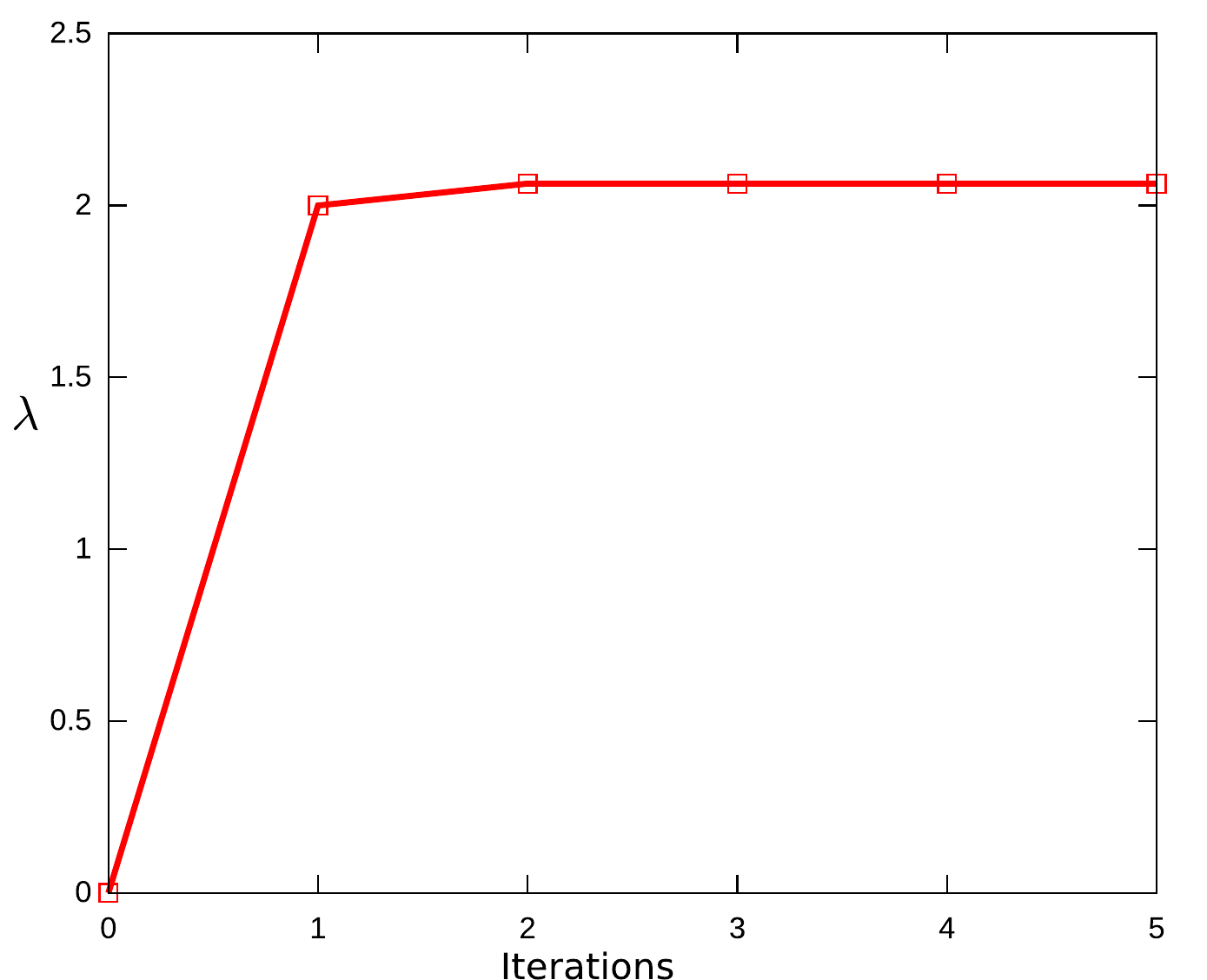} &
\includegraphics[width=.32\textwidth]{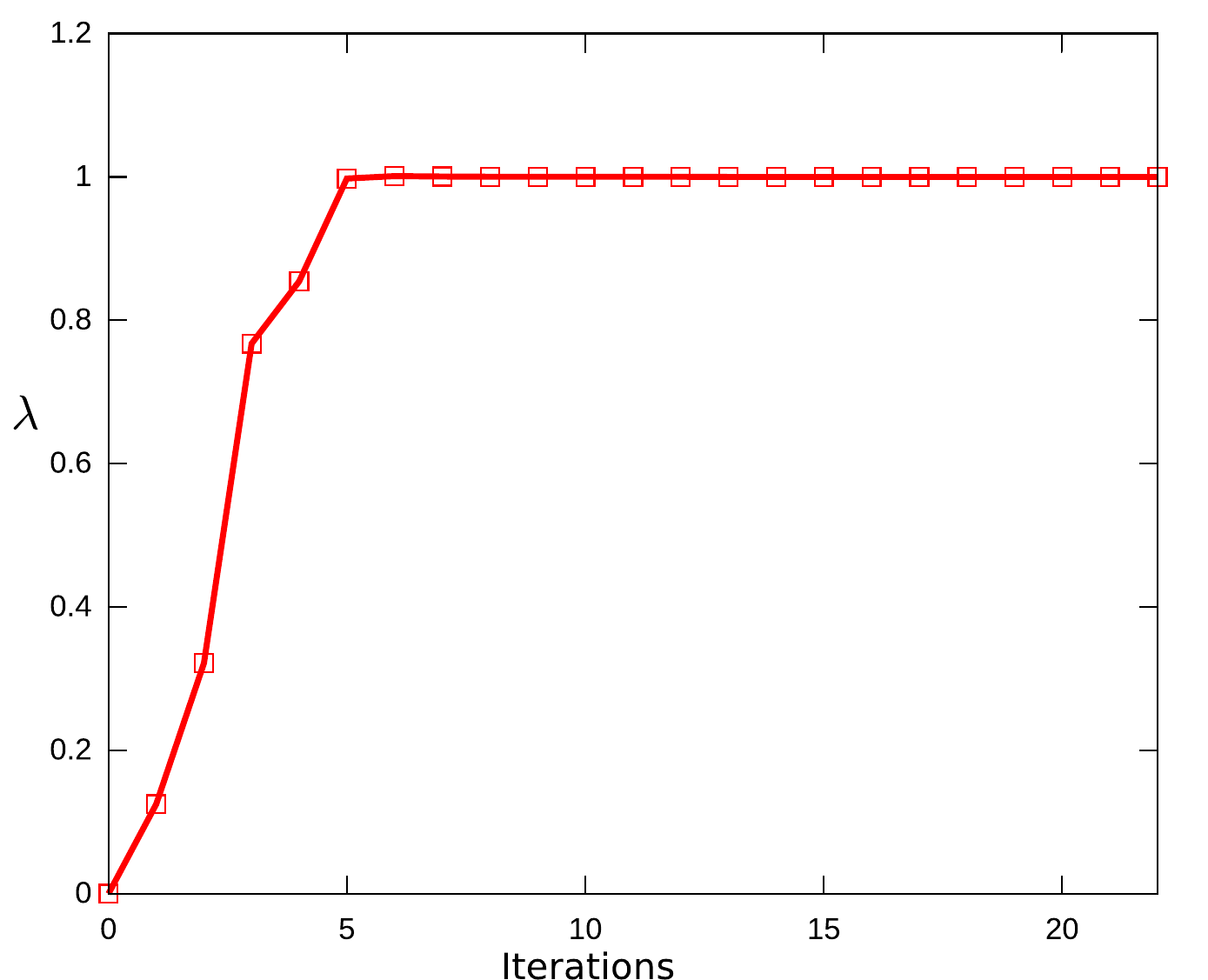} \\
(a)&(b)
\end{tabular}
\end{center}
\caption{Convergence of the method: $lambda$ vs number of iterations for $p=2$ (a) and $p=0.5$ (b).}\label{lambda-vs-iterations}
\end{figure}

\begin{figure}[h!]
 \begin{center}
 \begin{tabular}{cc}
\includegraphics[width=.3\textwidth]{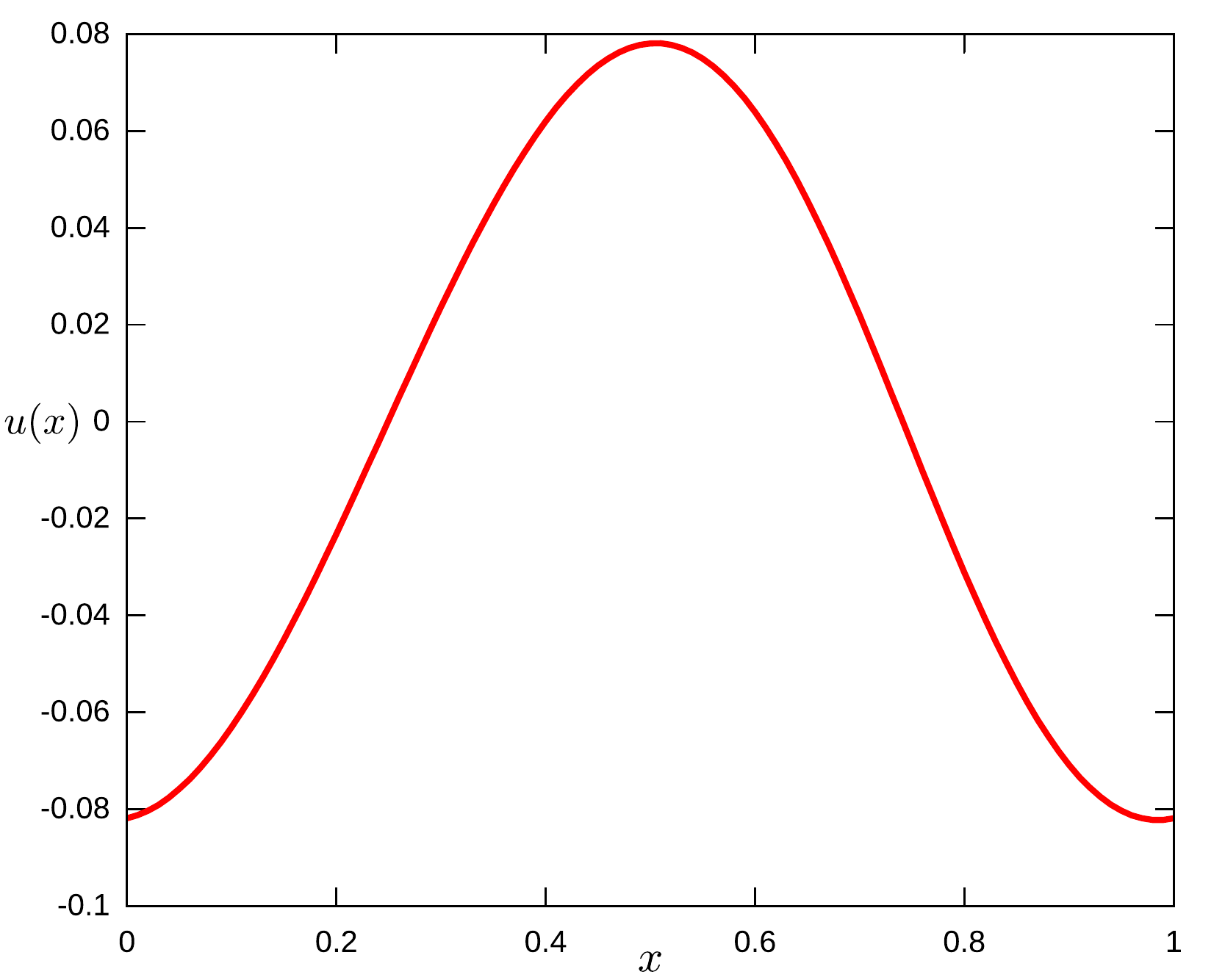} &
\includegraphics[width=.3\textwidth]{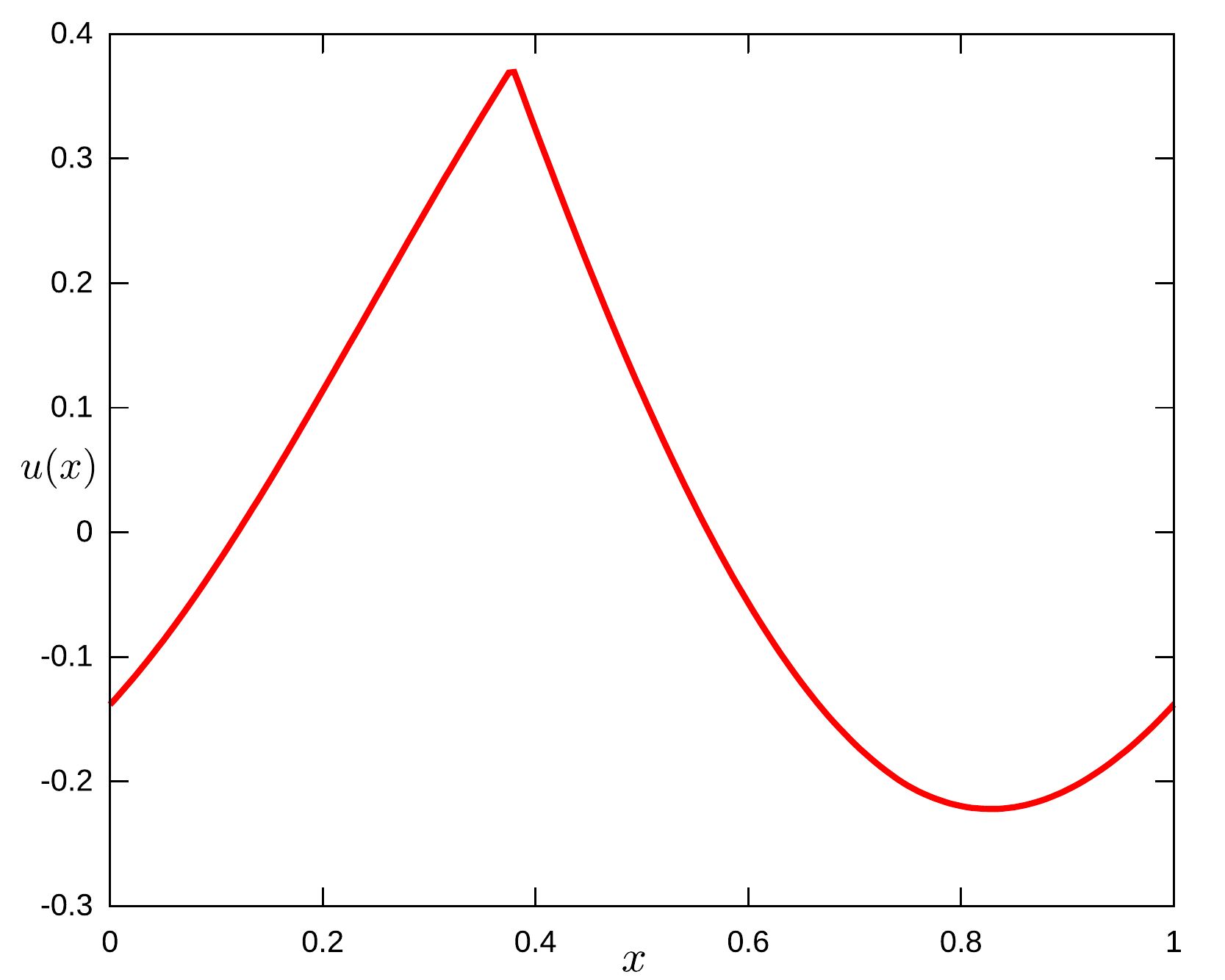} \\
(a)&(b)
\end{tabular}
\end{center}
\caption{Viscosity solution of the cell problem for $p=2$ (a) and $p=0.5$ (b).}\label{correctors-eikonal-1d}
\end{figure}

\begin{figure}[h!]
 \begin{center}
 \begin{tabular}{cc}
\includegraphics[width=.35\textwidth]{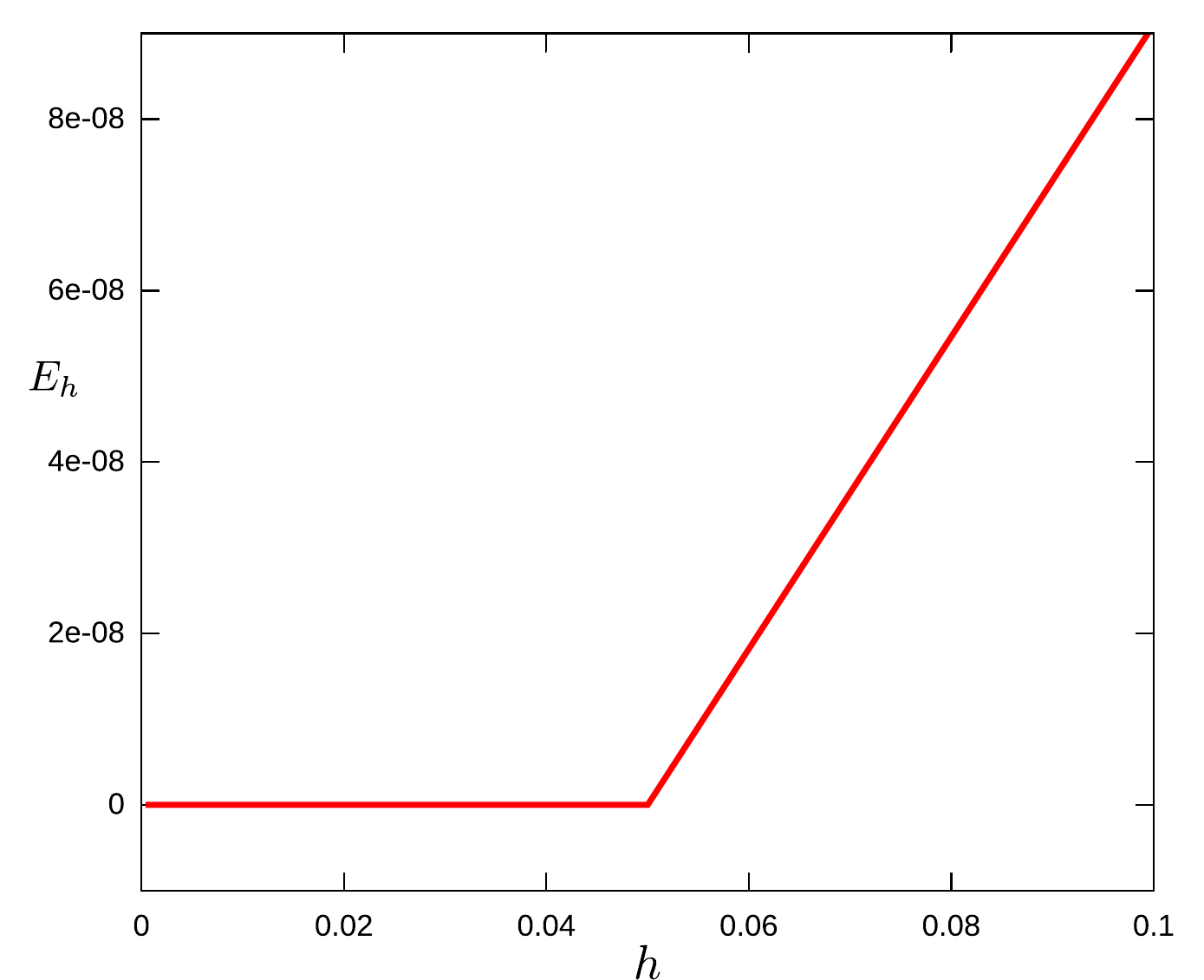} &
\includegraphics[width=.35\textwidth]{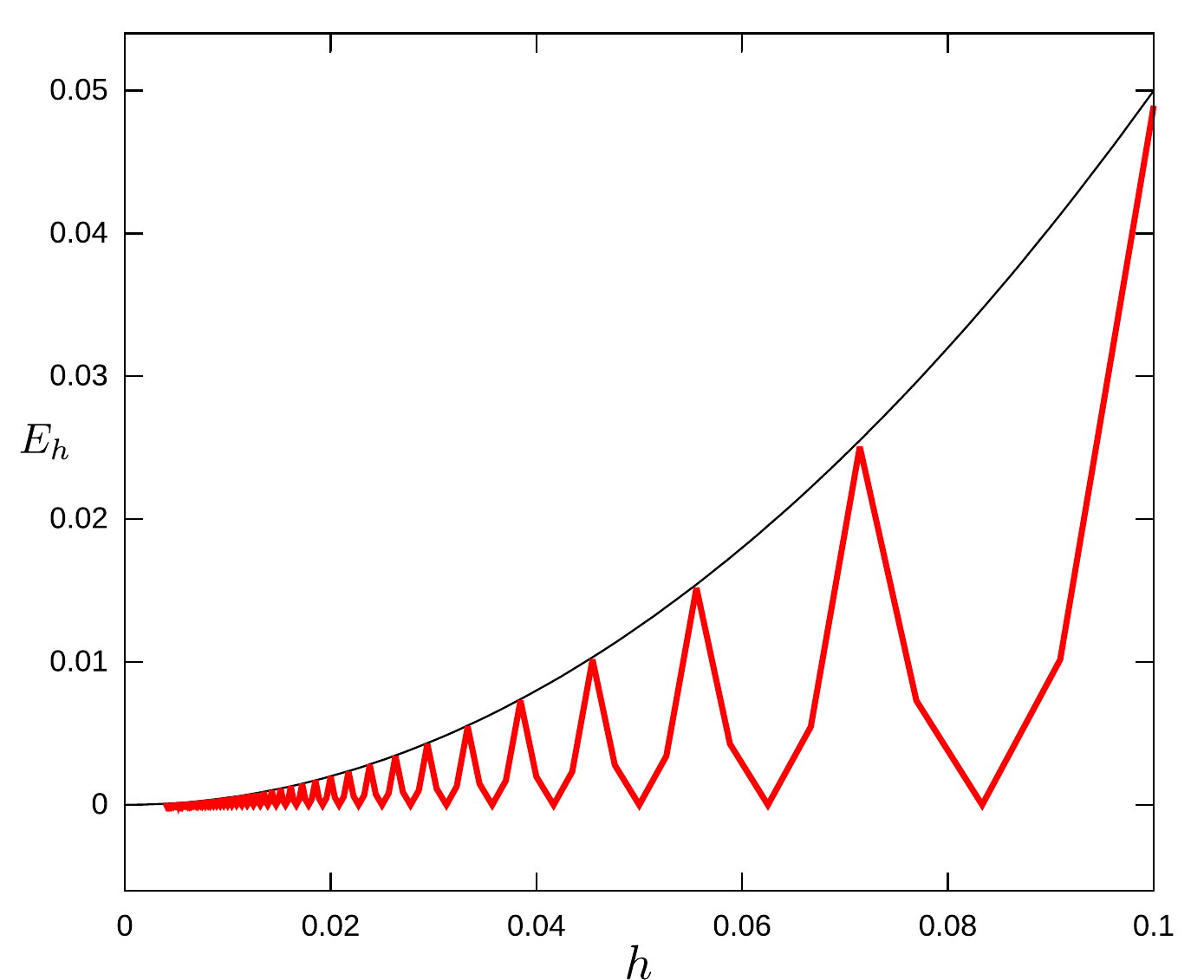} \\
(a)&(b)
\end{tabular}
\end{center}
\caption{Convergence under grid refinement for $p=2$ (a) and $p=0.5$ (b).}\label{convergence-eikonal-1d}
\end{figure}

Following \cite{Q}, in our first test we choose $V(x)=\sin(2\pi x)$, for which $\min V=-1$ and $p_c=4/\pi$.
As initial guess we always choose $(U,\Lambda)=(0,0)$ and we set to $\varepsilon=10^{-6}$ the tolerance for the stopping criterion of the algorithm.
Figure \ref{lambda-vs-iterations} shows the computed $\lambda$ as a function of the number of iterations to reach convergence.
\\
In particular, Figure \ref{lambda-vs-iterations}a
corresponds to $p=2$, which is outside the plateau of $\bar H$, whereas Figure \ref{lambda-vs-iterations}b corresponds to $p=0.5\in[-p_c,p_c]$.
In both cases the torus $\mathbb{T}^1$ is discretized with $100$ nodes and the computation occurs in {\em real time}.
Note that the convergence is very fast, much more in the first case. We observe that this depends on the fact that the corresponding {\em corrector} (i.e., the viscosity
solution $u$ of the cell problem) is smooth in the first case and only Lipschitz in the second, so that the algorithm needs some further iteration to 
compute the kinks, as shown in Figure \ref{correctors-eikonal-1d}.
\\
In Figure \ref{convergence-eikonal-1d} we plot the error in the approximation of $\lambda$ under grid refinement, for both $p=2$ and $p=0.5$.
The behavior of the error is quite surprising. For $p$ outside the plateau of $\bar H$ and $h$ sufficiently small, the error seems to be independent of the grid size,
close to machine precision. On the other hand, for $p$ in the plateau, the error is at worst quadratic in the space step
(this is typical for first order Newton-like methods),
but with some strange oscillations, corresponding to particular choices of the grid, for which the error is close to machine precision. The reason
of this phenomenon of super-convergence is currently unclear to us.

\begin{figure}[h!]
 \begin{center}
 \begin{tabular}{cc}
\includegraphics[width=.35\textwidth]{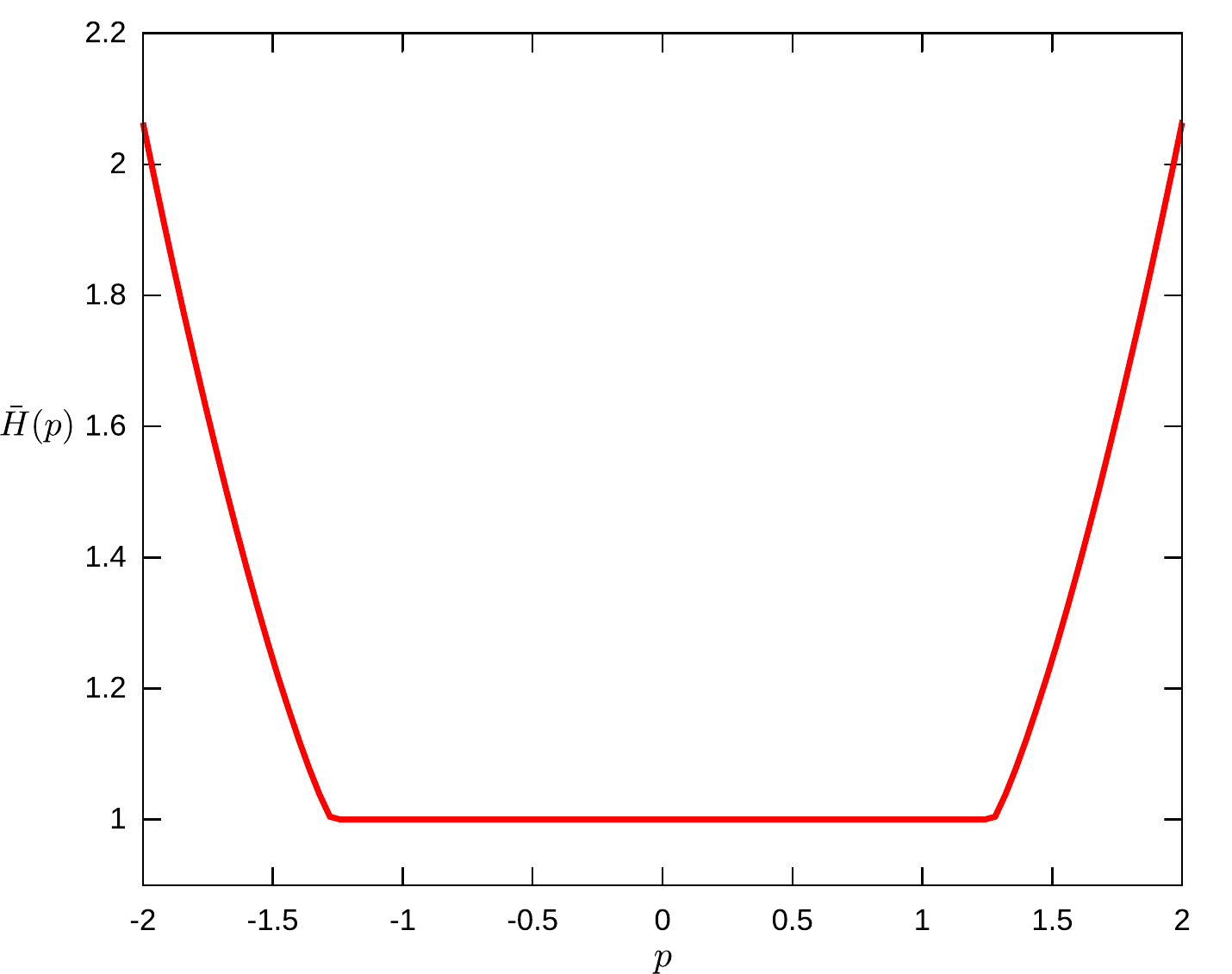} &
\includegraphics[width=.35\textwidth]{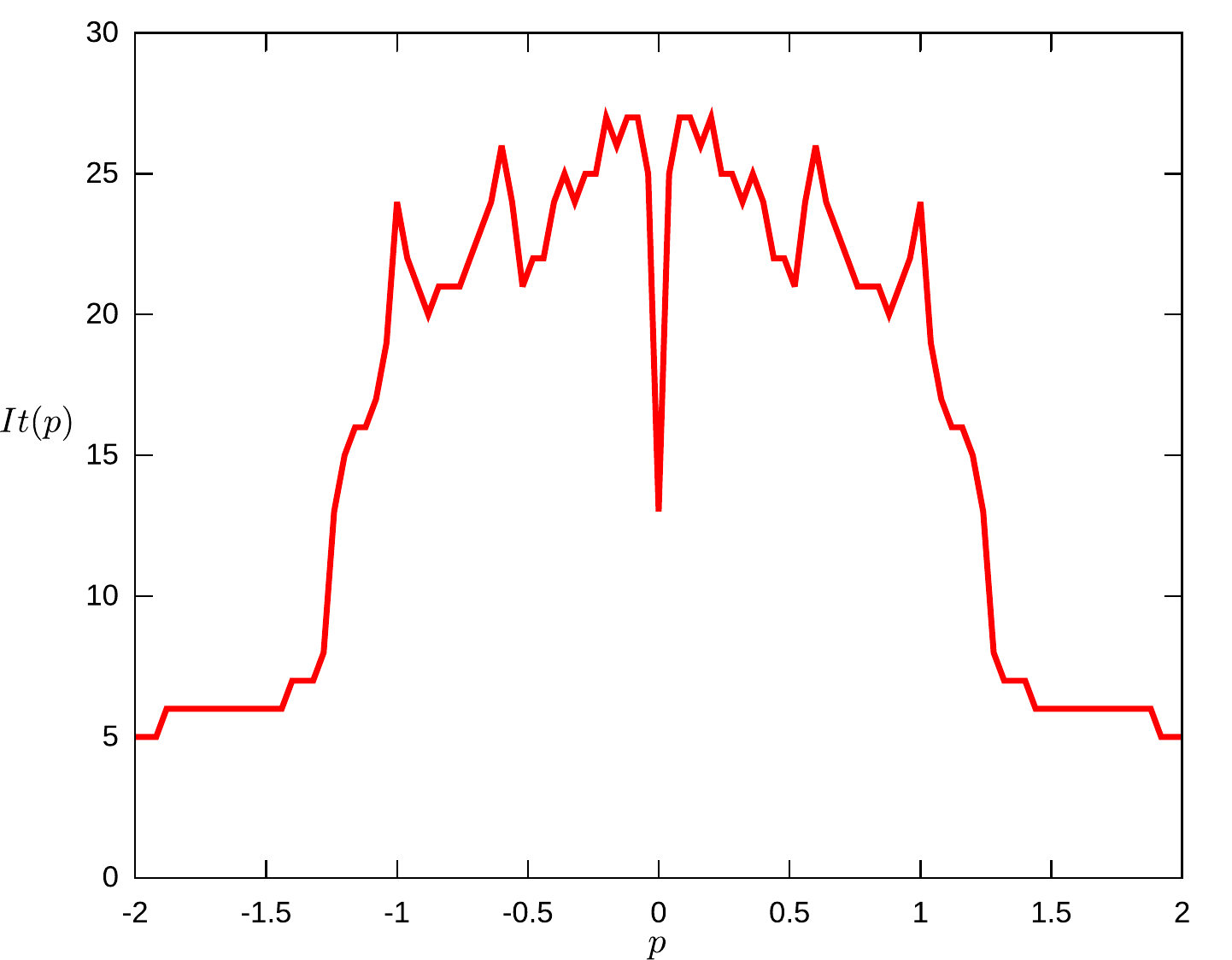} \\
(a)&(b)
\end{tabular}
\end{center}
\caption{Effective Hamiltonian for $p\in[-2,2]$ (a) and the corresponding iterations (b).}\label{Hbar-eikonal-1d}
\end{figure}

In the next test we compute the effective Hamiltonian for $p$ in the interval $[-2,2]$, discretized with $101$ uniformly distributed nodes. The mesh size for
the torus $\mathbb{T}^1$ is again $100$ nodes. In Figure \ref{Hbar-eikonal-1d}a we show the graph of $\bar H$ and in Figure \ref{Hbar-eikonal-1d}b
the corresponding number of iterations to reach convergence as a function of $p$.
\\
We remark again how the number of iterations increases in the plateau $[-p_c,p_c]$.
The total computational time for this simulation is just $1.18$ seconds. \\

Now let us consider the same problem in dimension two, namely
$$
\frac12|Du+p|^2-V(x_1,x_2)=\lambda\qquad \mbox{in }\mathbb{T}^2\,,
$$
where $p\in\R^2$, $\l\in\R$ and $V$ is again a  $1\hbox{-}$periodic potential.
As a first example we choose
\begin{equation}\label{ex2d-1}
V(x_1,x_2)=\cos(2\pi x_1)+\cos(2\pi x_2)\,,
\end{equation}
corresponding to the celebrated problem of two uncoupled penduli.
In this case the effective Hamiltonian is separable, i.e., it is just the sum of the two effective Hamiltonians associated to the one dimensional potential
$V(x)=cos(2\pi x)$:
\begin{equation}\label{hbarsep}
\bar H(p)=\bar H_{1d}(p_1)+\bar H_{1d}(p_2)\,,\qquad p=(p_1,p_2)\,.
\end{equation}

\begin{figure}[h!]
 \begin{center}
 \begin{tabular}{cc}
\includegraphics[width=.35\textwidth]{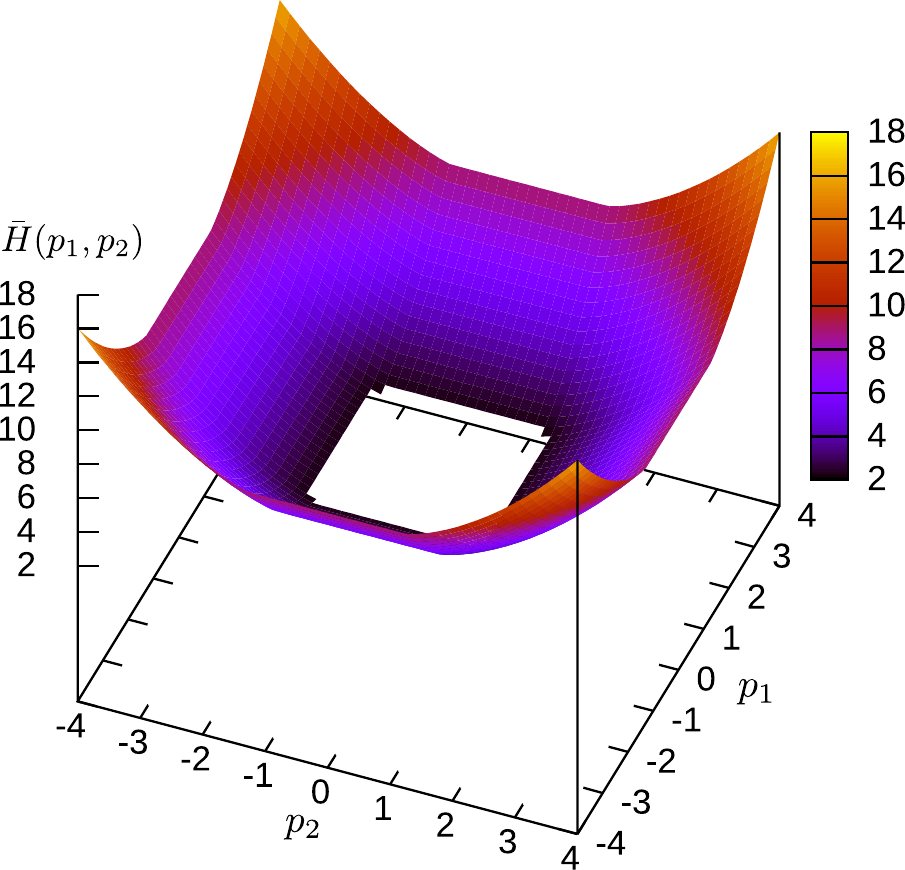} &
\includegraphics[width=.3\textwidth]{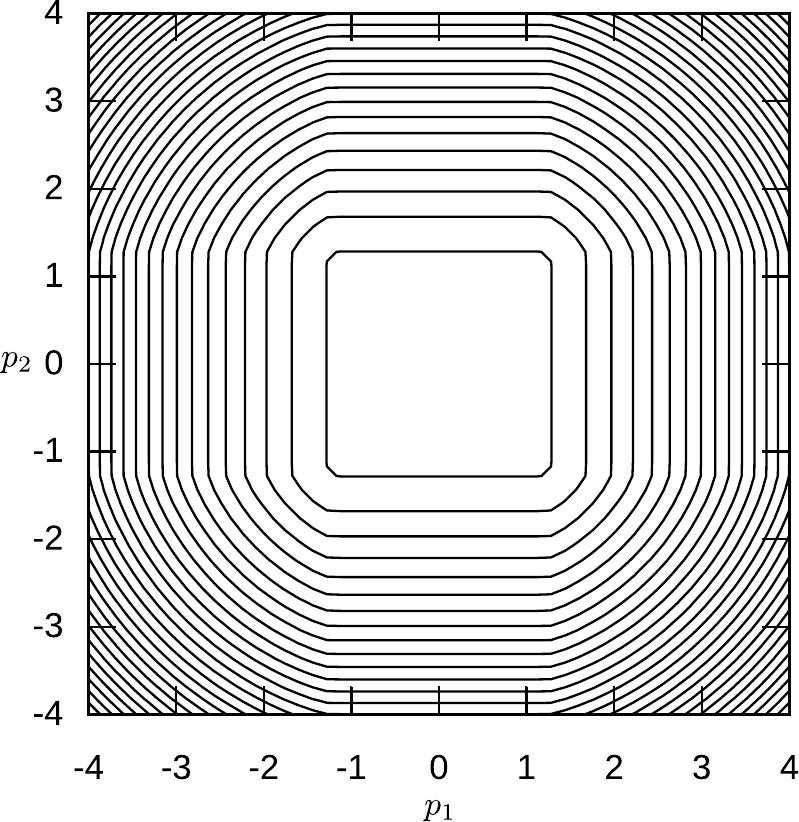} \\
(a)&(b)
\end{tabular}
\end{center}
\caption{Effective Hamiltonian for $p\in[-4,4]^2$: (a) surface and (b) level sets.}\label{Hbar-2d-Q52}
\end{figure}

\noindent Nevertheless, we perform a full 2D computation of $\bar H(p)$ for $p\in[-4,4]^2$ discretized with $51\times 51$ uniformly distributed nodes. The mesh size for
the torus $\mathbb{T}^2$ is $25\times 25$ nodes. Figure \ref{Hbar-2d-Q52} shows the effective Hamiltonian surface and its level sets.
We removed the color filling the level set $2$ to better
appreciate the plateau $\{p\in\R^2\,:\,\bar H(p)=2\}$.\\
For a single point $p$, the average number of iterations is about $16$ and the average computational time is $0.4$ seconds,
whereas the total time is $970.45$ seconds. \\

\begin{figure}[h!]
 \begin{center}
 \begin{tabular}{ccc}
\includegraphics[width=.3\textwidth]{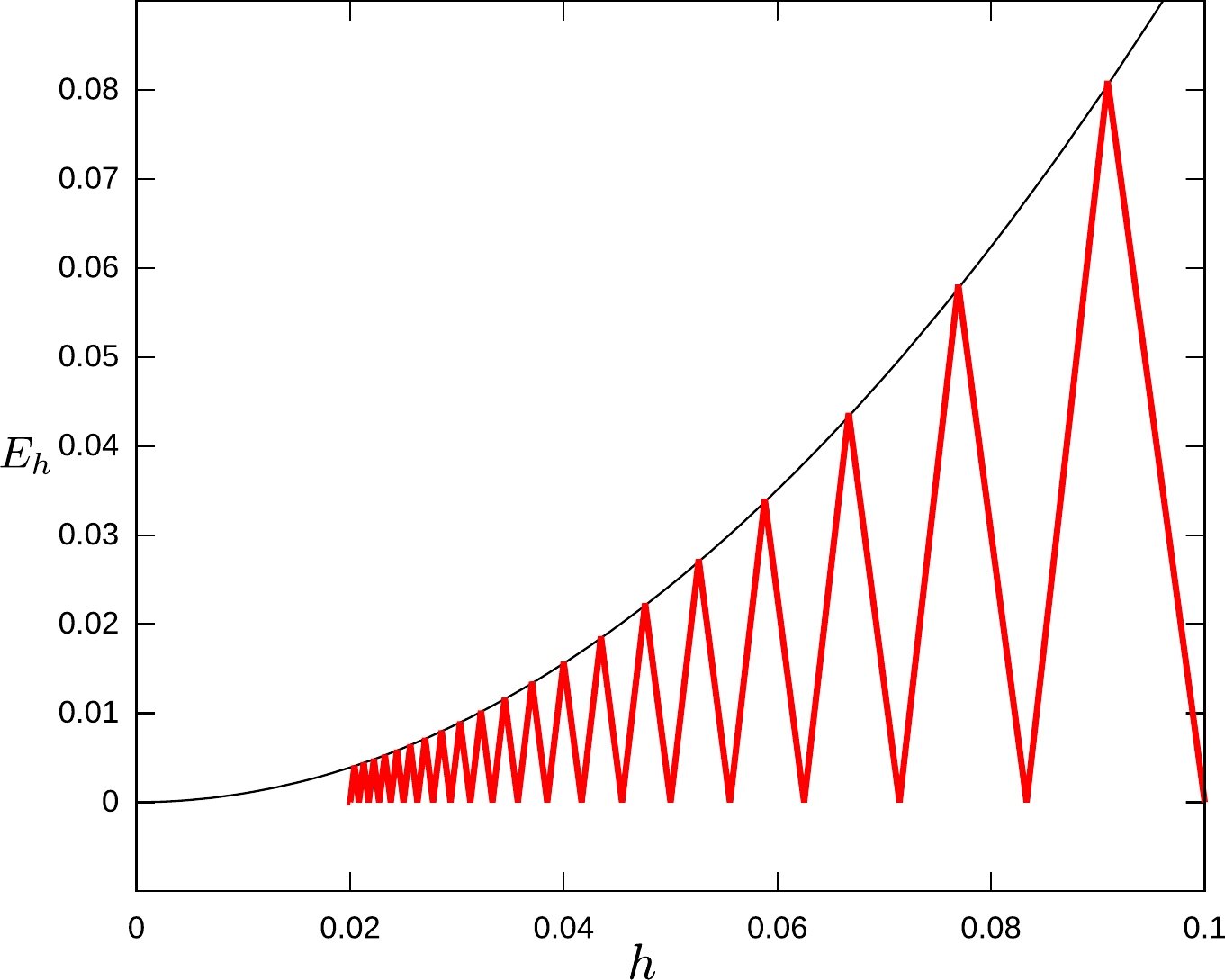} &
\includegraphics[width=.3\textwidth]{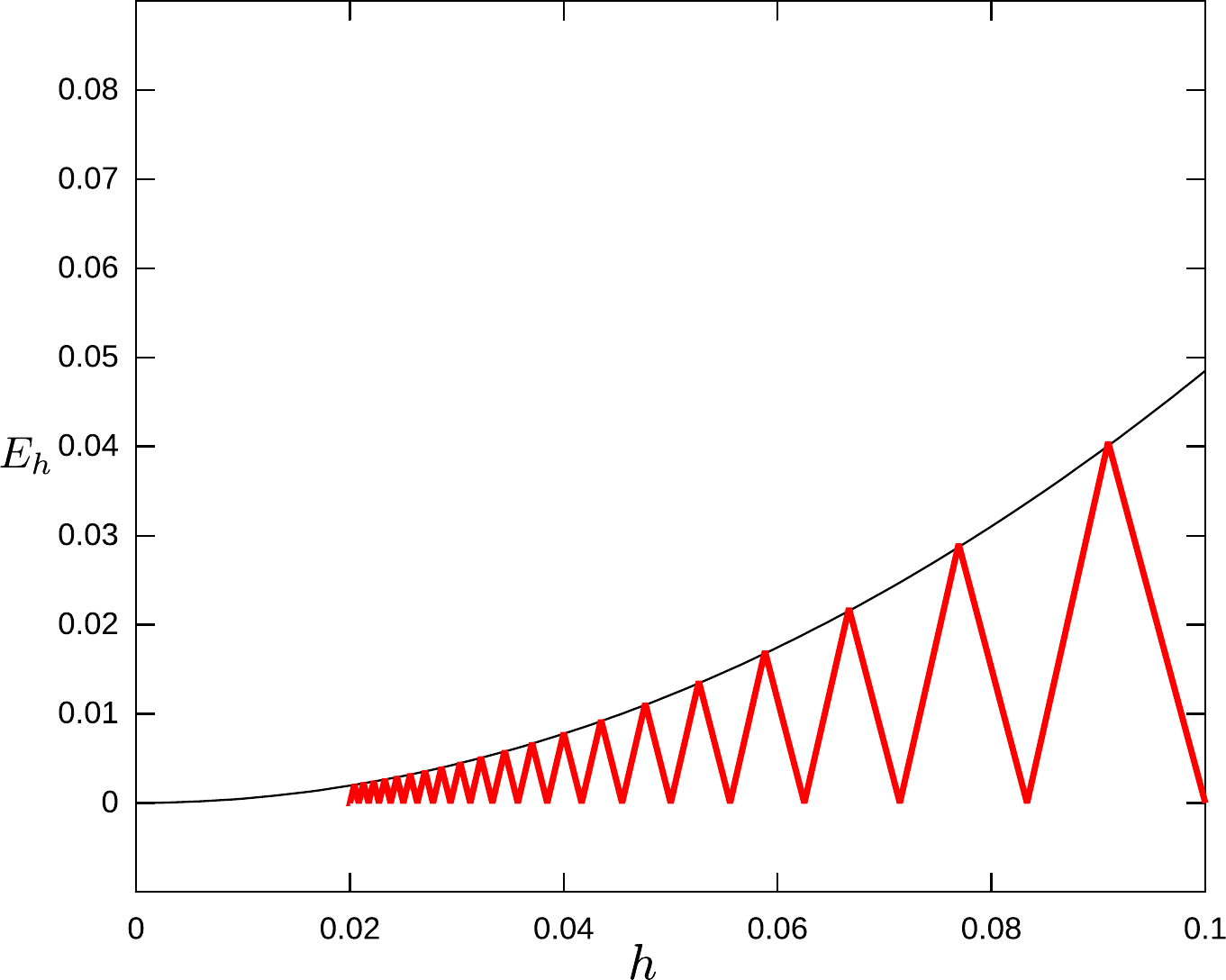} &
\includegraphics[width=.285\textwidth]{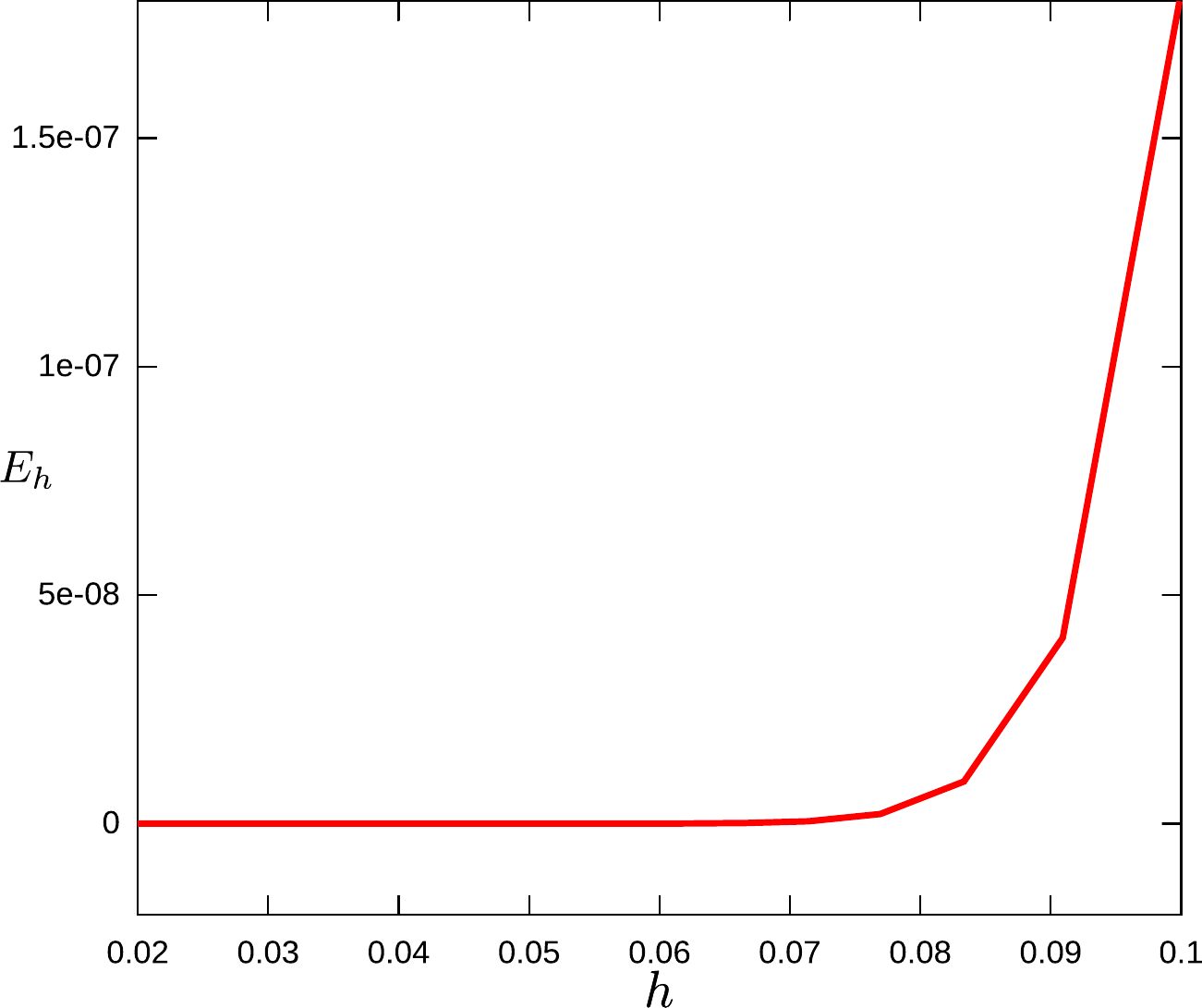} \\
(a)&(b)&(c)
\end{tabular}
\end{center}
\caption{Convergence under grid refinement for $p=(0,0)$ (a), $p=(2,0)$ (b), $p=(2,2)$ (c).}\label{convergence-2d-Q52}
\end{figure}

We now perform, as before, a convergence analysis of the algorithm under grid refinement at fixed $p$. We choose the three points
$p^A=(0,0)$, $p^B=(2,0)$ and $p^C=(2,2)$. In the first case both components $p^A_1$ and $p^A_2$ of $p^A$ are in the plateau of the corresponding
one dimensional effective Hamiltonian. In the second case we have the first component $p^B_1$ outside and the second component $p^B_2$ inside the plateau.
In the third case both components of $p^C$ are outside.
The error is obtained using as correct value for $\bar H$, the formula \eqref{hbarsep} computed by the 1D code for $h=0.0002$.
Figure \ref{convergence-2d-Q52} shows results similar to the one dimensional case.
Indeed, we observe again an experimental convergence at least of order $2$, and even higher outside the plateau.

\begin{figure}[h!]
 \begin{center}
 \begin{tabular}{ccc}
\includegraphics[width=.3\textwidth]{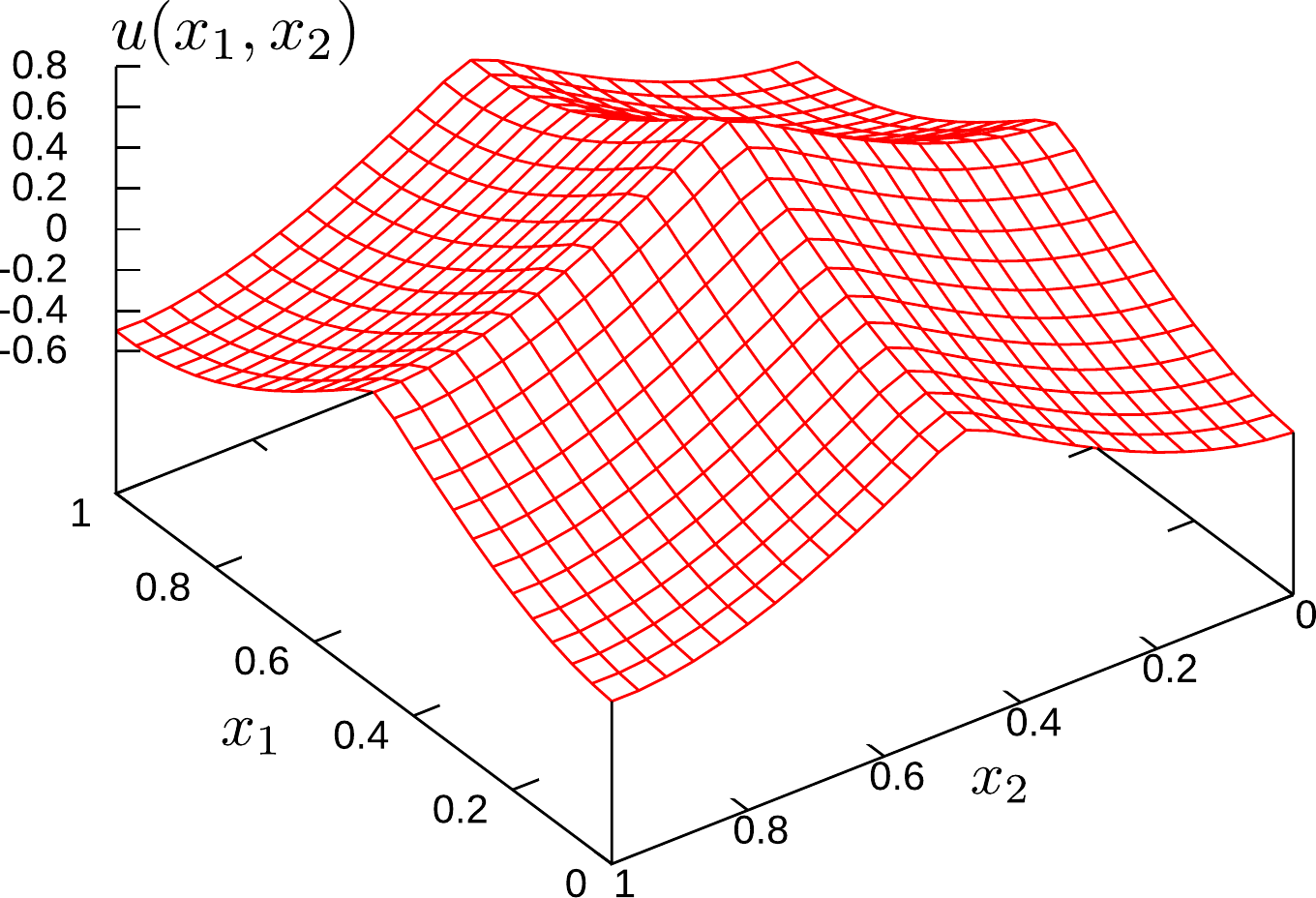} &
\includegraphics[width=.3\textwidth]{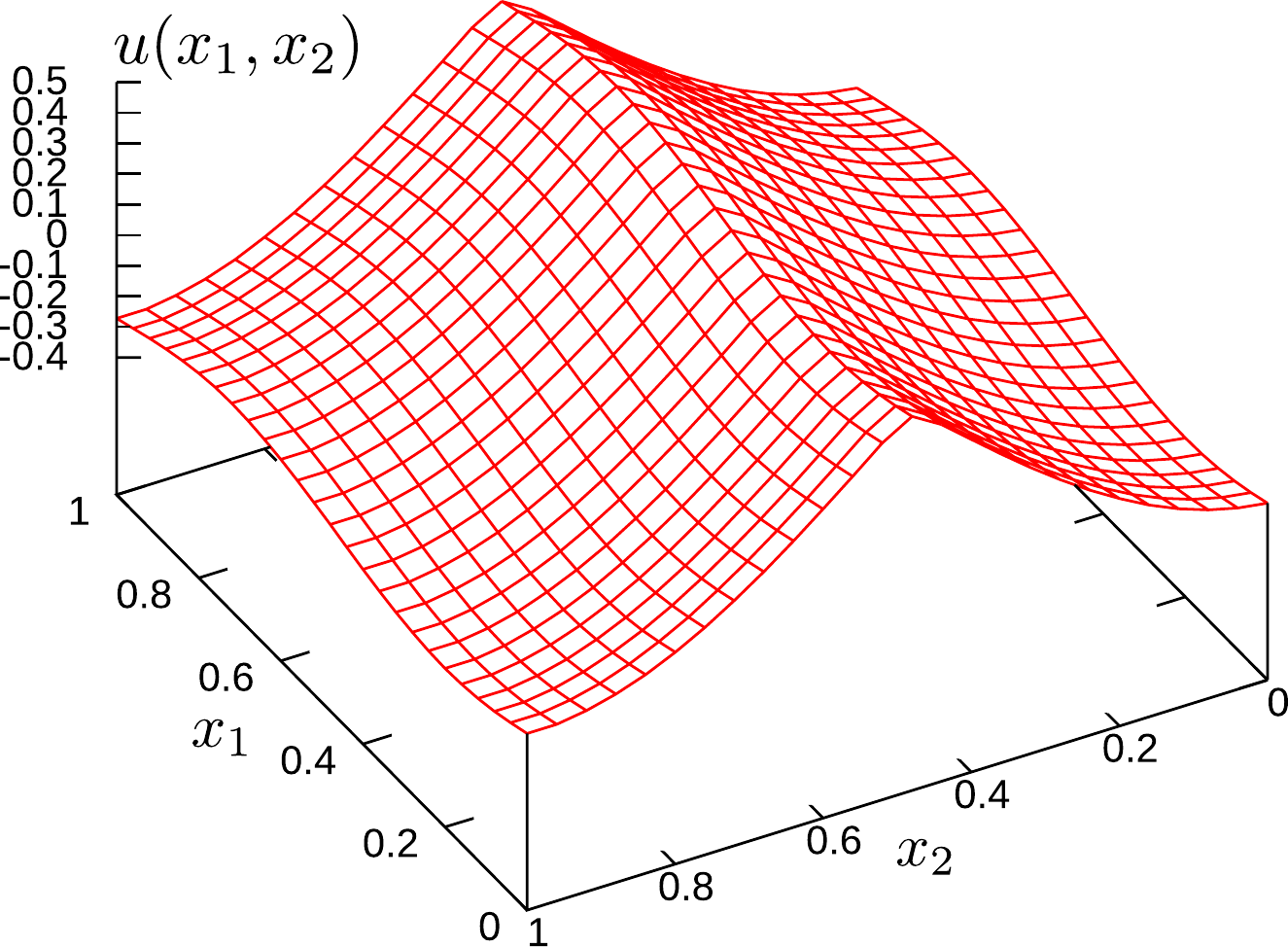} &
\includegraphics[width=.3\textwidth]{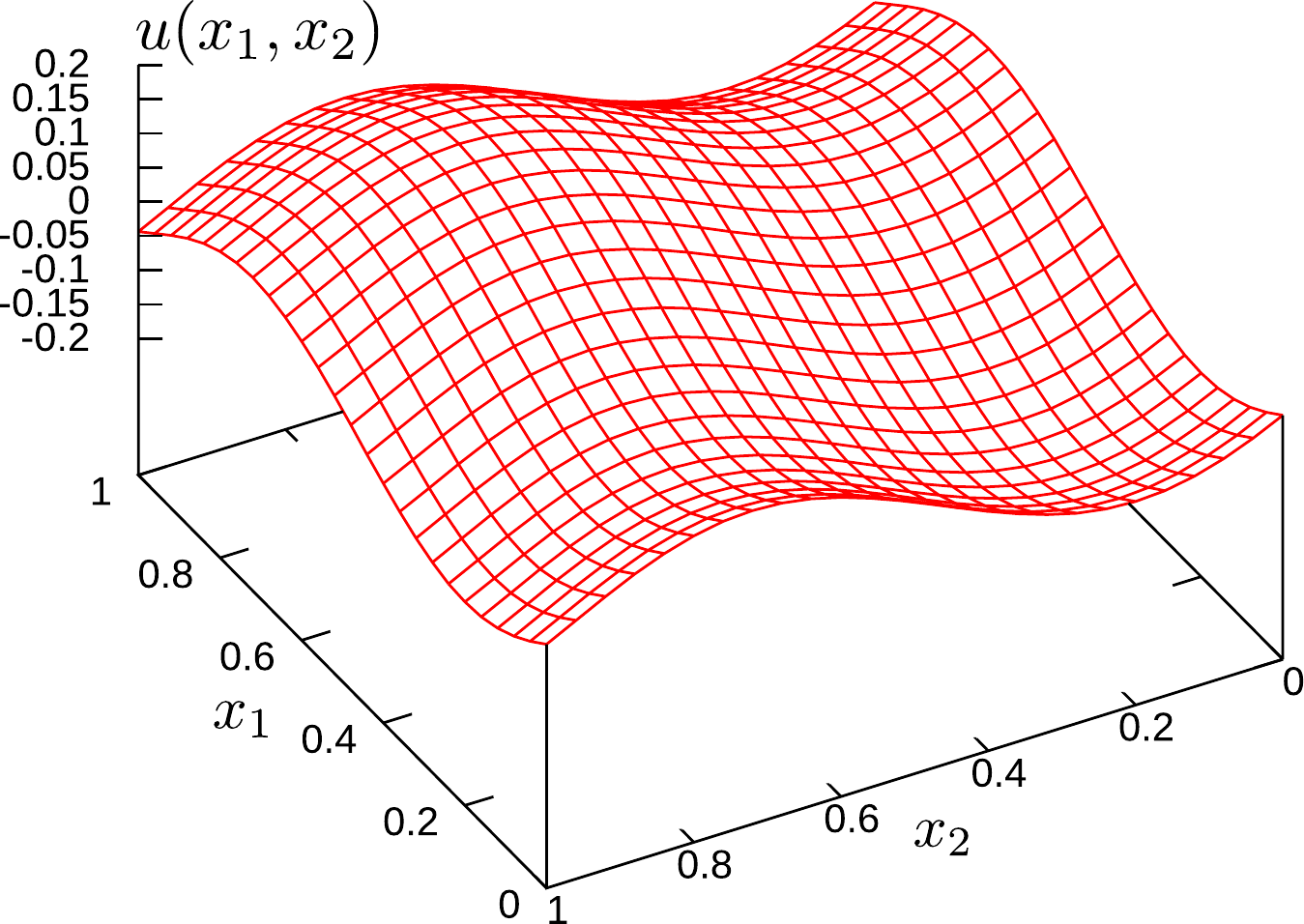} \\
(a)&(b)&(c)
\end{tabular}
\end{center}
\caption{Viscosity solution of the cell problem for $p=(0,0)$ (a), $p=(2,0)$ (b), $p=(2,2)$ (c).}\label{correctors-2d-Q52}
\end{figure}

In Figure \ref{correctors-2d-Q52} we show the correctors corresponding to $p^A$, $p^B$ and $p^C$. It is interesting to note that
the smoothness of the viscosity solution in the directions $x_1$ and $x_2$ precisely depends on what component of $p$ belongs or not to the plateau
of $\bar H$. These results are qualitatively in agreement with those obtained in \cite{R06}.

\begin{figure}[h!]
 \begin{center}
 \begin{tabular}{cc}
\includegraphics[width=.3\textwidth]{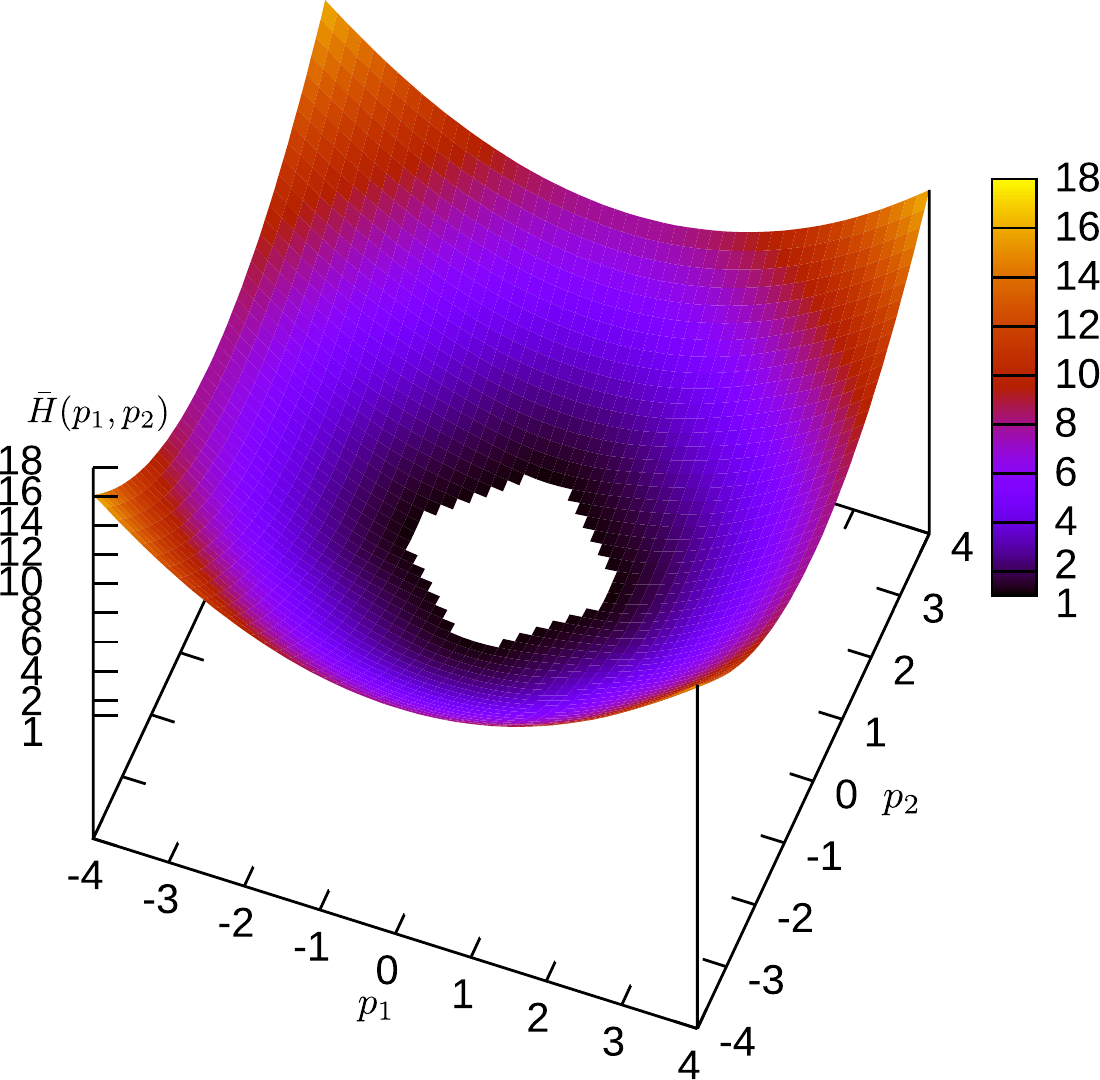} &
\includegraphics[width=.25\textwidth]{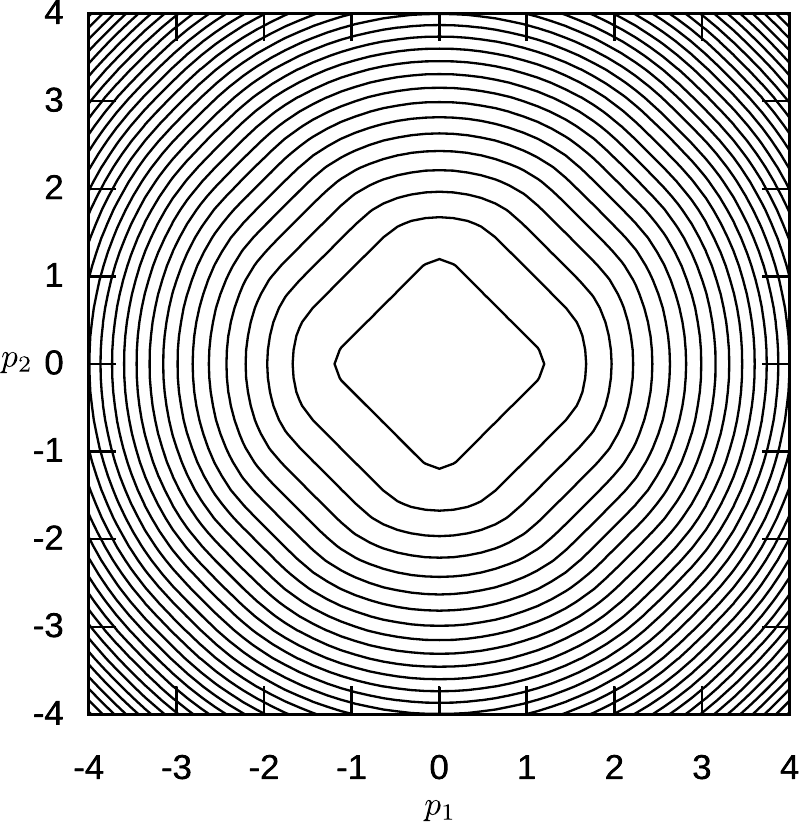} \\
(a)&(b)
\end{tabular}
\end{center}
\caption{Effective Hamiltonian for $p\in[-4,4]^2$: (a) surface and (b) level sets.}\label{Hbar-2d-Q53}
\end{figure}

\begin{figure}[h]
 \begin{center}
 \begin{tabular}{cc}
\includegraphics[width=.3\textwidth]{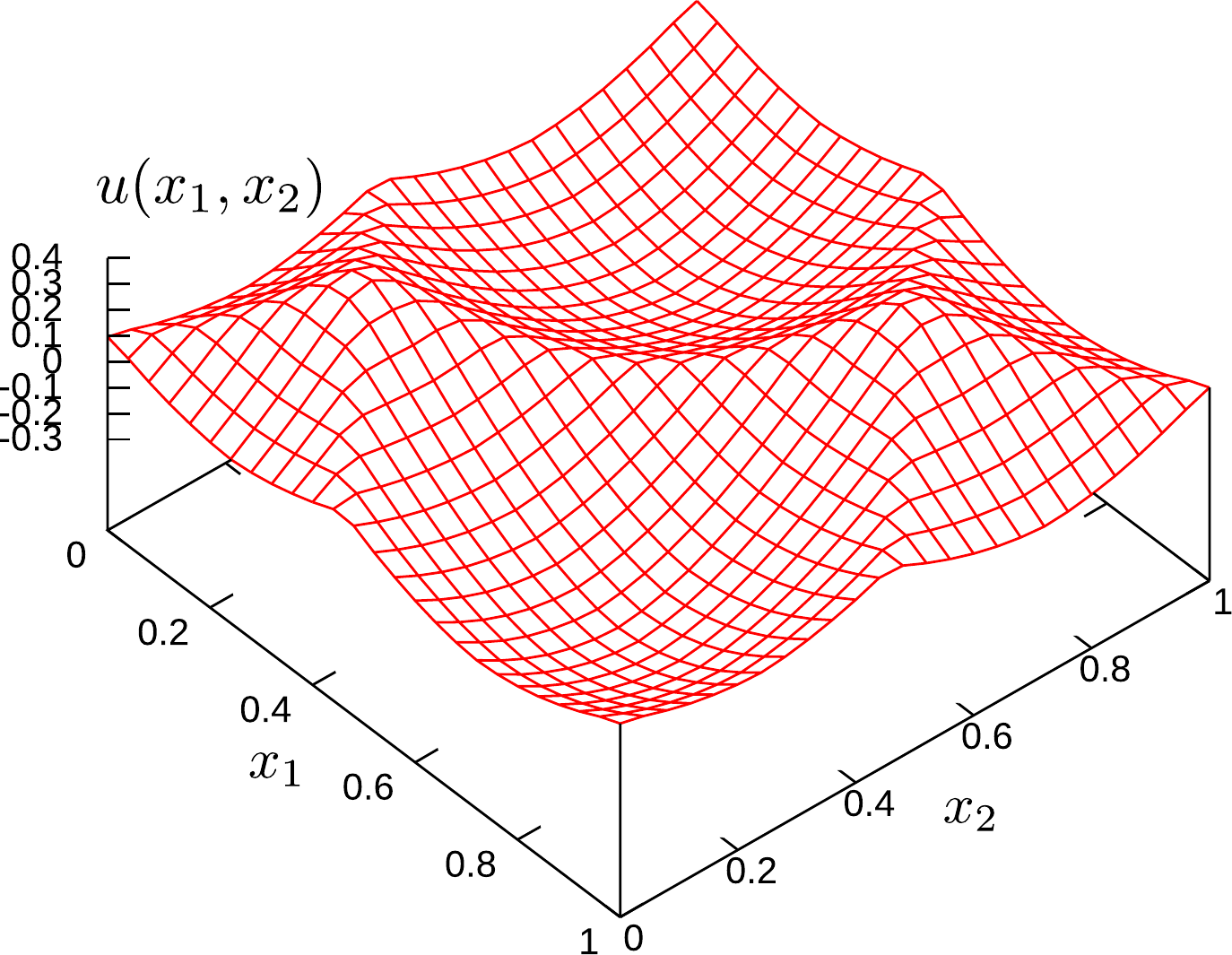} &
\includegraphics[width=.3\textwidth]{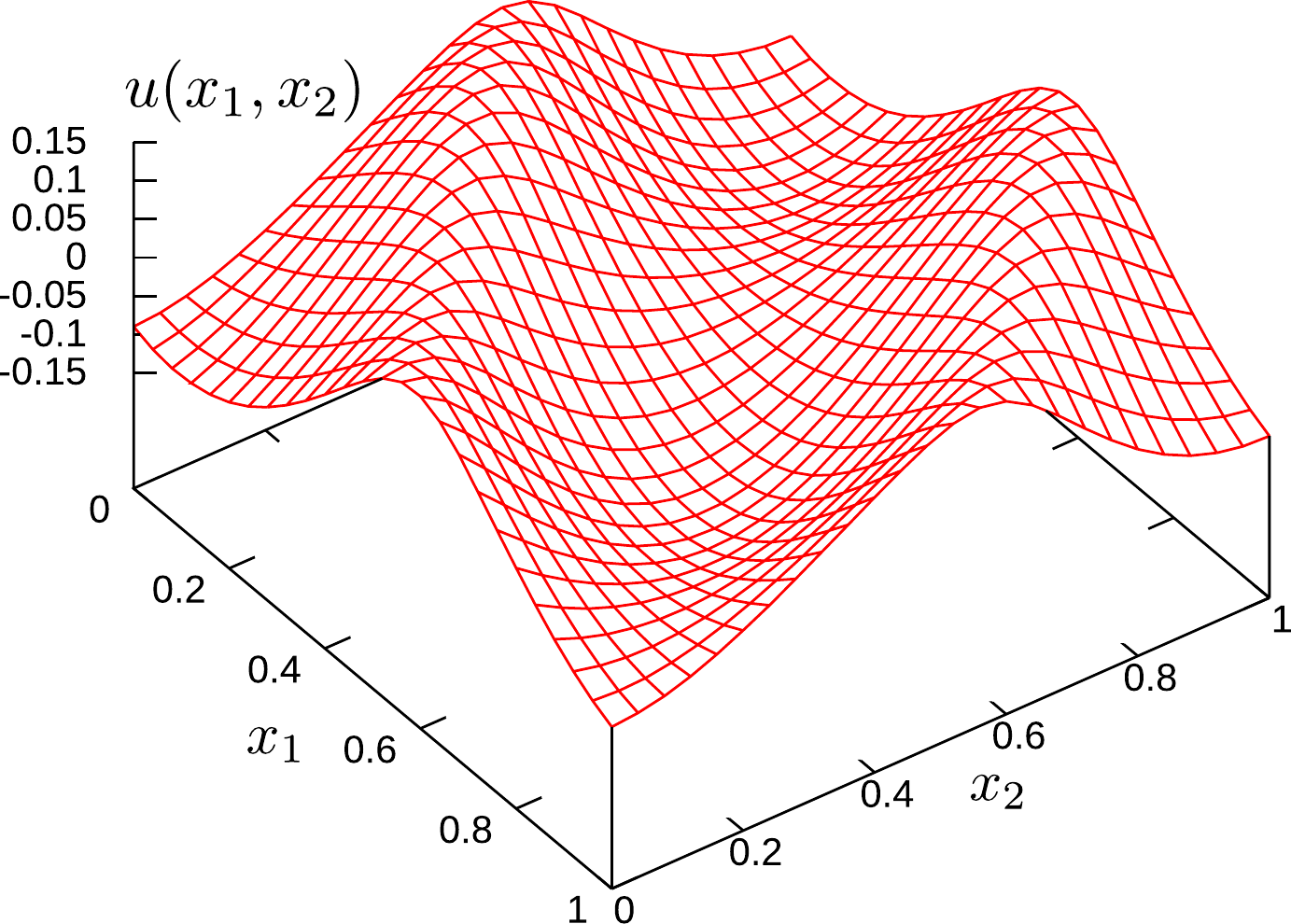} \\
(a)&(b)
\end{tabular}
\end{center}
\caption{Viscosity solution of the cell problem for $p=(0,0)$ (a), $p=(2,1.5)$ (b).}\label{correctors-2d-Q53}
\end{figure}

We proceed with the second example by choosing
\begin{equation}\label{ex2d-2}
V(x_1,x_2)=\sin(2\pi x_1)\sin(2\pi x_2)\,.
\end{equation}
In this case formula \eqref{hbarsep} no longer holds. Figure \ref{Hbar-2d-Q53} shows the effective Hamiltonian surface and its level sets. 
The computed value inside the plateau is $\bar H(p)=1$.
\\
In Figure \ref{correctors-2d-Q53} we show a couple of correctors for different values of $p$.
\\
 For a single point $p$, the average number of iterations is about $7$ and
the average computational time  is about $0.2$ seconds, whereas the total time is $480.75$ seconds.

\begin{figure}[h!]
 \begin{center}
 \begin{tabular}{cc}
\includegraphics[width=.3\textwidth]{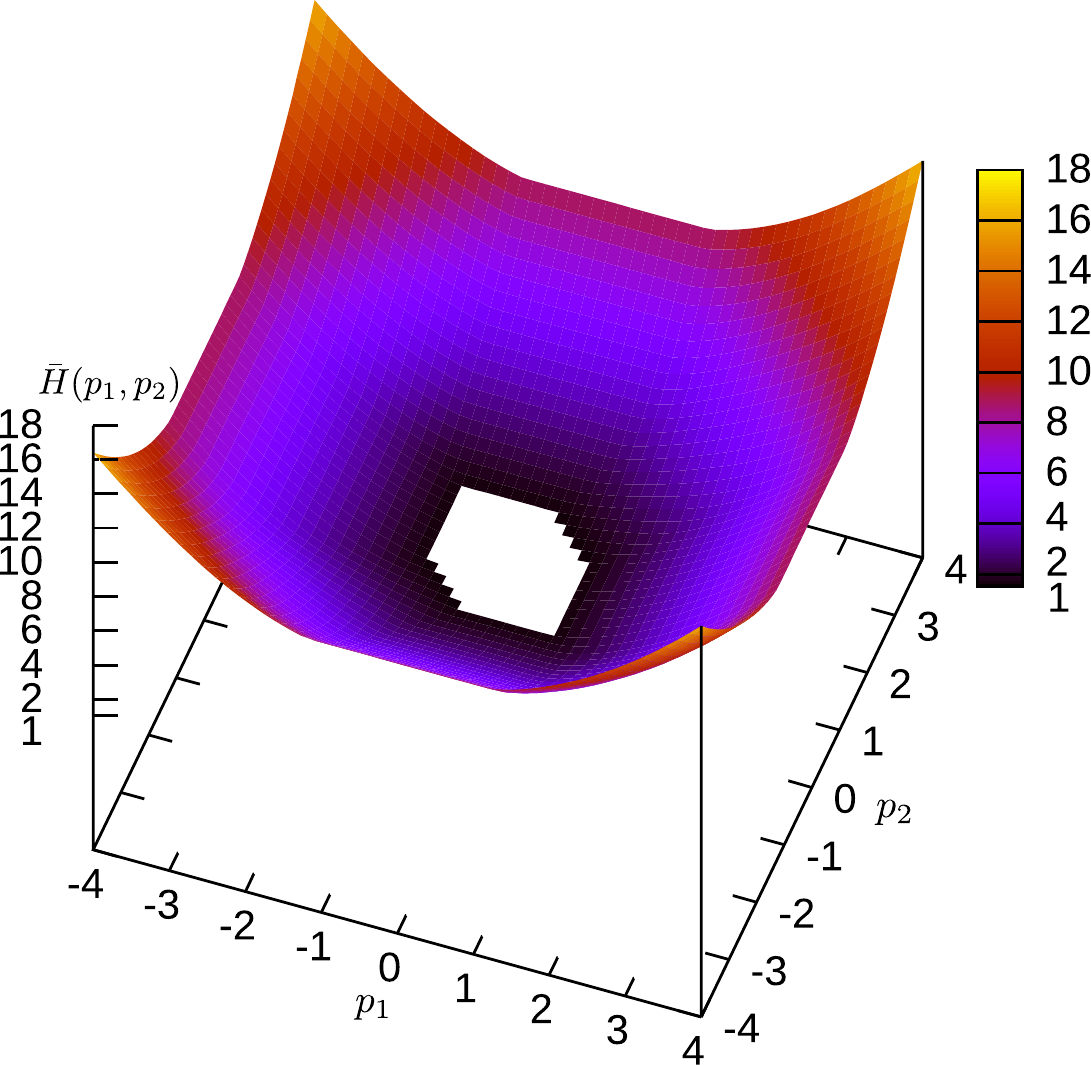} &
\includegraphics[width=.25\textwidth]{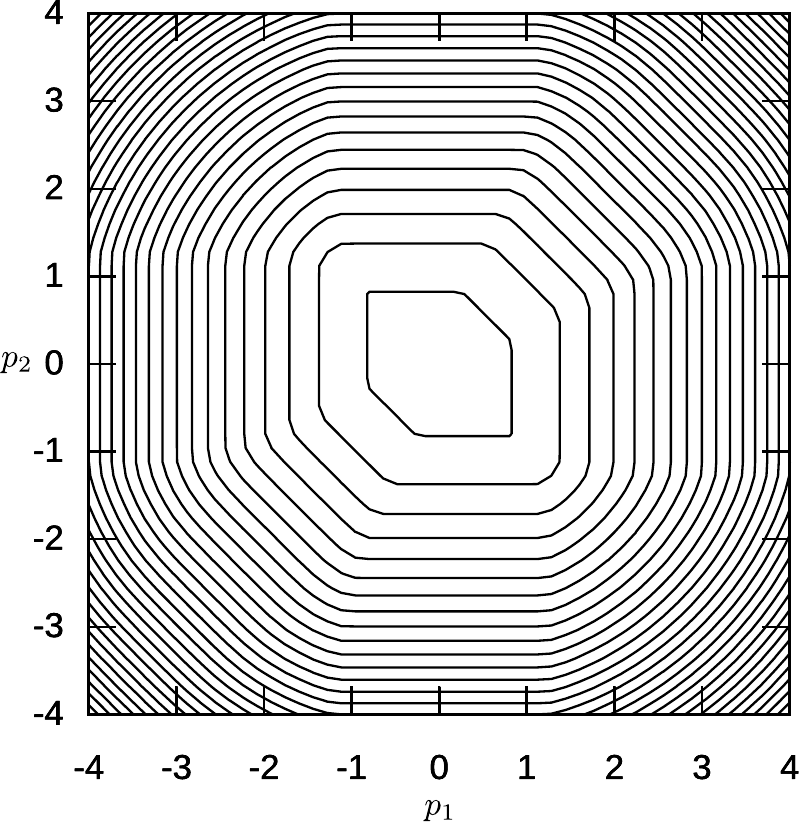} \\
(a)&(b)
\end{tabular}
\end{center}
\caption{Effective Hamiltonian for $p\in[-4,4]^2$: (a) surface and (b) level sets.}\label{Hbar-2d-Q54}
\end{figure}

\begin{figure}[h!]
 \begin{center}
 \begin{tabular}{cc}
\includegraphics[width=.3\textwidth]{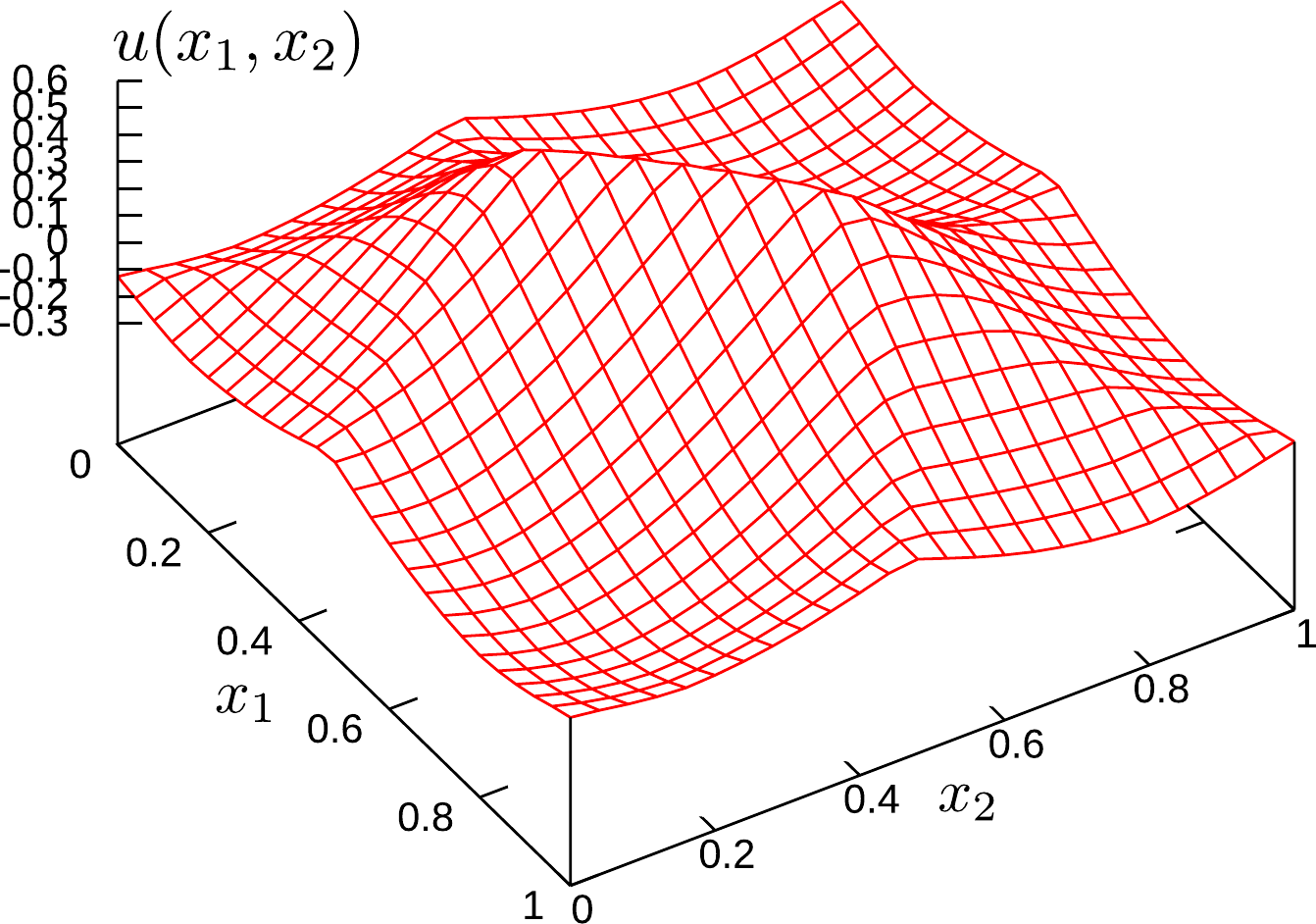} &
\includegraphics[width=.3\textwidth]{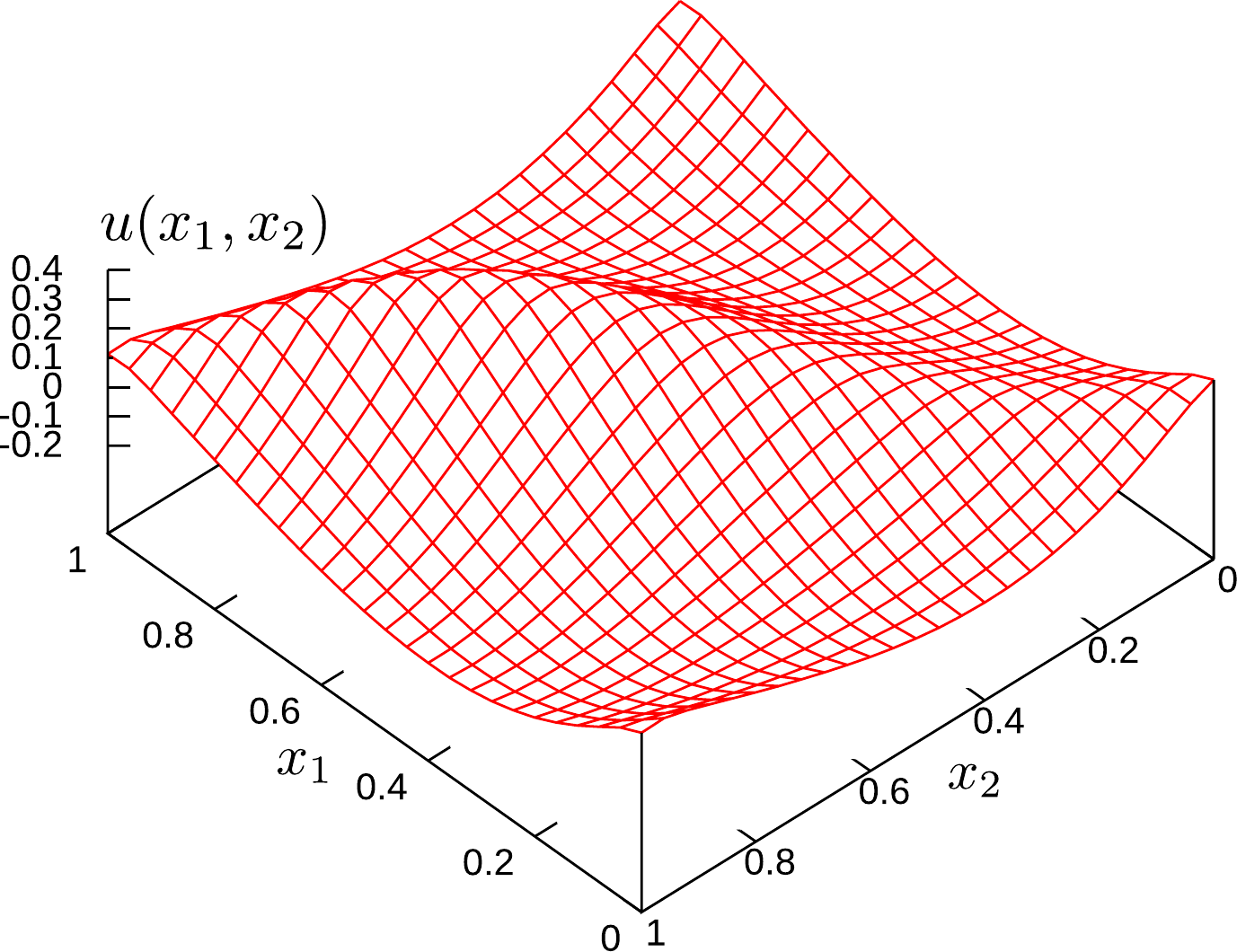} \\
(a)&(b)
\end{tabular}
\end{center}
\caption{Viscosity solution of the cell problem for $p=(0,0)$ (a), $p=(2,2)$ (b).}\label{correctors-2d-Q54}
\end{figure}

We finally consider the case
\begin{equation}\label{ex2d-3}
V(x_1,x_2)=\cos(2\pi x_1)+\cos(2\pi x_2)+\cos\left(2\pi (x_1-x_2)\right)\,.
\end{equation}
Figure \ref{Hbar-2d-Q54} shows the effective Hamiltonian surface and its level sets.
The computed value inside the plateau $\bar H(p)=1.488983$.
For a single point $p$, the average number of iterations is about $10$ and the average computational time  is about $0.24$ seconds,
whereas the total time is $630.77$ seconds. In Figure \ref{correctors-2d-Q54} we show a couple of correctors for different values of $p$.

We conclude this section summarizing the results in Table \ref{table-2d-norm2}, where we report
the average time and the average number of iterations for a single point $p$, and the total computational time for the whole simulation.
\begin{table}[!ht]
 \centering
  \begin{tabular}{|c|c|c|c|}
    \hline
    $V(x_1,x_2)$ & Av. CPU (secs) per $p$ & Av. Iterations per $p$ & Total CPU (secs)\\
    \hline
    \eqref{ex2d-1}&0.4 &16 &970.45\\\hline
    \eqref{ex2d-2}&0.2 & 7 &480.75\\\hline
    \eqref{ex2d-3}&0.24 &10 &630.77\\\hline
 \end{tabular}
  \caption{Performance for the 2D eikonal case.}\label{table-2d-norm2}
\end{table}\\
We observe that the perfomance of the method mainly depends on the extension of the plateau,
since it requires more iterations (and CPU time) to compute the corresponding corrector.

\section{First order convex Hamiltonians}\label{qnorms}
We consider for $q\ge 1$ the   cell problem on $\mathbb{T}^n$ (for $n=1,2$)
$$\frac{1}{q}|Du+p|^q-V(x)=\lambda\,,$$
where, as discussed in Section \ref{newtonlike}, the singularity at the origin of the derivative of $|\cdot|^q$ for $1\leq q<2$ is handled by choosing,
in a nonsmooth-Newton fashion, an element of the sub-differential. Here, we simply choose $0$ if $Du+p=0$ at some point.

\begin{figure}[h!]
 \begin{center}
 \begin{tabular}{ccc}
\includegraphics[width=.3\textwidth]{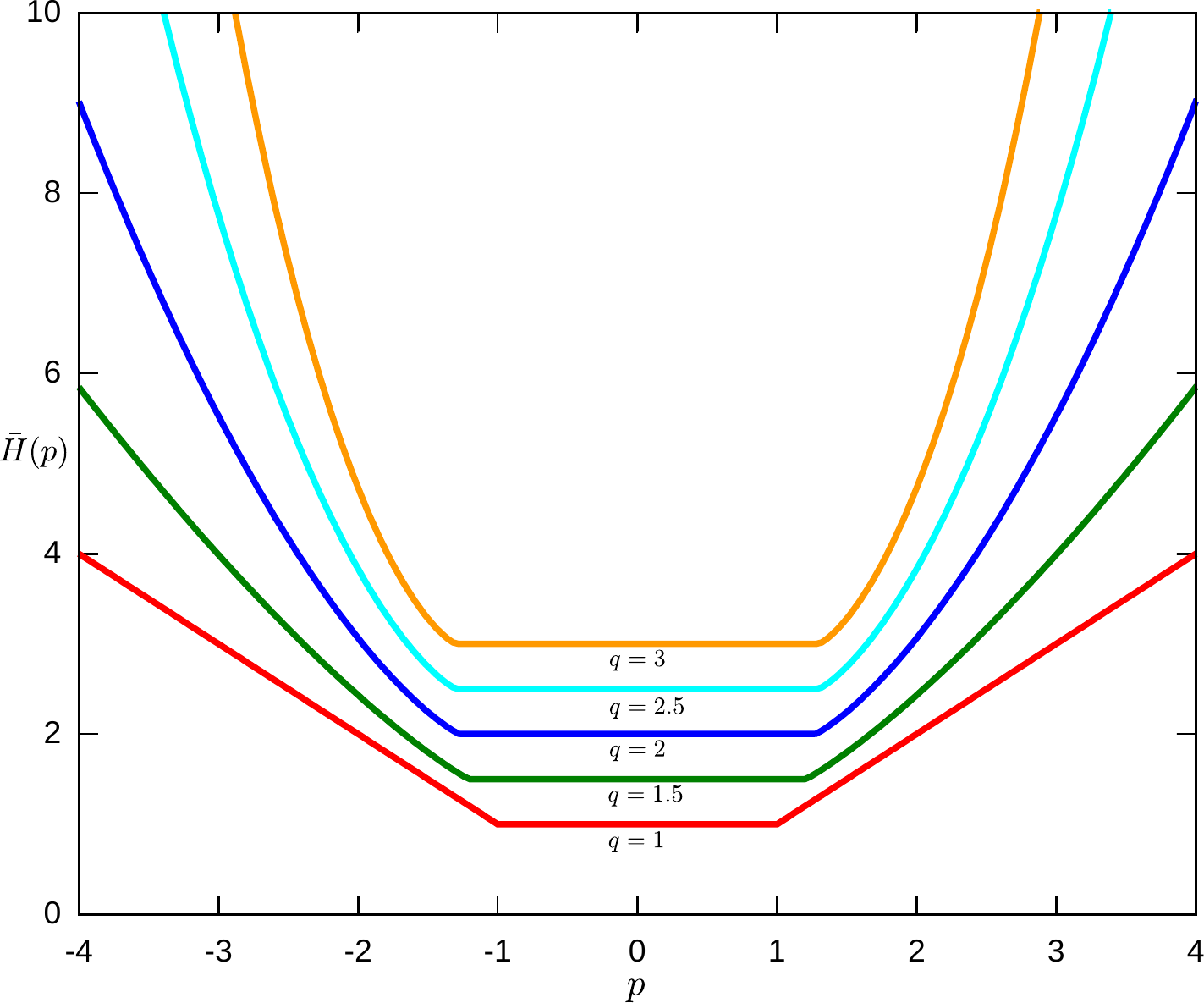} &
\includegraphics[width=.3\textwidth]{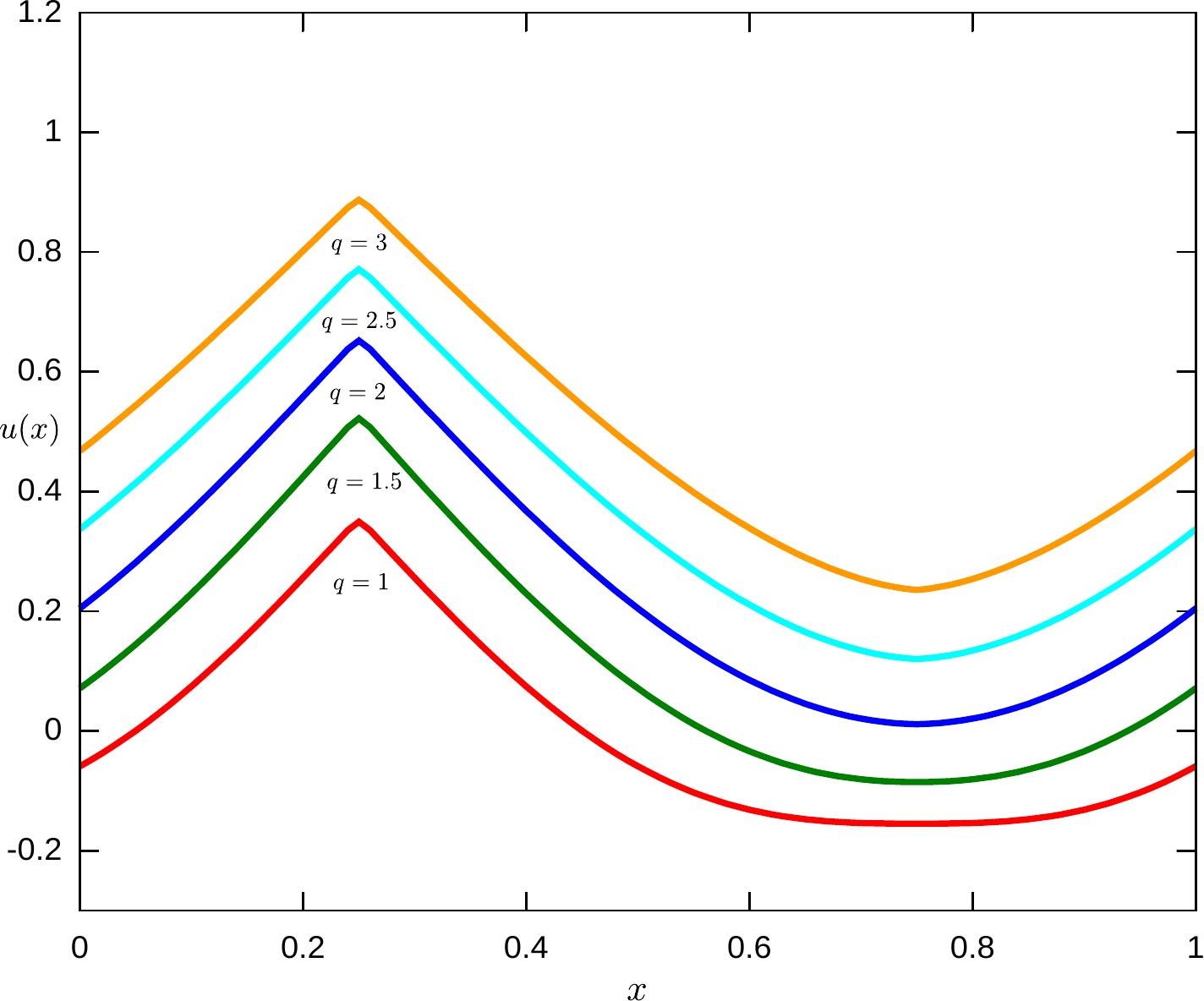} &
\includegraphics[width=.3\textwidth]{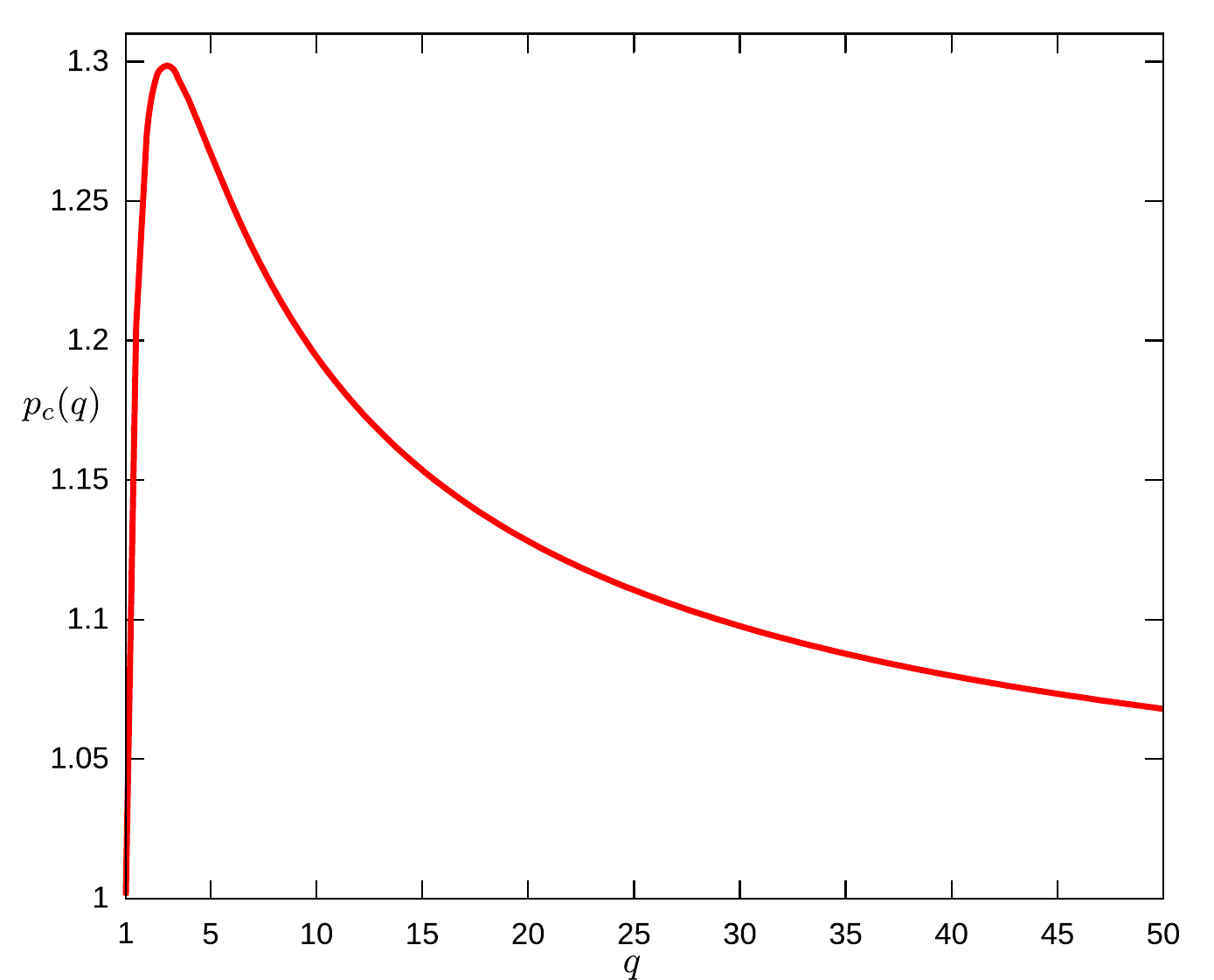} \\\vspace{5pt}
(a)&(b)&(c)
\end{tabular}
\end{center}
\caption{Effective Hamiltonians for $p\in[-4,4]$ and different $q$ (a) and corresponding correctors for $p=0$ (b). Extension of the plateau of $\bar H$: $p_c$ vs $q$ (c).}\label{Hbar-q-1d}
\end{figure}

In the one dimensional case, we consider again the potential $V(x)=\sin(2\pi x)$ and we compute the effective Hamiltonian for different values of $q$.
For each $q$ the computation takes about $2.5$ seconds and is performed choosing $p\in[-4,4]$ discretized with $201$ nodes.
Figure \ref{Hbar-q-1d}a shows the graphs of $\bar H$ and Figure \ref{Hbar-q-1d}b some corresponding correctors.
All the effective Hamiltonians are equal to $1$ in their
plateau, but we translated each graph of a fixed value in order to avoid overlapping and better appreciate the plateau itself.
We observe an interesting feature: the extension of the plateau of $\bar H$, or equivalently (in this case) the value $p_c=\inf\{p\ge 0\,:\,\bar H(p)>1\}$,
is not monotone with respect to $q$. This is confirmed by the graph in Figure \ref{Hbar-q-1d}c, showing $p_c$ as a function of $q$.
Each value $p_c$ is obtained by a simple bisection in $p$ applied to $\bar H(p)-1$ for $q$ ranging in $[1,50]$.
The maximum is achieved at $q^*=2.865$ with $p_c(q^*)=1.29876458$.\\

\begin{figure}[h!]
 \begin{center}
 \begin{tabular}{ccc}
\includegraphics[width=.25\textwidth]{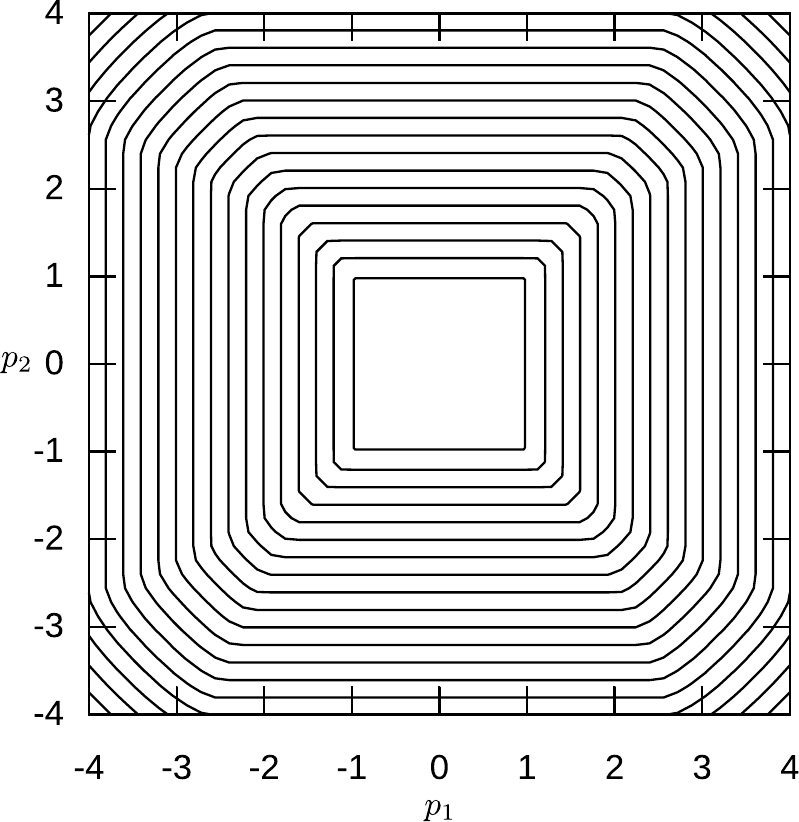} &
\includegraphics[width=.25\textwidth]{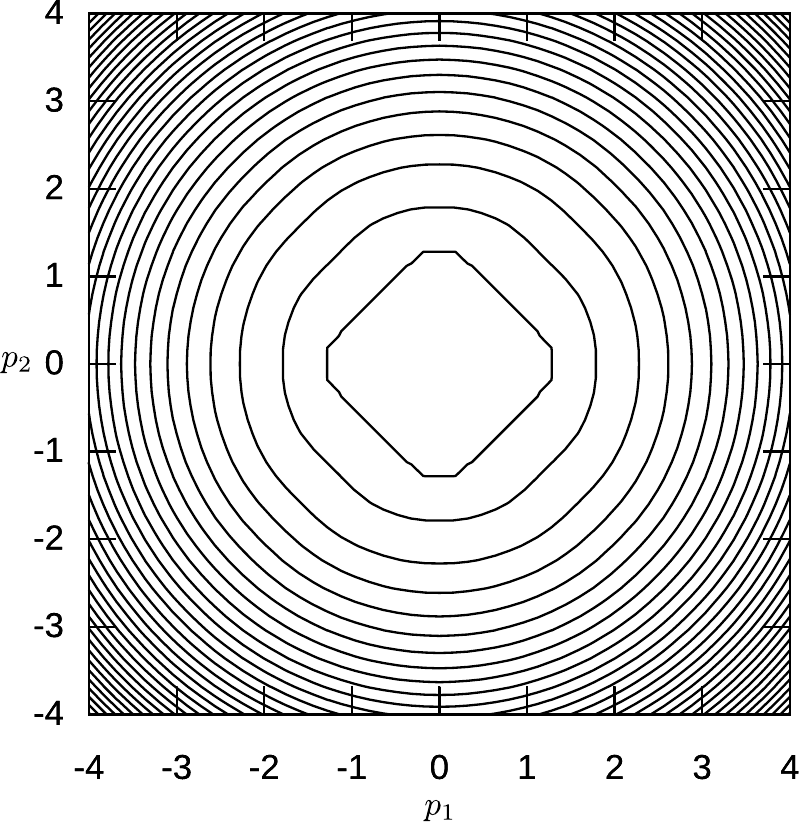} &
\includegraphics[width=.25\textwidth]{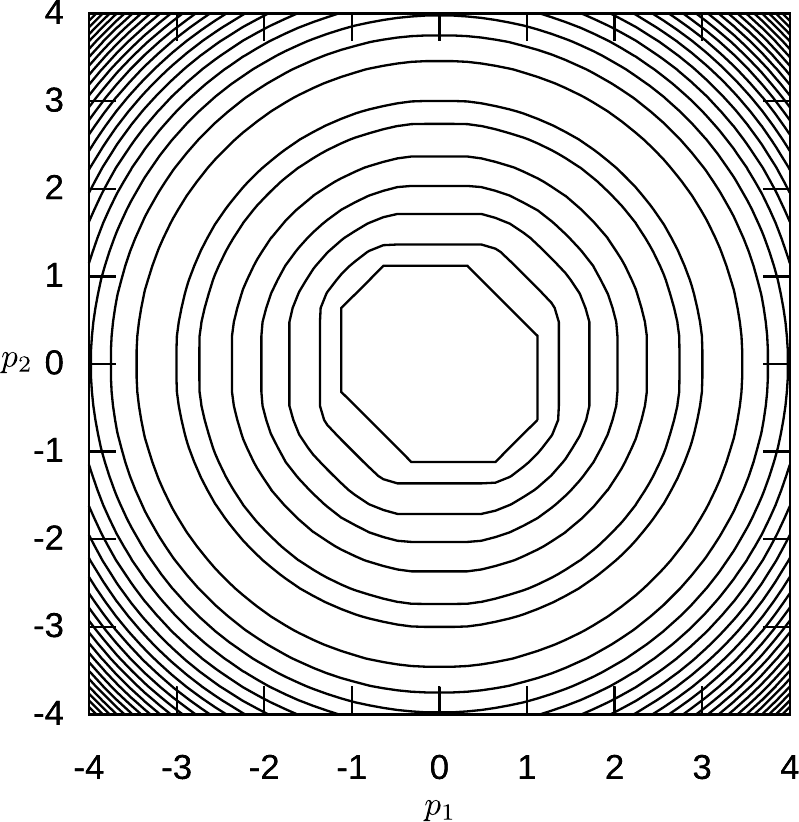} \\
\quad(a)&\quad(b)&\quad(c)
\end{tabular}
\end{center}
\caption{Level sets of the effective Hamiltonians for problems $(a)$, $(b)$, $(c)$.}\label{Hbar-2d-qnorms}
\end{figure}

We now consider a two dimensional case, choosing different values of $q$ and the same 2D potentials of the previous section. More precisely, we take\\
\begin{enumerate}
 \item[$(a)$] $q=1$ and $V(x_1,x_2)=\cos(2\pi x_1)+\cos(2\pi x_2)$\,,\\
 \item[$(b)$] $q=3$ and $V(x_1,x_2)=\sin(2\pi x_1)\sin(2\pi x_2)$\,,\\
 \item[$(c)$] $q=5$ and $V(x_1,x_2)=\cos(2\pi x_1)+\cos(2\pi x_2)+\cos\left(2\pi (x_1-x_2)\right)$\,.
\end{enumerate}
\medskip
\noindent We choose as before $p\in[-4,4]^2$ discretized with $51\times 51$ uniformly distributed nodes, while the mesh size for
the torus $\mathbb{T}^2$ is $25\times 25$ nodes.
In Figure \ref{Hbar-2d-qnorms} we show the level sets of the effective Hamiltonians.
Note the differences with respect to the analogous tests in the case $q=2$, already shown in Figure \ref{Hbar-2d-Q52}, Figure \ref{Hbar-2d-Q53} and Figure \ref{Hbar-2d-Q54} respectively.

Finally, in Table \ref{table-2d-qnorms} we report the results for the three tests, including again, for a single point $p$,
the average time and the average number of iterations, and the total computational time for the whole simulation.
Note that the most expensive case corresponds to $q=1$, in which the non-smoothness of the problem plays a crucial role.\vskip-2pt
\begin{table}[!ht]
 \centering
  \begin{tabular}{|c|c|c|c|}
    \hline
    Problem & Av. CPU (secs) per $p$ & Av. Iterations per $p$ & Total CPU (secs)\\
    \hline
    $(a)$&0.6 &28 &1633.12\\\hline
    $(b)$&0.2 & 9 & 522.37\\\hline
    $(c)$&0.4 &18 &1042.18\\\hline
 \end{tabular}
  \caption{Performance for the 2D $q\hbox{-}$power Hamiltonian case.}\label{table-2d-qnorms}
\end{table}
\section{First order nonconvex Hamiltonians}\label{nonconvex}
We consider the case of a nonconvex Hamiltonian in dimension one, namely the cell problem on $\mathbb{T}^1$
$$\frac12(|u^\prime+p|^2-1)^2-V(x)=\lambda\,.$$
In this setting, as for the one dimensional eikonal equation, a formula for the effective Hamiltonian is still available
(see \cite{Q} for details):
$$
\bar{H}(p)=\left\{
\begin{array}{ll}
 -\min V & \mbox{if } |p|\le p_c\\
 \lambda & \mbox{if } |p|> p_c\quad\mbox{s.t. }|p|=\displaystyle\int_0^1 \sqrt{1+\sqrt{2(V(s)+\lambda)}}ds
\end{array}
\right. \vspace{-10pt}
$$
where the plateau is given by $p_c=\displaystyle\int_0^1\sqrt{1+\sqrt{2(V(s)-\min V))}}ds$.\\

\begin{figure}[h!]
 \begin{center}
 \begin{tabular}{cc}
\includegraphics[width=.3\textwidth]{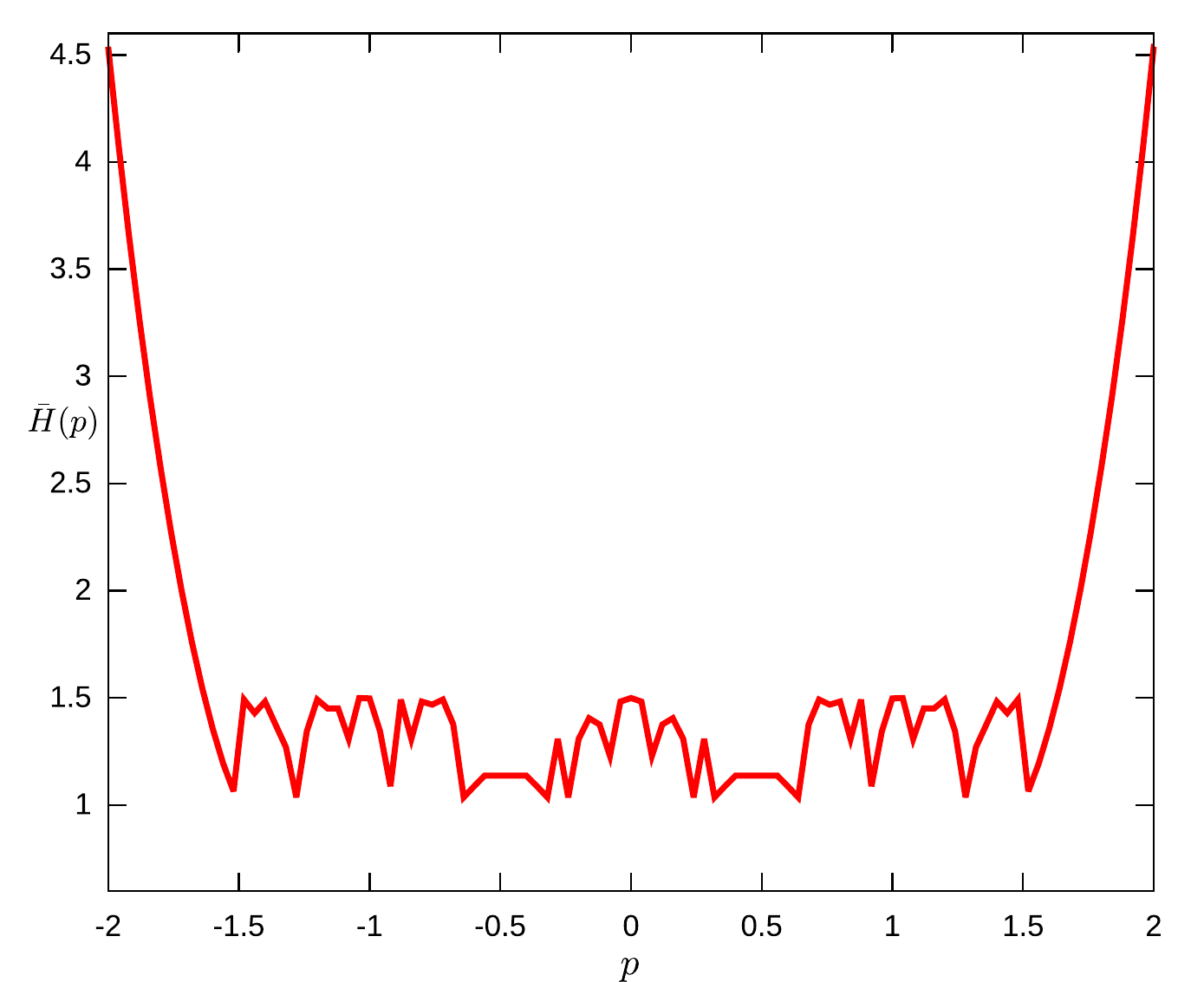} &
\includegraphics[width=.3\textwidth]{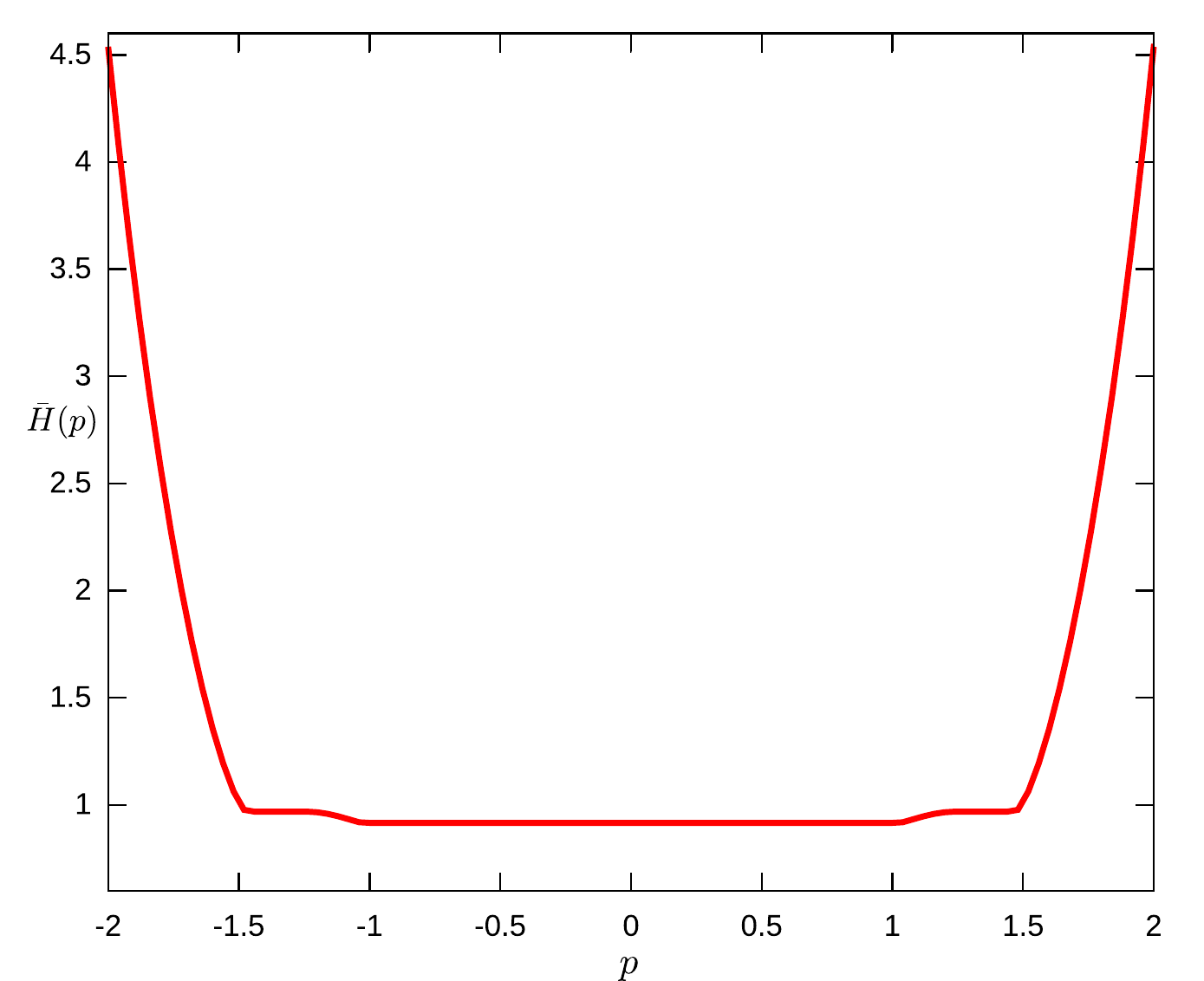} \\
(a)&(b)
\end{tabular}
\end{center}
\caption{Effective Hamiltonian for $p\in[-2:2]$: Engquist-Osher (a) and global Lax-Friedrichs (b).}\label{Hbar-nonconvex-1d}
\end{figure}

The crucial point of the nonconvex case is that viscosity solutions are allowed to have kinks pointing both upward and downward,
but it is known that the standard Engquist-Osher numerical Hamiltonian is not able to select them correctly as in the convex case.
Nevertheless, we run the algorithm and we show in Figure \ref{Hbar-nonconvex-1d}a the effective Hamiltonian computed in this situation for the potential $V(x)=\sin(2\pi x)$.
Our method still works, i.e., it converges to a solution $(U,\Lambda)$ with zero residual, but the result is completely wrong in the plateau ($p_c\sim 1.4918$),
where solutions with kinks are expected.
On the other hand, outside the plateau, the viscosity solution is smooth and the result is correct.
To proceed, we employ the well known global Lax-Friedrichs numerical Hamiltonian, which is very easy to implement
and can handle general Hamiltonians. The only drawback
is a loss of accuracy introduced by the global artificial viscosity term.

Figure \ref{Hbar-nonconvex-1d}b shows the computed effective Hamiltonian, which is qualitatively in agreement with that computed in \cite{Q}.
The effect of the artificial viscosity is evident in the plateau of $\bar H$ where a small bump appears.
\\
Both simulations above were performed
for $p$ in the interval $[-2,2]$ discretized with $101$ uniformly distributed nodes, while the mesh size for
the torus $\mathbb{T}^1$ is again $100$ nodes.
The first simulation toke $3.15$ seconds with an average number of iterations equal to $38$, whereas the second simulation toke
$8.97$ seconds with an average number of iterations equal to $126$.

\section{Second order Hamiltonians}\label{2ndorder}
We consider the  following cell problem for the homogenization of fully nonlinear second order Hamilton-Jacobi equations:
\begin{equation}\label{CP2order}
     H(x, p , D^2u +s)=\l,\qquad\qquad x\in \T^n,
\end{equation}
for $p\in\R^n$, $\l\in\R$ and $s\in \mathcal S^n$, where $\mathcal S^n$ is the space of symmetric $n\times n$ matrices.
Assuming that  $H$ is continuous  and uniformly  elliptic, then  there exists a unique $\l=\bar H(p,s)$ and a unique (up to a constant)
$u$ such that the cell  problem admits a viscosity solution (see \cite{G} for details).
A finite difference approximation of cell problem \eqref{CP2order} is discussed in \cite{CM}, where a convergence result
for the approximation of the effective Hamiltonian $\bar H(p,s)$ is given.

\begin{figure}[h!]
 \begin{center}
 \begin{tabular}{ccc}
\includegraphics[width=.25\textwidth]{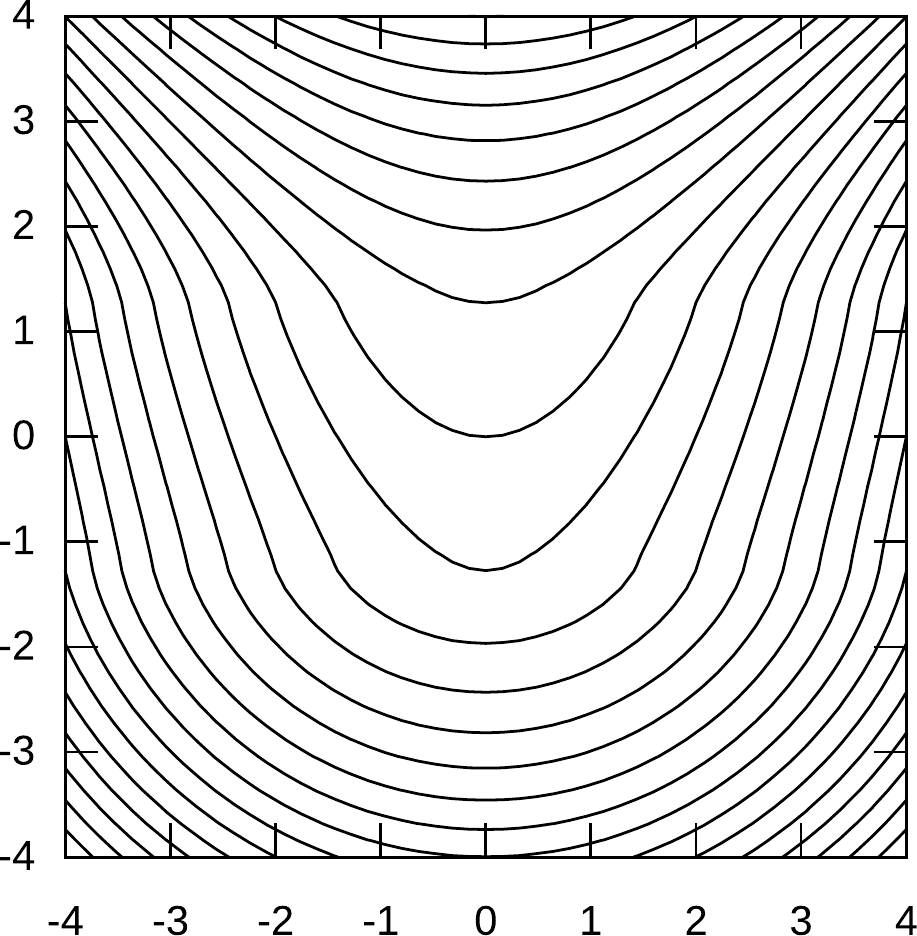} &
\includegraphics[width=.25\textwidth]{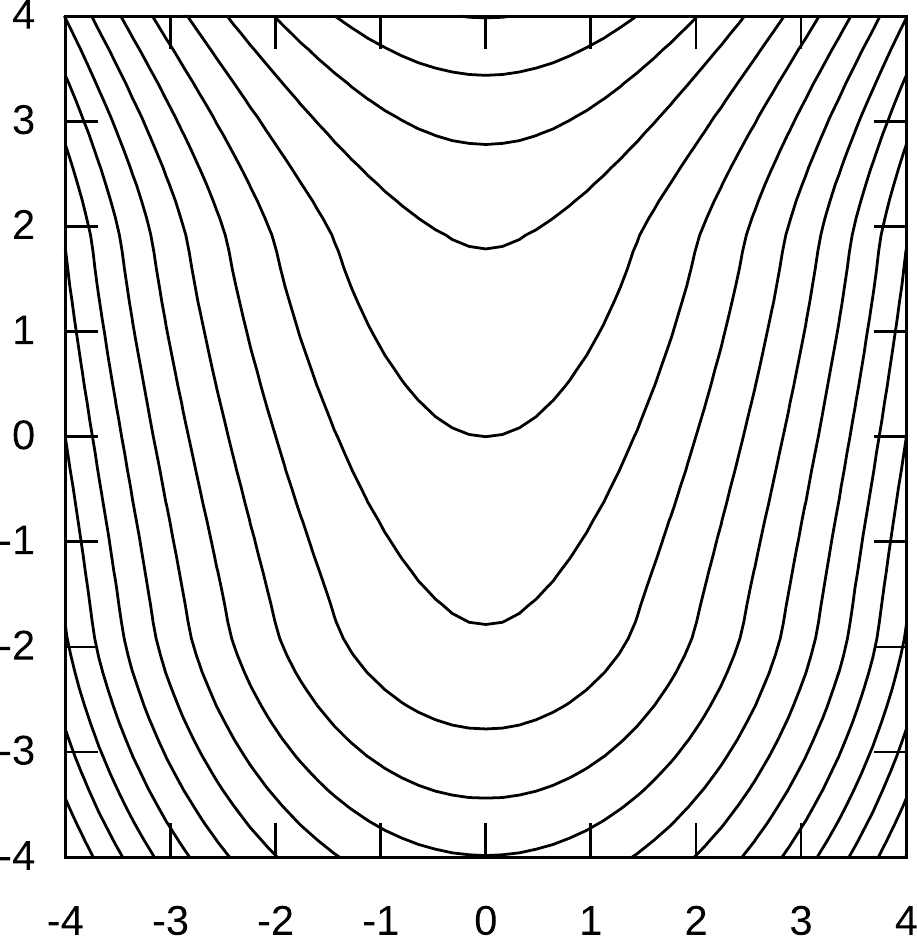} &
\includegraphics[width=.25\textwidth]{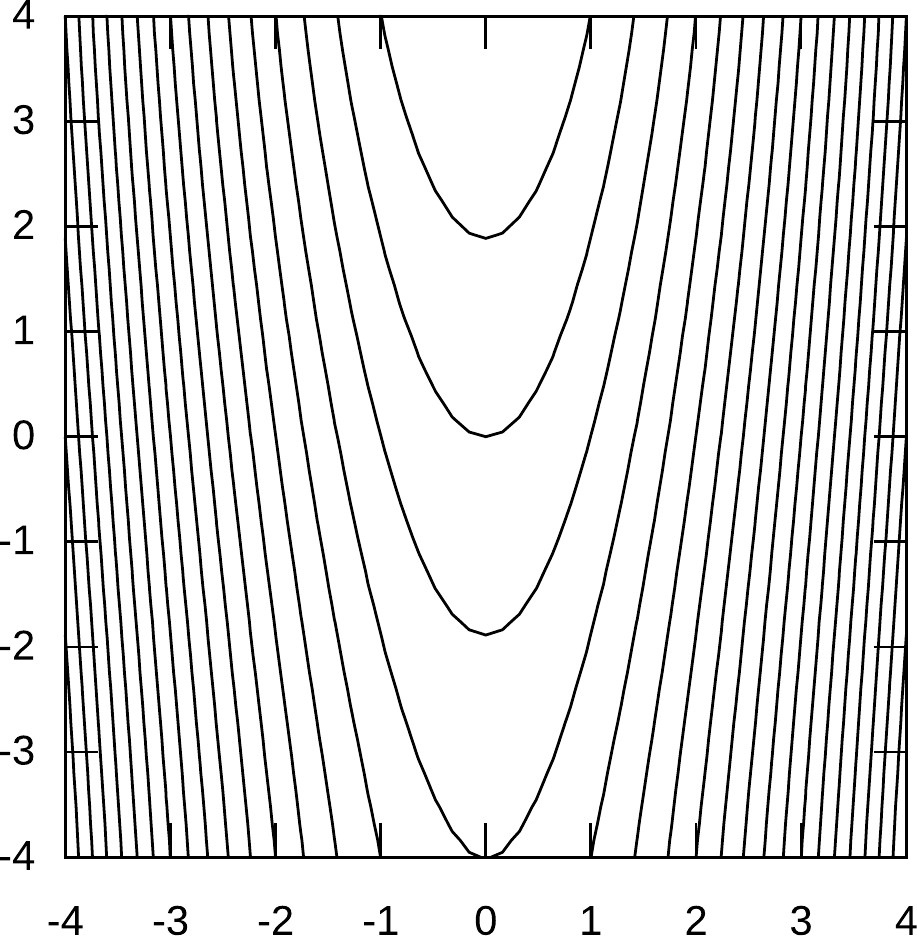} \\
(a)&(b)&(c)
\end{tabular}
\end{center}
\caption{Effective Hamiltonian level sets for $(p,s)\in[-2,2]^2$: $\alpha=1$ (a), $\alpha=\frac12$ (b), $\alpha=\frac{1}{10}$ (c).}\label{hbar1d2ndorder-levels}
\end{figure}

\begin{figure}[h!]
 \begin{center}
 \begin{tabular}{ccc}
\includegraphics[width=.25\textwidth]{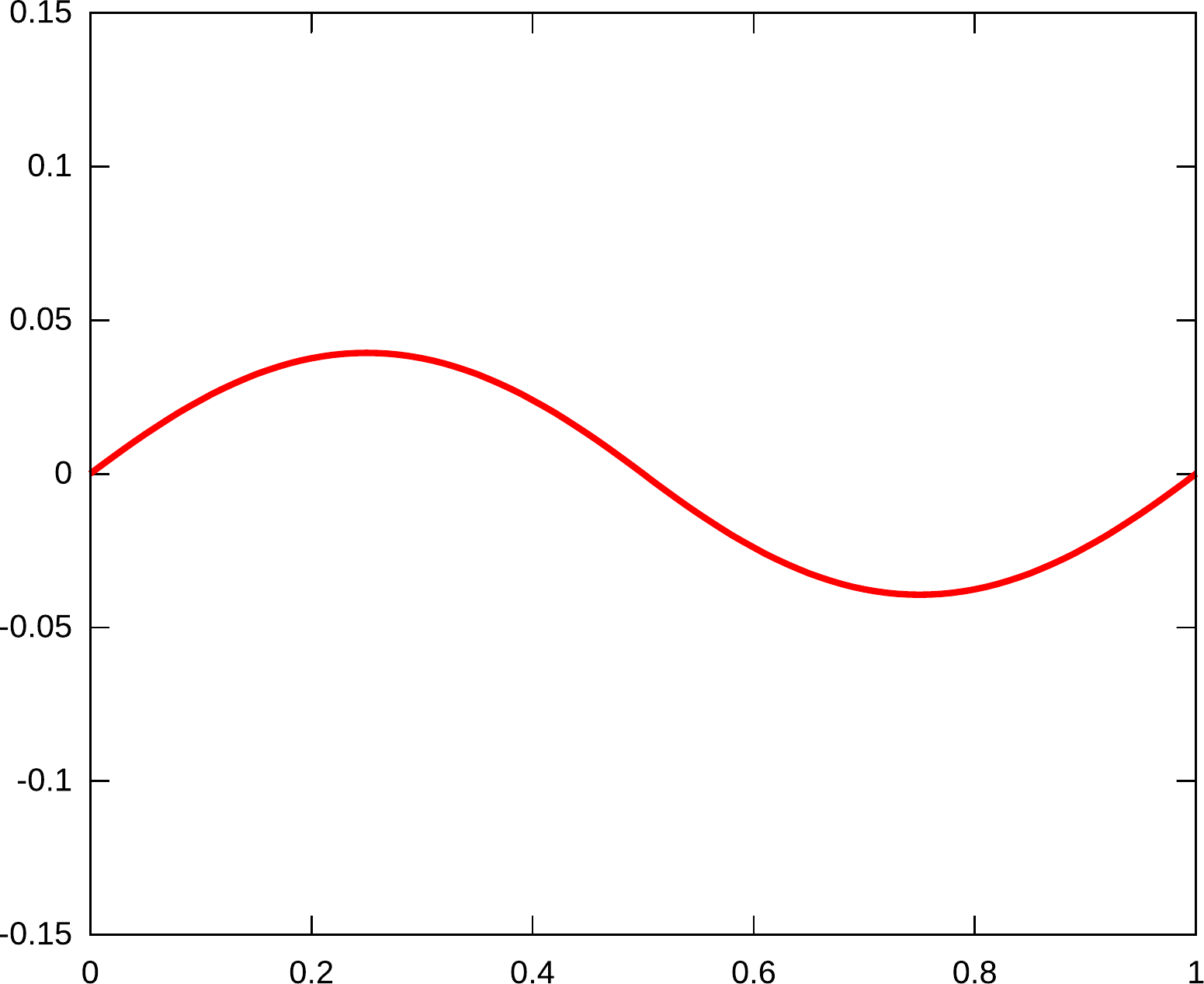} &
\includegraphics[width=.25\textwidth]{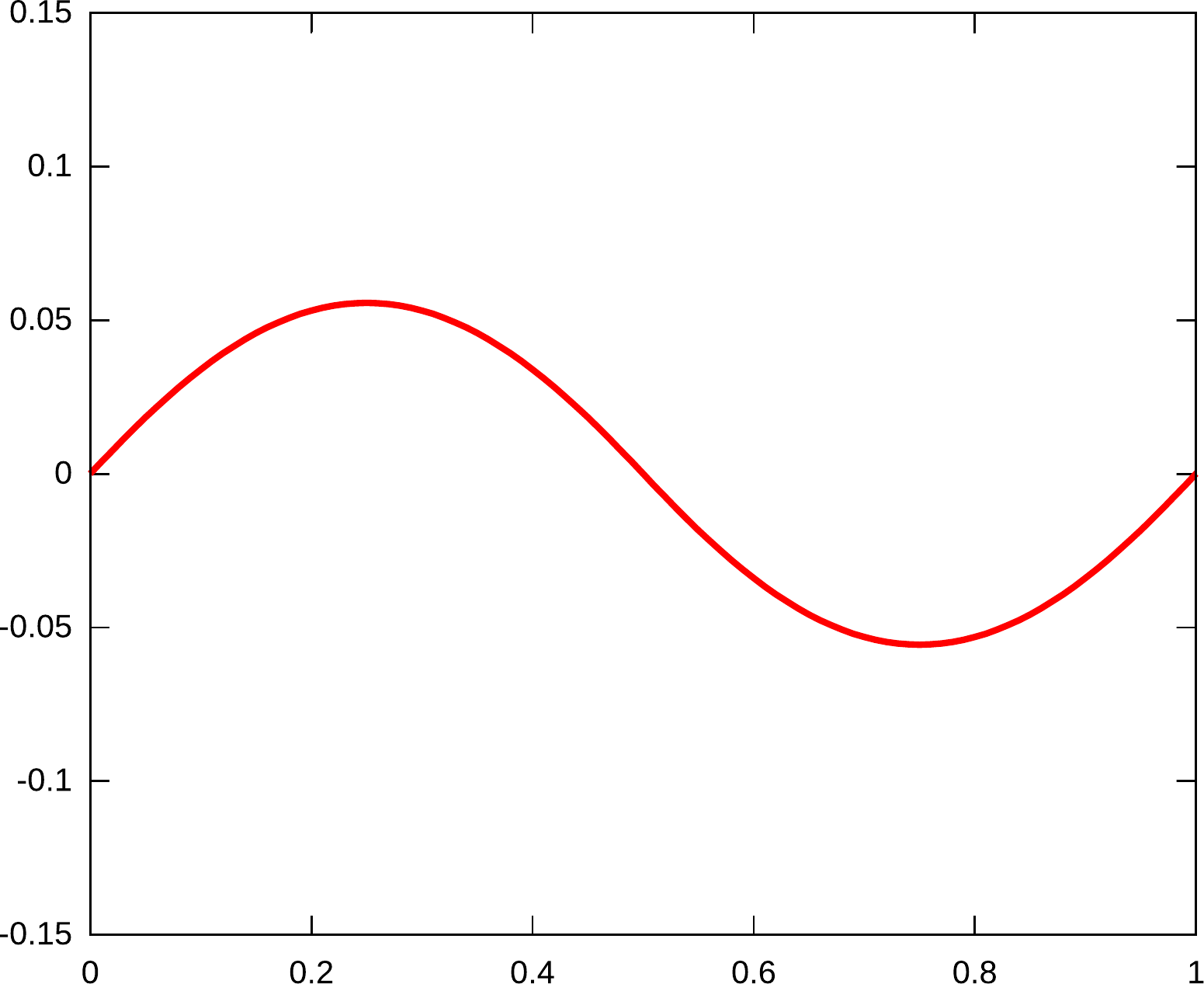} &
\includegraphics[width=.25\textwidth]{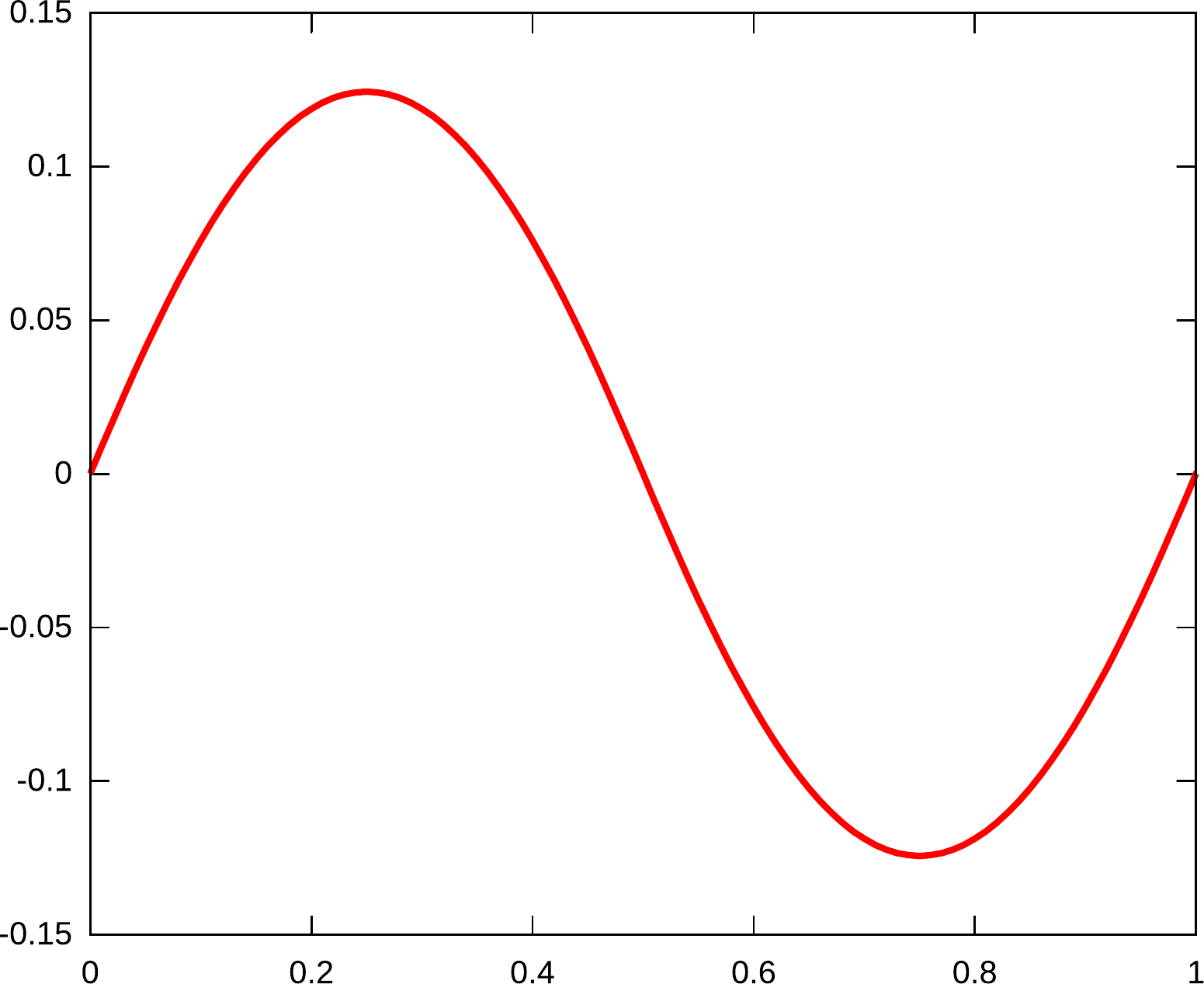} \\
(a)&(b)&(c)
\end{tabular}
\end{center}
\caption{Solutions of the cell problem for $(p,s)=(0,0)$: $\alpha=1$ (a), $\alpha=\frac12$ (b), $\alpha=\frac{1}{10}$.}\label{hbar1d2ndorder-correctors}
\end{figure}

Here we consider for simplicity the case of a fully nonlinear second order  Hamiltonian in dimension one, namely the cell problem on $\mathbb{T}^1$
$$-\alpha|u^{\prime\prime}+s|(u^{\prime\prime}+s)+\frac12|p|^2-V(x)=\lambda\,,$$
where $p,s\in\R$ and $\alpha> 0$. We choose $(p,s)\in[-4,4]^2$, discretized with $51\times 51$ uniformly distributed nodes,
while the mesh size for the torus $\mathbb{T}^1$ is $100$ nodes.
\begin{figure}[h!]
 \begin{center}
 \begin{tabular}{ccc}
\includegraphics[width=.3\textwidth]{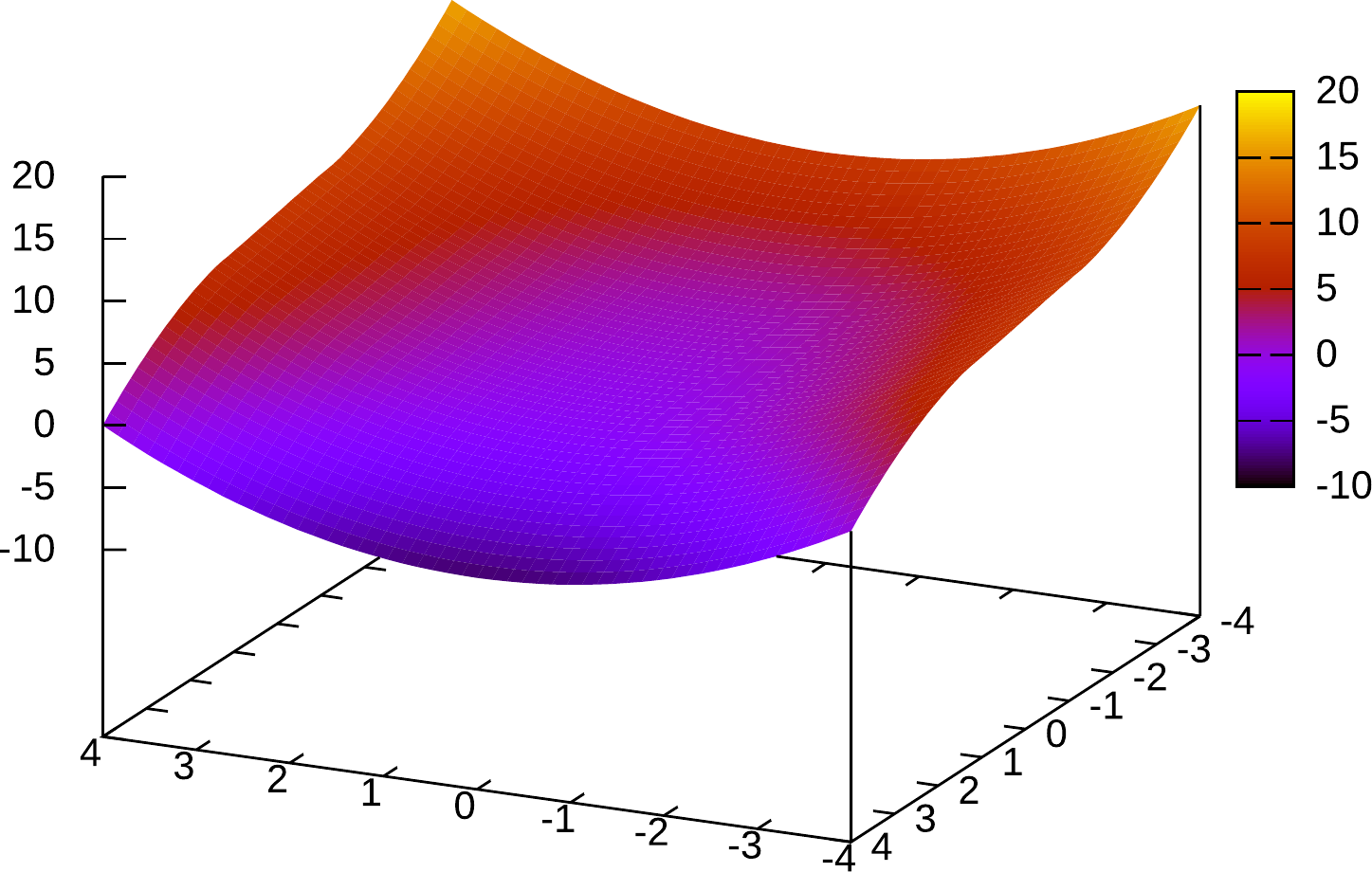} &
\includegraphics[width=.3\textwidth]{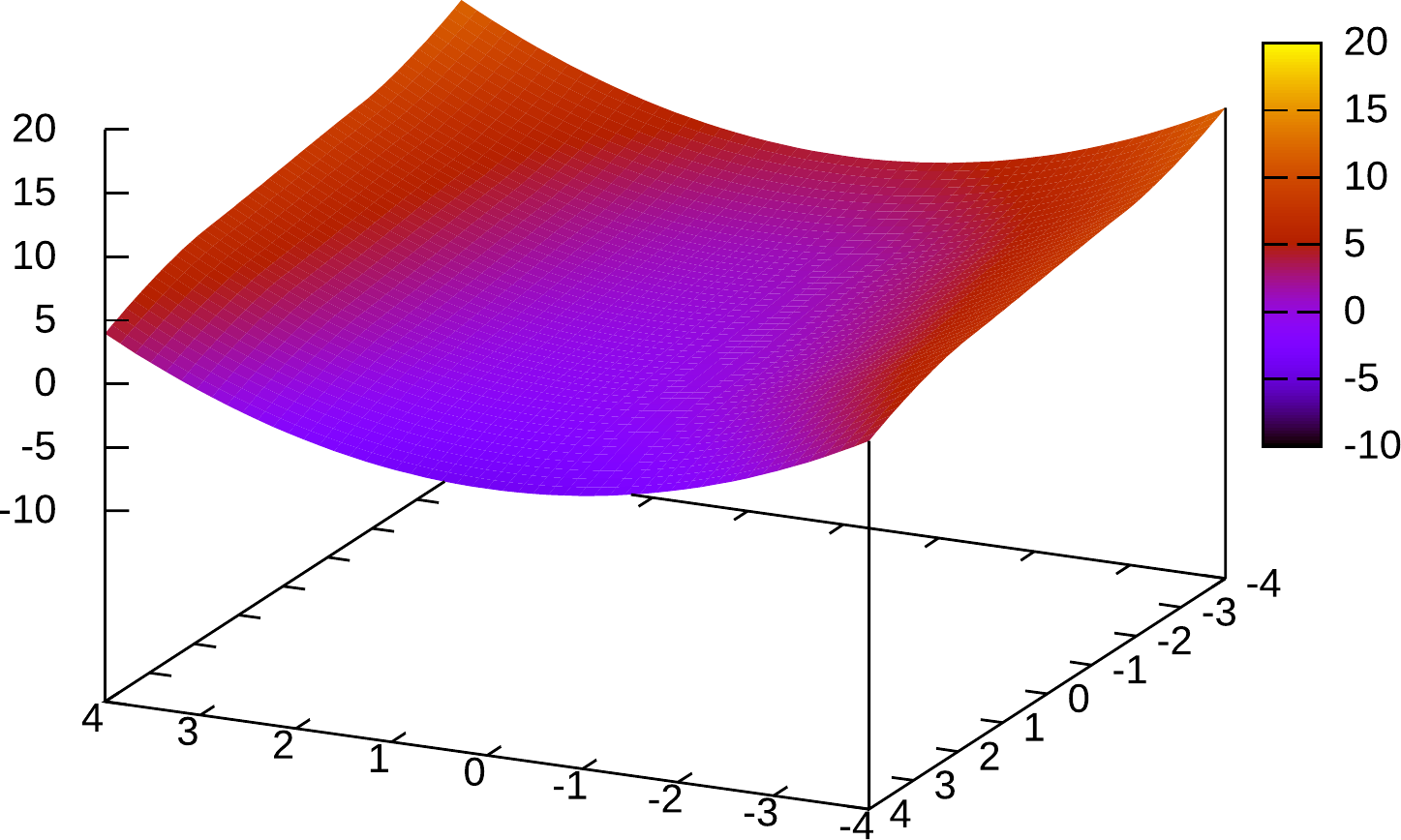} &
\includegraphics[width=.3\textwidth]{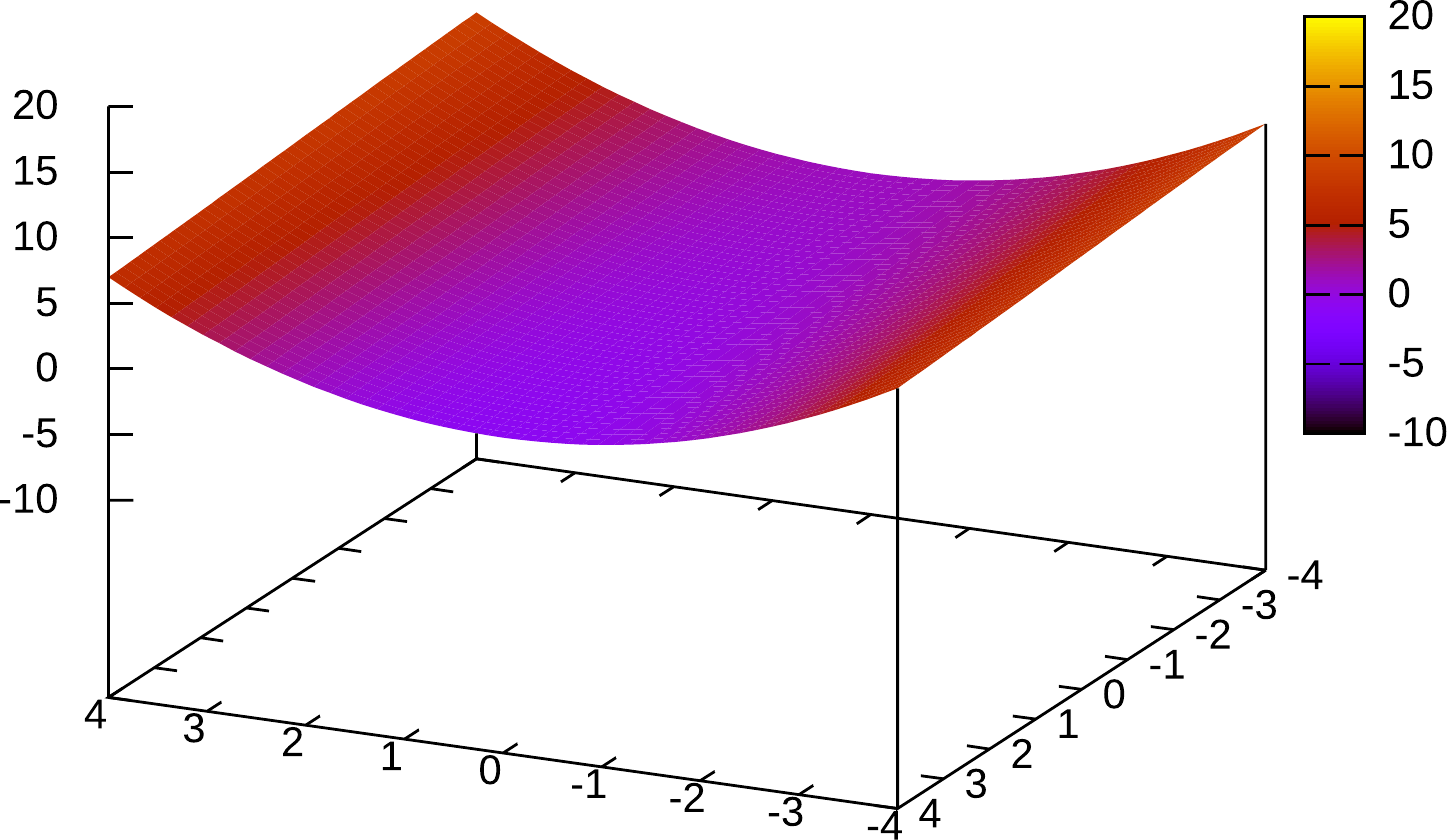} \\
(a)&(b)&(c)
\end{tabular}
\end{center}
\caption{Effective Hamiltonian surface for $(p,s)\in[-2,2]^2$: $\alpha=1$ (a), $\alpha=\frac12$ (b), $\alpha=\frac{1}{10}$ (c).}\label{hbar1d2ndorder-plot}
\end{figure}\\
In Figure \ref{hbar1d2ndorder-plot} we show the surfaces of the computed effective Hamiltonians for different values of $\alpha$, while
in Figure \ref{hbar1d2ndorder-levels} we show the corresponding level sets and in Figure \ref{hbar1d2ndorder-correctors} the corresponding
correctors for $(p,s)=(0,0)$.
\\
Finally, in Table \ref{table-1d-2nd} we report the performance of the method depending on $\alpha$.
\begin{table}[!ht]
 \centering
  \begin{tabular}{|c|c|c|c|}
    \hline
    $\alpha$ & Av. CPU (secs) per $(p,s)$ & Av. Iterations per $(p,s)$ & Total CPU (secs)\\
    \hline
    $1$&0.009 &7 &25.29\\\hline
    $0.5$&0.009 & 8 &25.26\\\hline
    $0.1$&0.011 &10 &30.78\\\hline
   \end{tabular}
  \caption{Performance for the 1D second order case.}\label{table-1d-2nd}
\end{table}

\section{Weakly coupled first order systems}\label{WCS}
In homogenization and long-time behavior of weakly coupled systems of Hamilton-Jacobi equations, the cell problem is given by
\begin{equation}\label{CPsystem}
    H_i(x,  Du_i+p)+C(x) u=\l,\qquad\qquad x\in \T^n,\quad i=1\dots,M
\end{equation}
where $u=(u_1,\dots,u_M)$ is a vector function, $p\in\R^n$, $\l\in\R$, the Hamiltonians $H_i$ are continuous and coercive and the $M\times M$  coupling matrix $C(x)=\{c_{ij}(x)\}_{i,j}$ is continuous, irreducible and satisfies
\[
c_{ij}(x)\le 0\,\,\text{for $j\neq i$},\qquad \sum_{j=1}^M c_{ij}(x)=0,\quad i=1,\dots,M.
\]
For a complete study of \eqref{CPsystem} we refer to \cite{DZ}. It is possible to prove that
also in this case there exists a unique $\l$ such that \eqref{CPsystem} admits a viscosity solution  $u$, which is  in general not unique.

For simplicity, we consider here the case of only two weakly coupled eikonal equations, namely the following cell problem on $\mathbb{T}^n$ ($n=1,2$)
$$\left\{\begin{array}{l}
\frac12|Du_1+p|^2-V_1(x)+c_1(x)(u_1-u_2)=\lambda\\\\
\frac12|Du_2+p|^2-V_2(x)+c_2(x)(u_2-u_1)=\lambda
\end{array}\right.
$$
for two $1\hbox{-}$periodic potentials $V_1$, $V_2$ and two nonnegative $1\hbox{-}$periodic functions $c_1$, $c_2$.

\begin{figure}[h!]
 \begin{center}
 \begin{tabular}{cc}
\includegraphics[width=.3\textwidth]{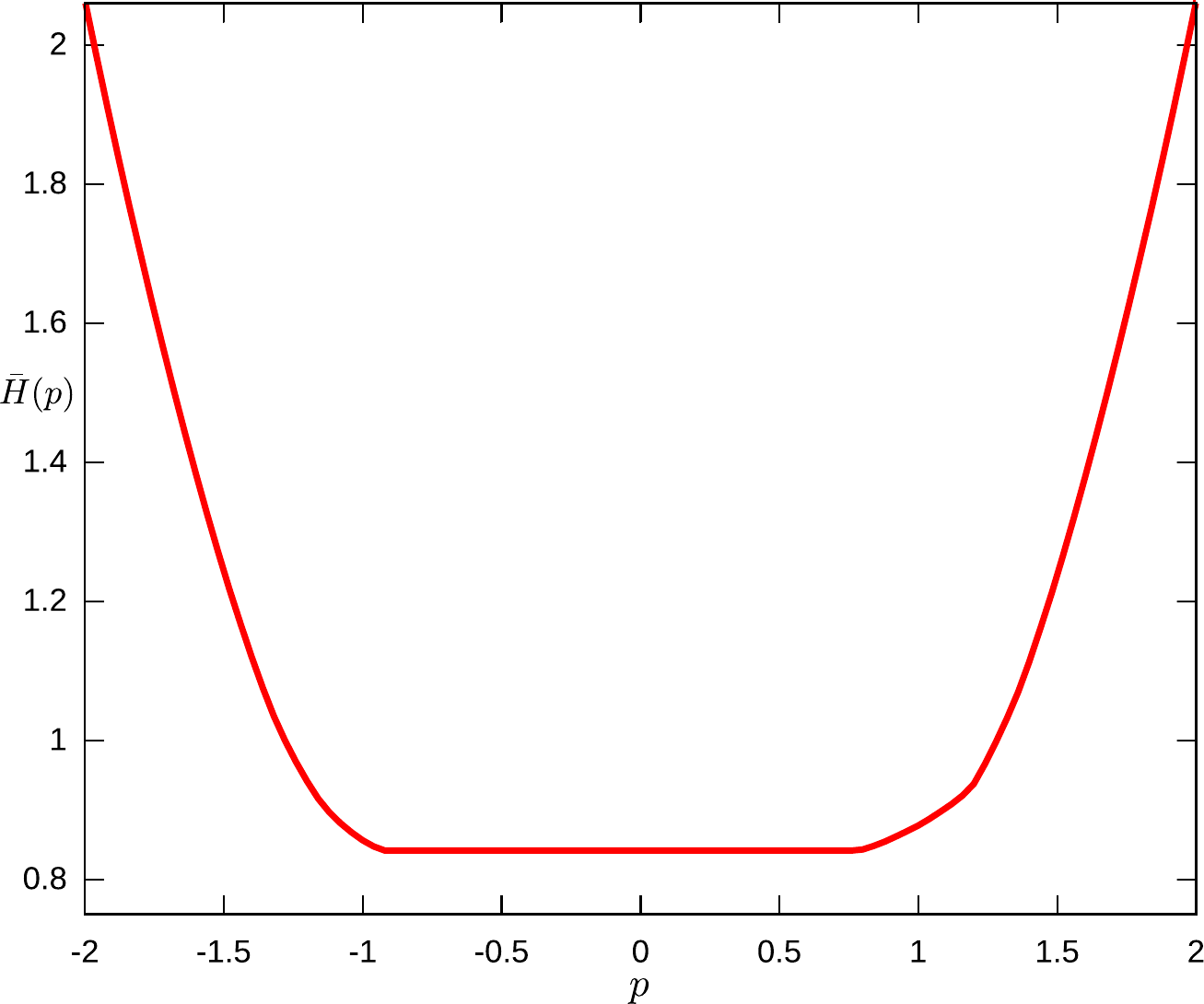} &
\includegraphics[width=.3\textwidth]{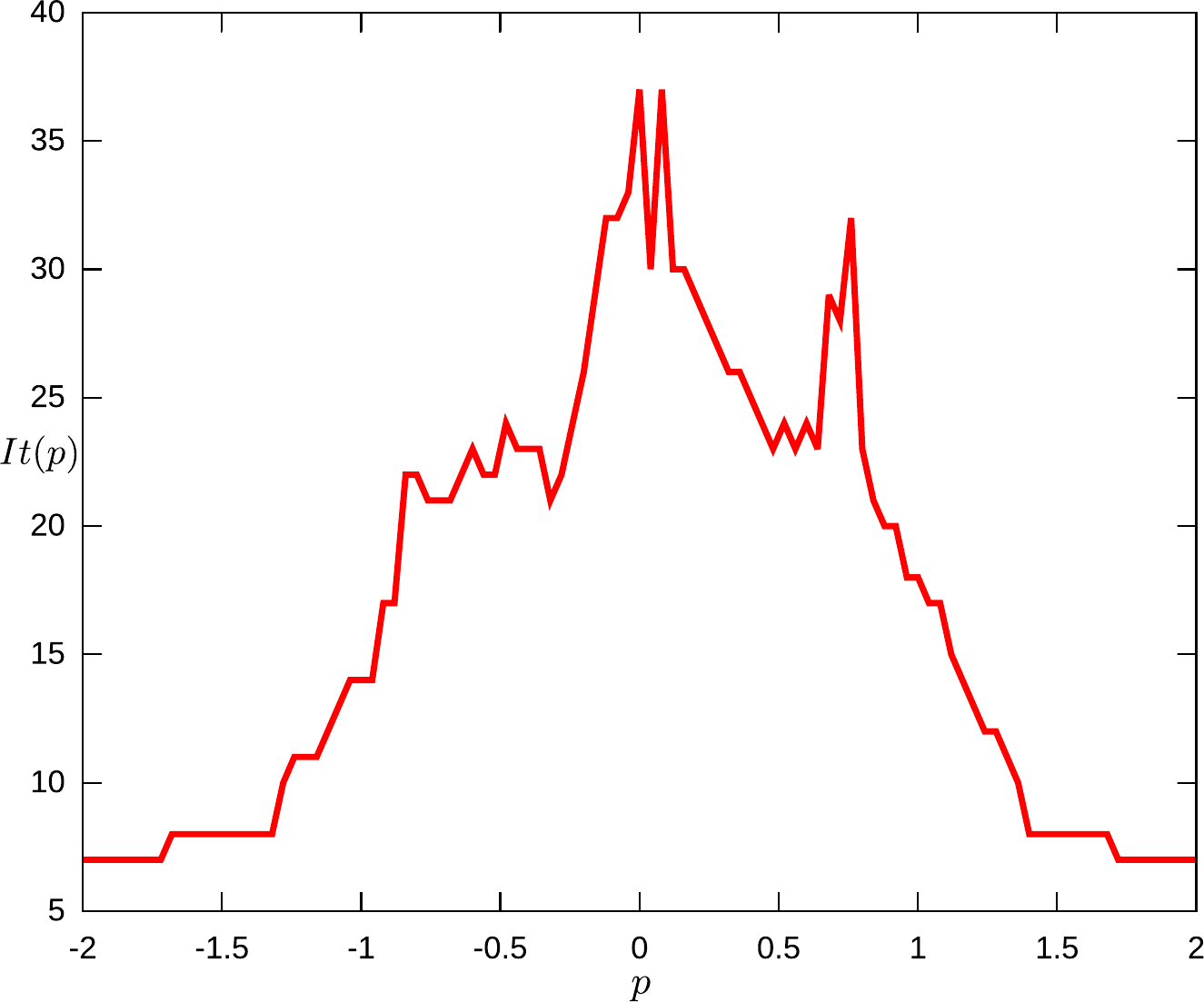} \\
(a)&(b)
\end{tabular}
\end{center}
\caption{Effective Hamiltonian for $p\in[-2,2]$ (a) and the corresponding iterations (b).}\label{Hbar-WCS-1d}
\end{figure}

\begin{figure}[h!]
 \begin{center}
 \begin{tabular}{ccccc}
\includegraphics[width=.165\textwidth]{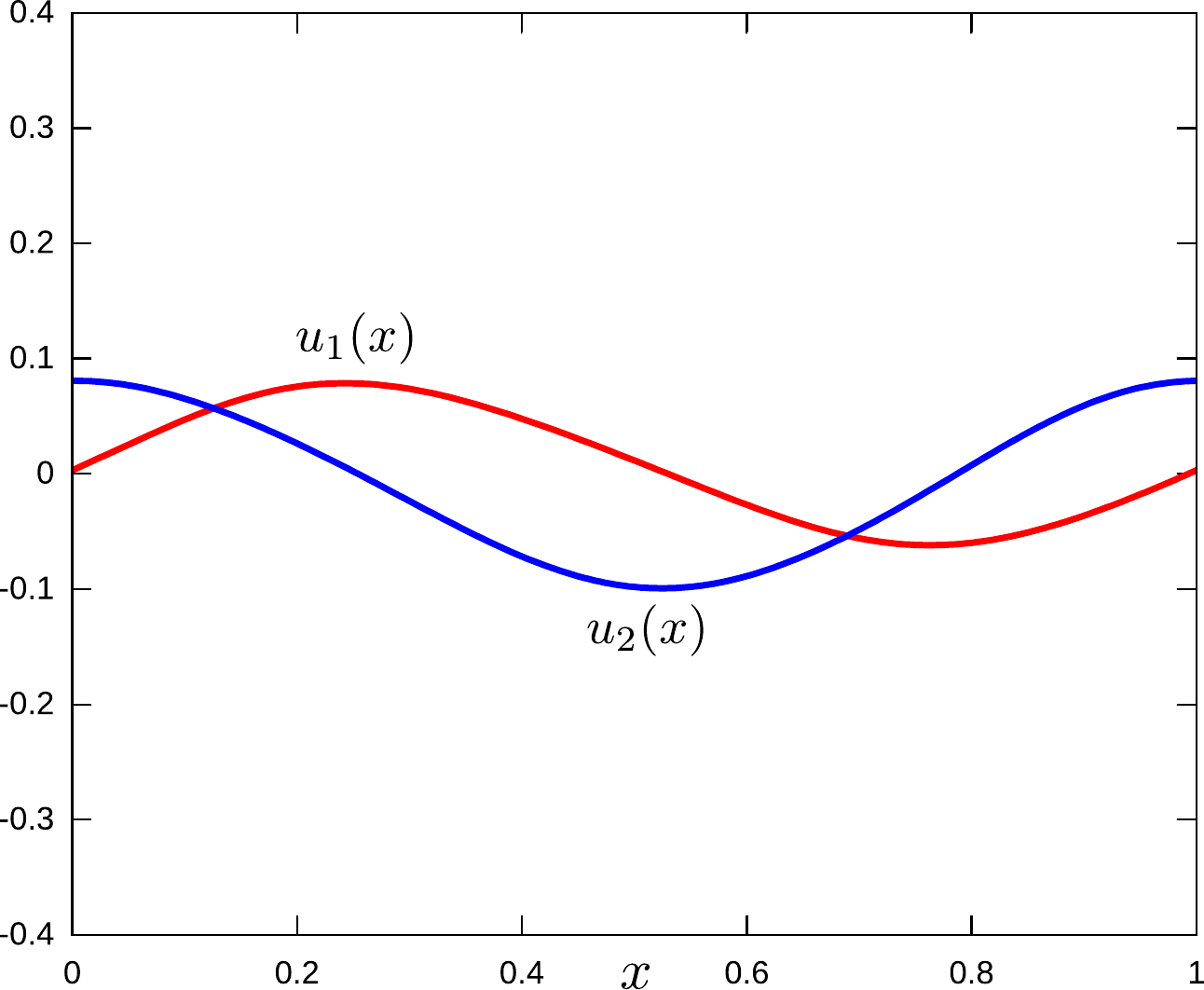} &
\includegraphics[width=.165\textwidth]{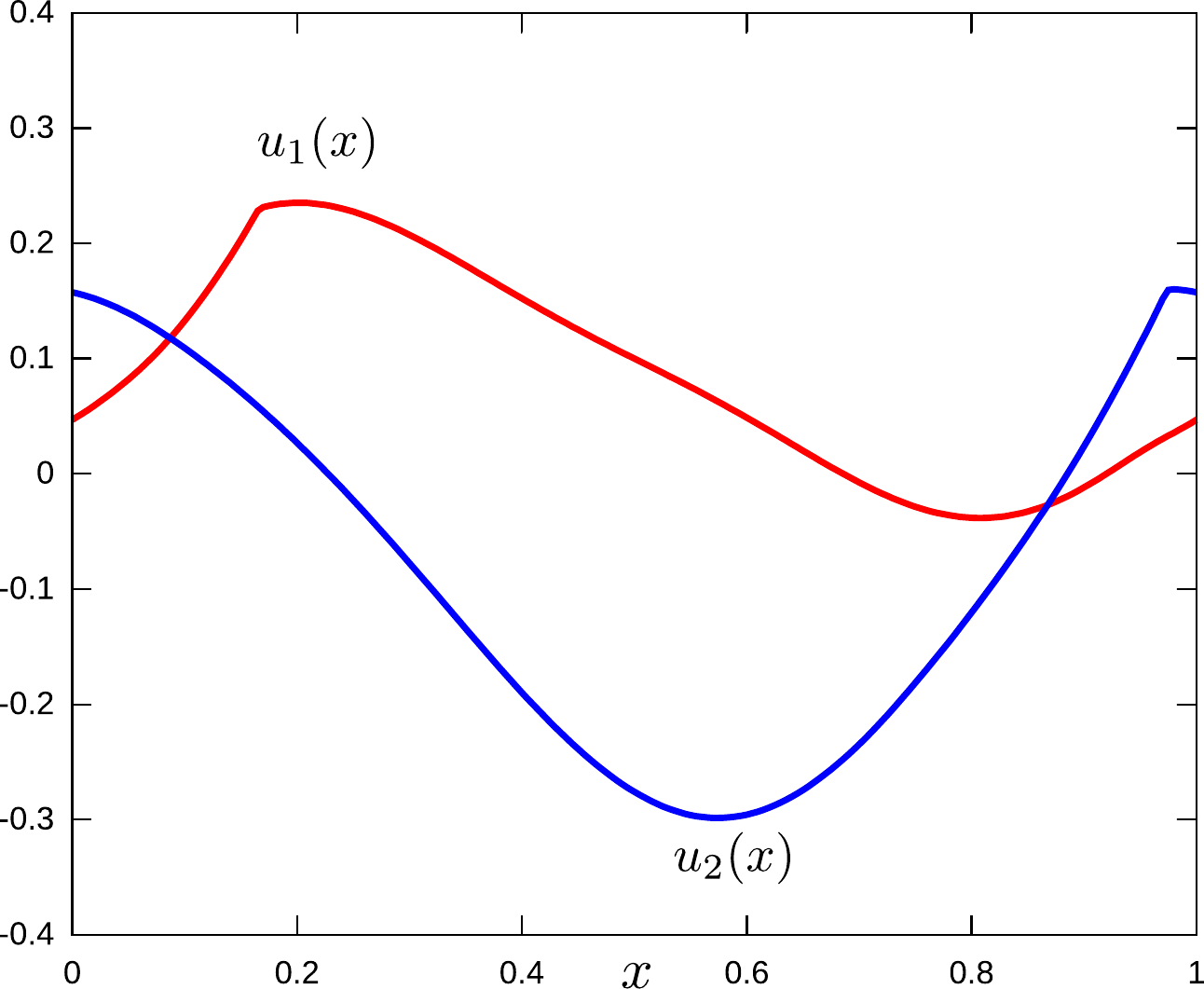} &
\includegraphics[width=.165\textwidth]{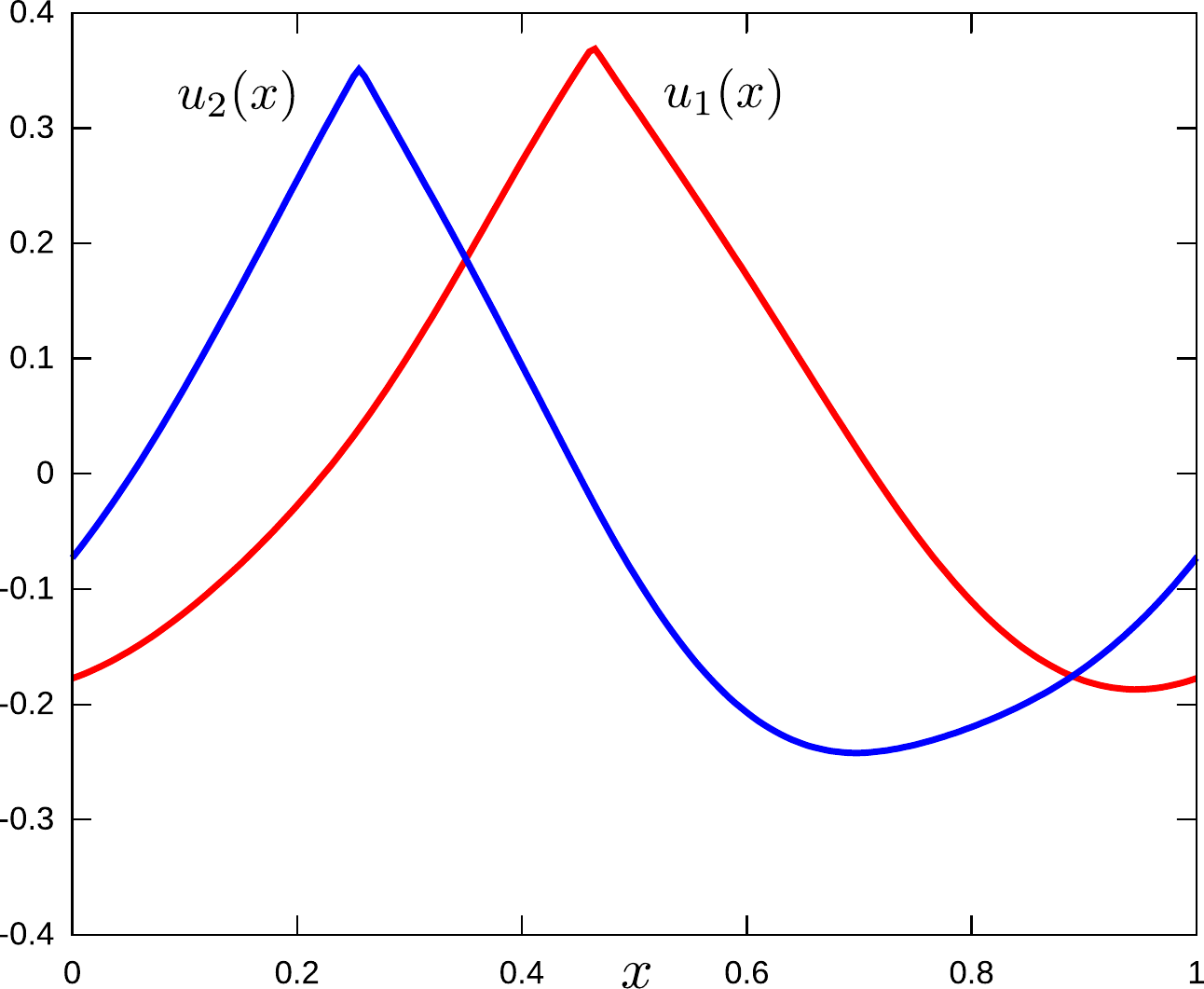} &
\includegraphics[width=.165\textwidth]{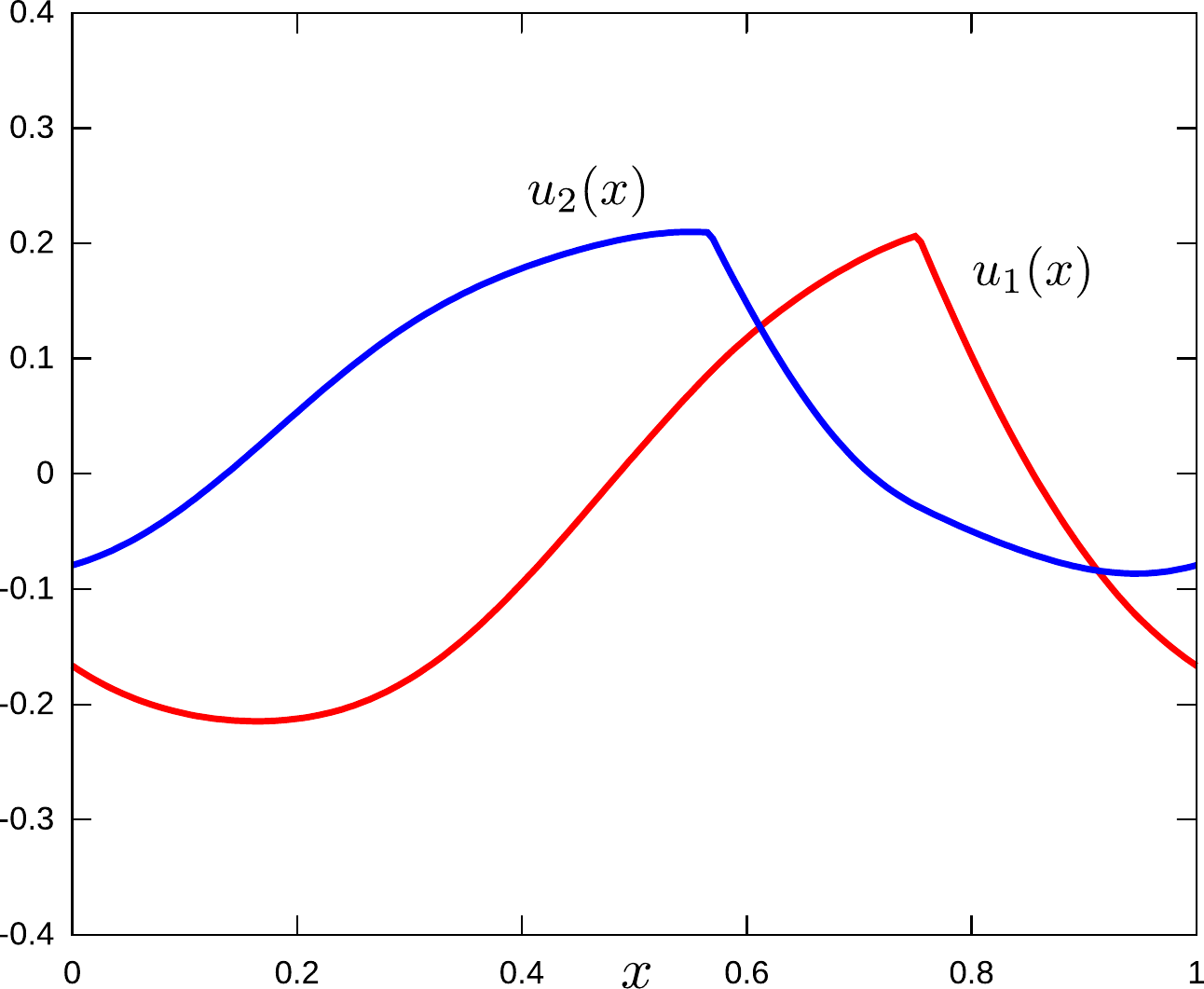} &
\includegraphics[width=.165\textwidth]{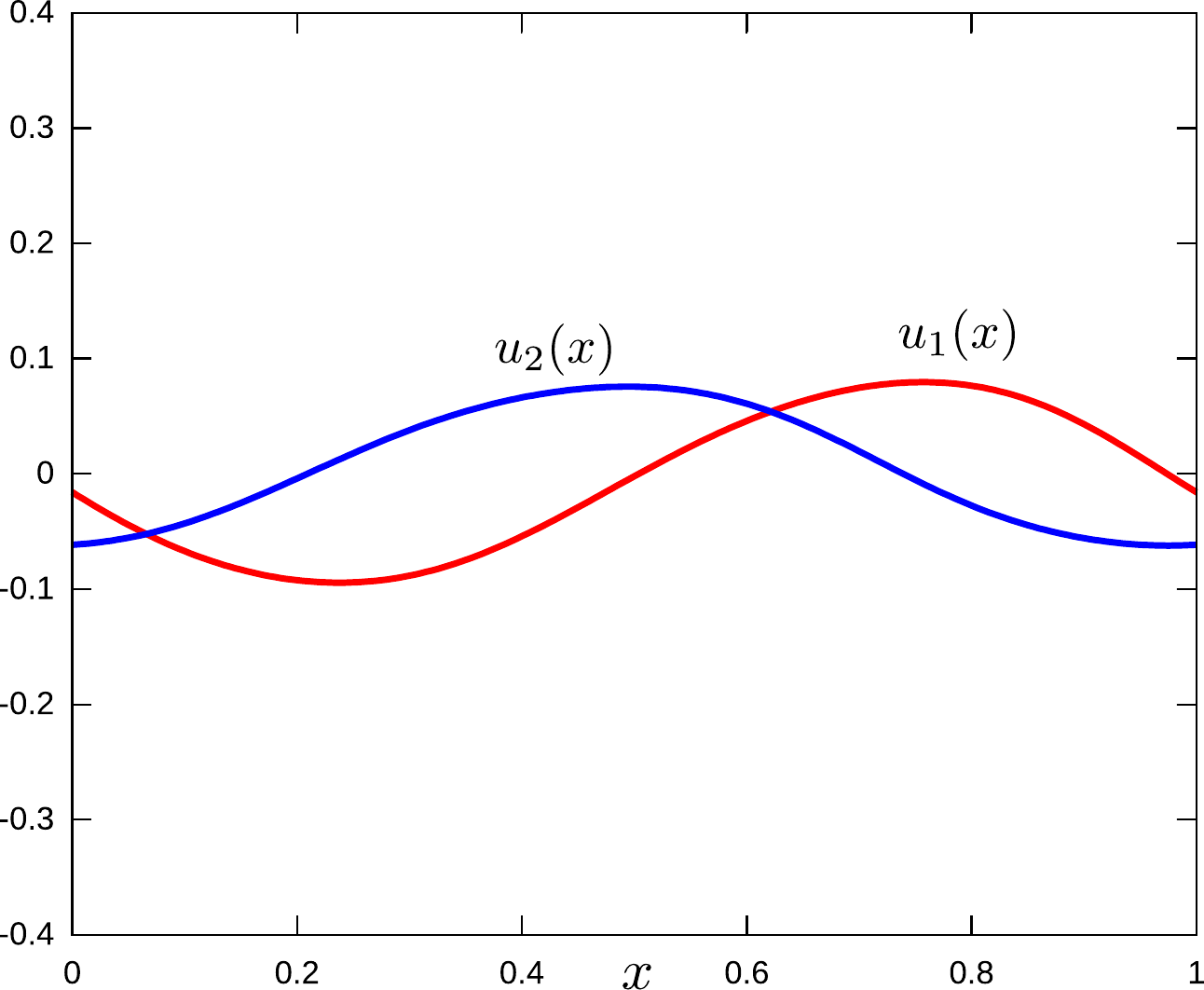} \\
(a)&(b)&(c)&(d)&(e)
\end{tabular}
\end{center}
\caption{Viscosity solution for $p=-2$ (a), $p=-1$ (b), $p=0$ (c), $p=1$ (d), $p=2$ (e).}\label{correctors-1d-wcs}
\end{figure}

In dimension $n=1$ we compute the effective Hamiltonian choosing $p\in[-2,2]$, discretized with $101$ nodes,
while the mesh size for the torus $\mathbb{T}^1$ is $100$ nodes. Moreover, we choose
$$V_1(x)=\sin(2\pi x)\,,\,\,\,\, V_2(x)=\cos(2\pi x)\,,\quad c_1(x)=1-\cos(4\pi x)\,,\quad c_2(x)=1+\sin(4\pi x)\,.$$
The total computational time is $3.19$ and the average number of iterations is $17$. Note that the dimension of the system ($200\times201$)
is doubled with
respect to the tests in the previous sections, since here the unknowns are $u_1$ and $u_2$ (plus $\lambda$).
In Figure \ref{Hbar-WCS-1d}a we show the graph of $\bar H$ and in Figure \ref{Hbar-WCS-1d}b
the corresponding number of iterations to reach convergence as a function of $p$.
\\
We readily observe some interesting features. First, the effective Hamiltonian is no longer symmetric with respect to $p$ (at least for small values of $p$) and
the computed plateau is $\{\bar H(p)=0.8417\}=[-0.925,0.788]$. This asymmetry is also evident looking at the graph of the number of iterations. On the other hand, we see that the number of iterations starts increasing consistently when $p$
enters the interval $[-1.29,1.36]$. As already remarked, this typically indicates that the corresponding corrector of the cell problem
starts developing kinks. Here it is quite surprising that this interval contains the plateau, differently from the scalar case where the two intervals coincide.
Hence we expect to find nonsmooth solutions also outside the
plateau. This is confirmed by Figure \ref{correctors-1d-wcs}, where we show pairs $(u_1,u_2)$ of correctors for different values of $p$.
These features deserve some theoretical investigation.\vspace{-5pt}
\\

Now we consider the same problem in dimension $n=2$. We compute the effective Hamiltonian choosing $p\in[-4,4]^2$, discretized with $51\times51$ nodes,
while the mesh size for the torus $\mathbb{T}^2$ is $25\times25$ nodes. Moreover, we choose
$$V_1(x_1,x_2)=\sin(2\pi x_1)\sin(2\pi x_2)\,,\quad V_2(x_1,x_2)=\cos(2\pi x_1)\cos(2\pi x_2)\,,$$
$$c_1(x_1,x_2)=1-\cos(4\pi x_1)\cos(4\pi x_2)\,,\quad c_2(x_1,x_2)=1+\sin(4\pi x_1)\sin(4\pi x_2)\,.$$
For a single point $p$, the average number of iterations is $12$ and  the average computational time  is $1.39$ seconds, while
the total computational time is $3640.38$.
Note that this CPU time is still reasonable, considering that the dimension of the system if four times the one in the scalar cases.

\begin{figure}[h!]
 \begin{center}
 \begin{tabular}{cc}
\includegraphics[width=.35\textwidth]{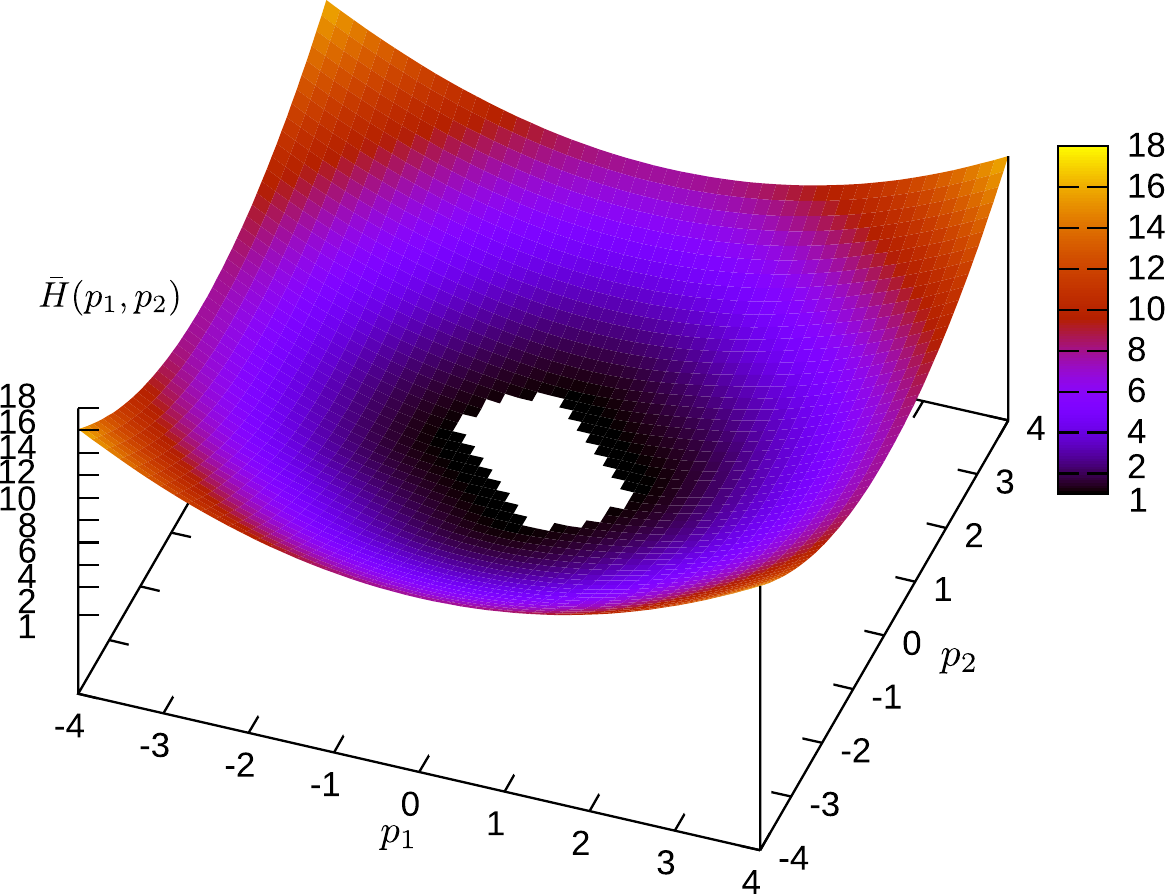} &
\includegraphics[width=.3\textwidth]{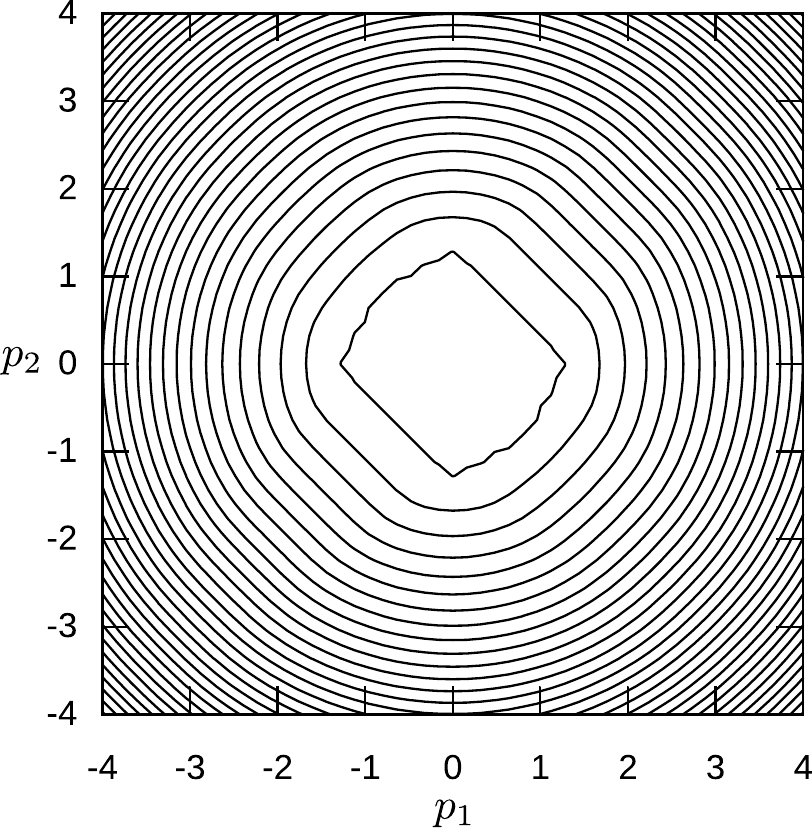} \\
(a)&(b)
\end{tabular}
\end{center}
\caption{Effective Hamiltonian for $p\in[-4,4]^2$: (a) surface and (b) level sets.}\label{Hbar-2d-WCS}
\end{figure}

\begin{figure}[h!]
 \begin{center}
 \begin{tabular}{cccc}
\includegraphics[width=.2\textwidth]{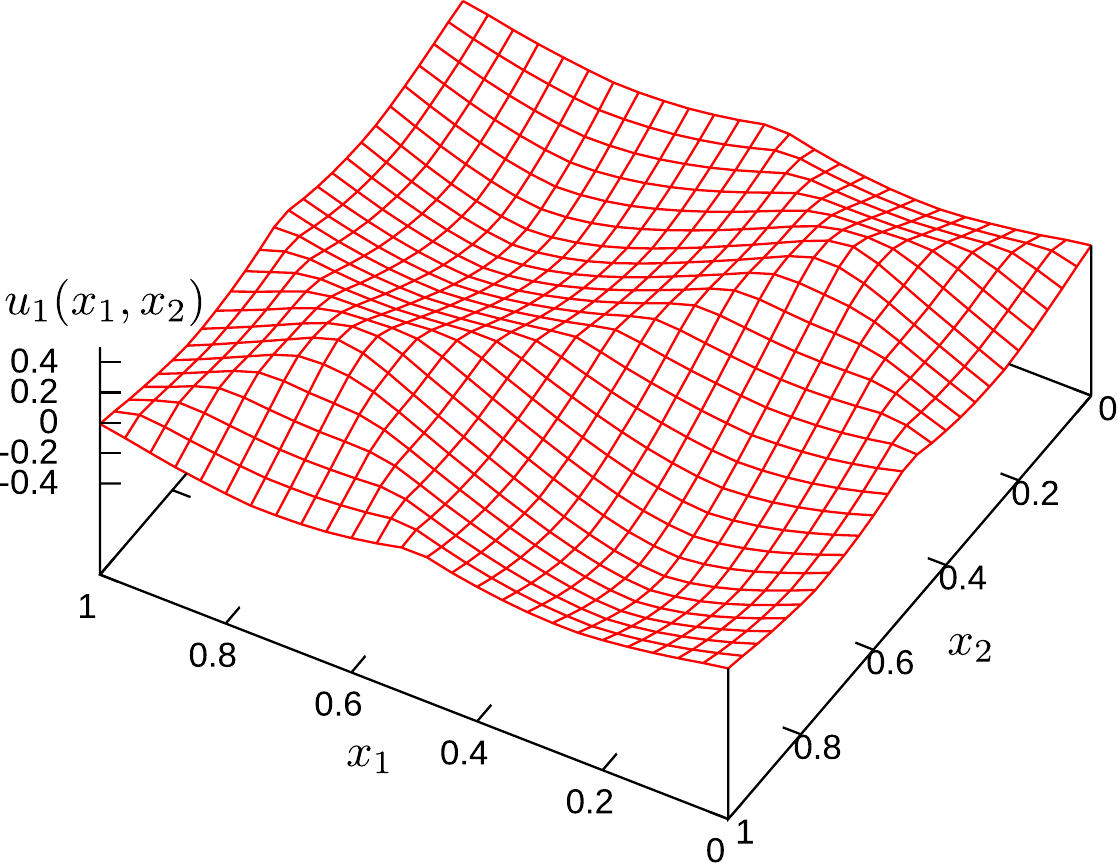} &
\includegraphics[width=.2\textwidth]{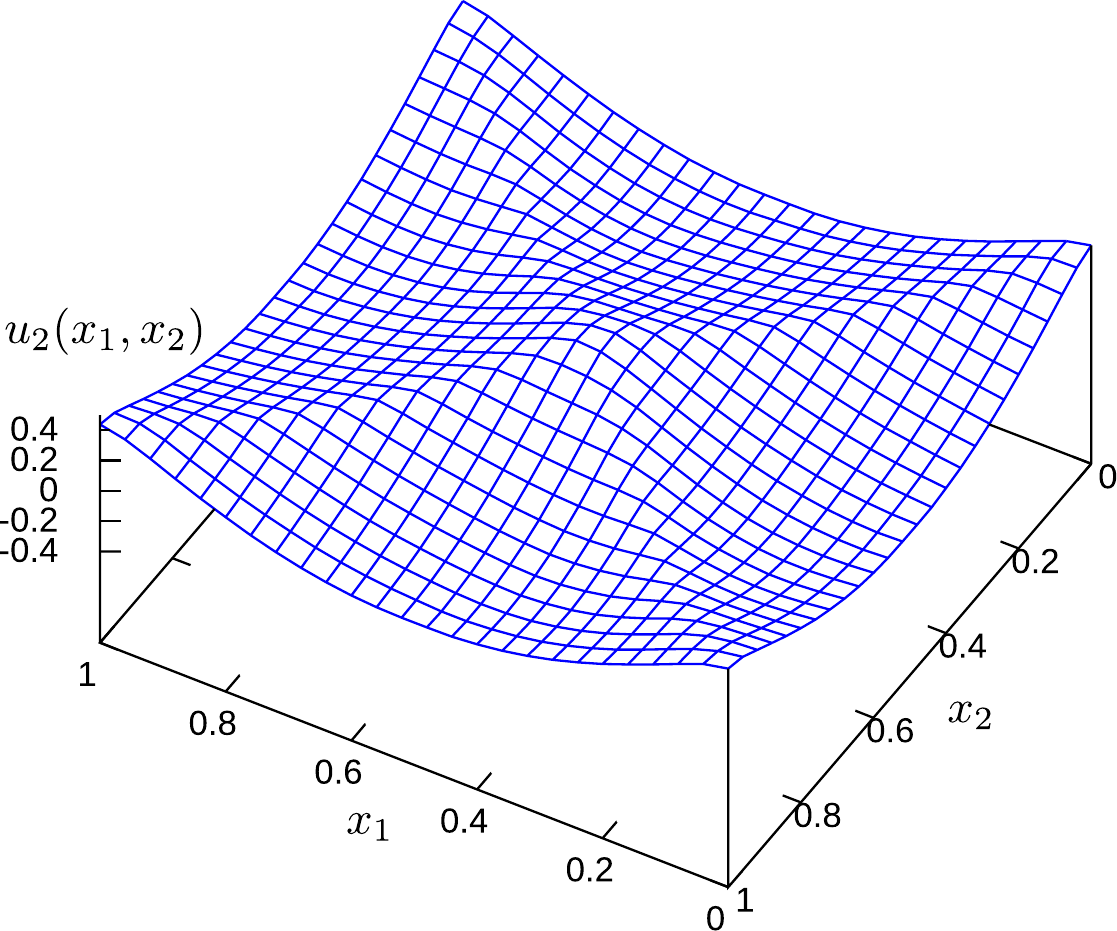} &
\includegraphics[width=.2\textwidth]{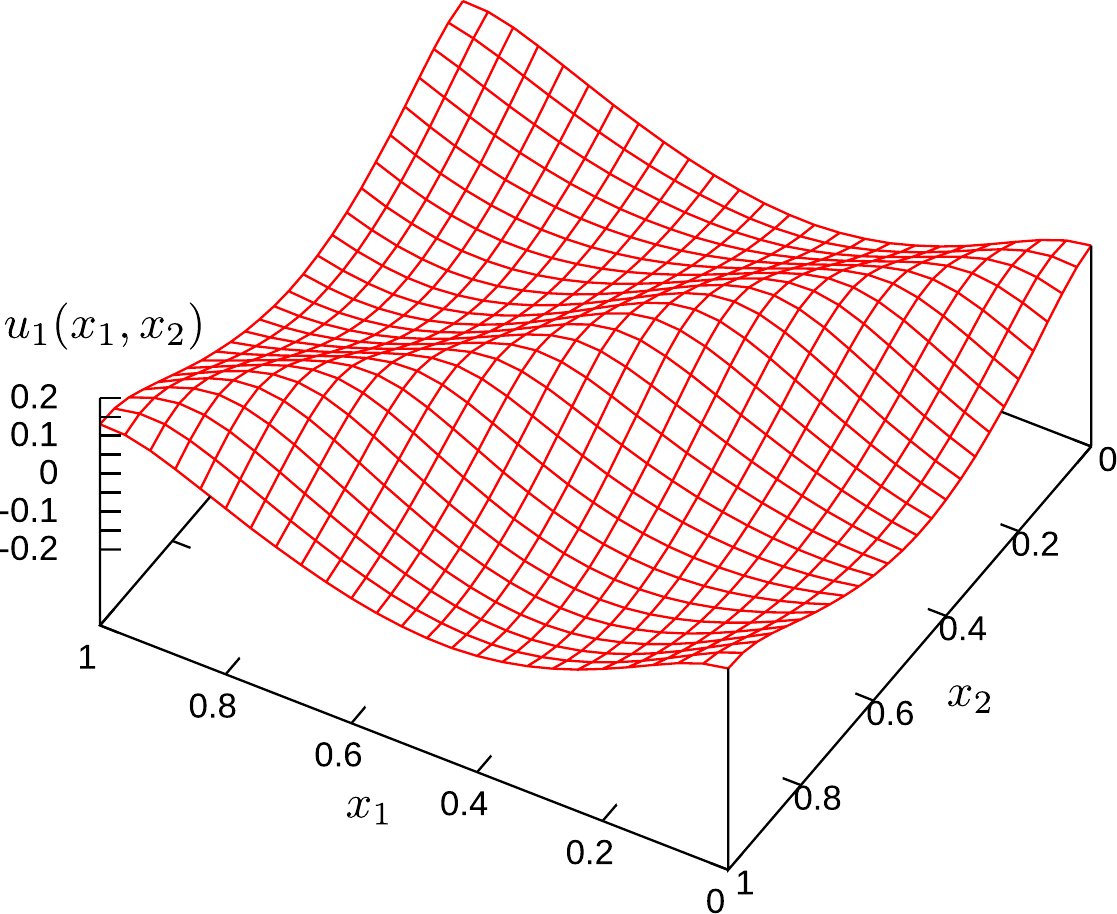} &
\includegraphics[width=.2\textwidth]{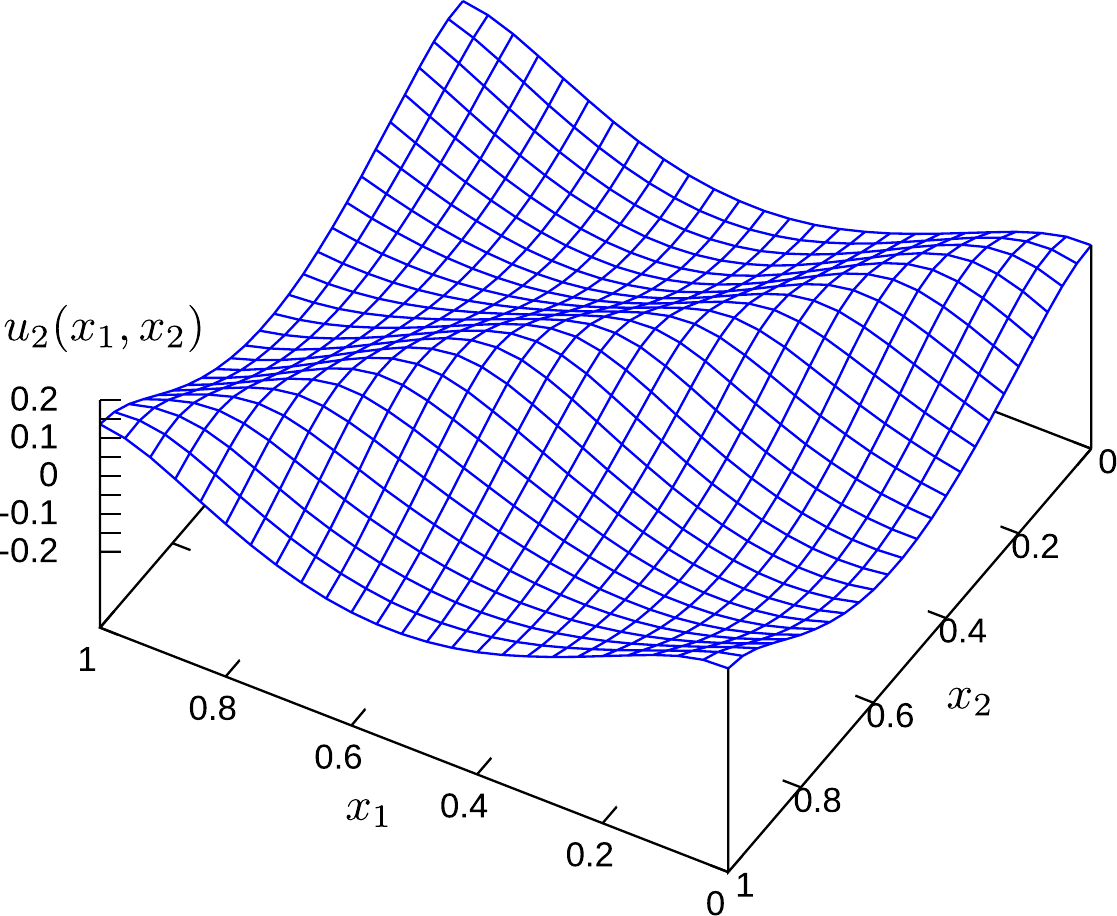} \\
(a)&(b)&(c)&(d)
\end{tabular}
\end{center}
\caption{Solutions of the cell problem for $p=(0,0)$ (a)-(b) and $p=(2,2)$ (c)-(d).}\label{correctors-2d-wcs}
\end{figure}

Figure \ref{Hbar-2d-WCS} shows the effective Hamiltonian surface and its level sets. The computed value inside the plateau is $\bar H(p)=1$.
\\
Finally, in Figure \ref{correctors-2d-wcs} we show two pairs of correctors for $p=(0,0)$ and $p=(2,2)$, respectively inside and outside the plateau.
\section{Dislocation dynamics}\label{dislocations}
Dislocations are line defects in the lattice structure of crystals, responsible for the plastic properties of the materials.
Several approaches have been proposed to study the effects of the interactions among these defects, including variational and PDE models.

In \cite{FIM}, a study on the homogenization of an ensemble of dislocations is presented,
in order to describe the effective macroscopic behavior of a system undergoing plastic deformations. This leads to consider a suitable cell problem
for a nonlocal Hamilton-Jacobi equation, that
we present here for simplicity in dimension one:\\ find $\l\in\R$ such that
\begin{equation}\label{CPdisl}
    \left(c_0(x)+L+M_p[u]\right) |Du+p |=\l,\qquad\qquad x\in\mathbb{T}^1
\end{equation}
admits a bounded and $1\hbox{-}$periodic viscosity solution $u$. The function $c_0$ is a $1$-periodic potential representing an obstacle to the motion of dislocations, 
$L\in\R$ is a drift representing a constant external stress and $M_p[u]$ is a nonlocal operator representing the interaction between dislocations, given by
$$
M_p[u](x)=\int_{\R}\mathcal{J}(z)\left\{E\left(u(x+z)-u(x)+pz\right)-pz\right\}dz\,,
$$
where $\mathcal{J}:\R\to\R^+$ is a nonnegative kernel satisfying
$$\mathcal{J}(-z)=\mathcal{J}(z)\quad \forall z\in\R\,,\qquad \mathcal{J}(z)\sim\frac{1}{z^2}\quad \mbox{for }|z|\gg 1$$
and $E:\R\to\R$ is an odd modification of the integer part, such that
$$E(\alpha)=\left\{\begin{array}{ll}
                    k & \mbox{if }\alpha=k\in\Z\,,\\
                    k+1/2 & \mbox{if } k<\alpha<k+1\,,\quad k\in\Z\,.
                   \end{array}
                   \right.
$$
In this setting the number $p\in\R$ represents the macroscopic density of dislocations. Moreover, the
dislocation lines are assumed to lie on a single plane, of which we only look at a cross section,
so that they are described as particle points corresponding to the integer level sets of $u(x)+px$.

\noindent The existence and uniqueness of $\lambda$ and the existence of a solution $u$ to the above cell problem is proved in \cite{FIM}, using a
suitable notion of viscosity solution for discontinuous HJ equations. Moreover, some
numerical approximations of \eqref{CPdisl} have been proposed in \cite{GHM} and in \cite{CCM}, both using the large-$t$ method.

Here we present the solution to this cell problem within our framework.
Following \cite{CCM},
we assume that only rational dislocation densities are  allowed, namely $p=P/Q$, with $P\in\Z$, $Q\in\mathbb{N}\setminus\{0\}$. Moreover, we fix
$N\in\mathbb{N}$ and choose the space step $h=1/N$. This allows to discretize the nonlocal term in \eqref{CPdisl} as
$$
M_p[u](x_i)=h\sum_{m=0}^{Q-1}\sum_{j=0}^{N-1}\left(J_{m,j}^Q E(u_{i+j}-u_i+p(x_j+m))- p\zeta_{m,j}^Q\right)\,,
$$
for $i=0,...,N-1$, with
$$
J_{m,j}^Q=\mathcal{J}^Q(x_j+m)\,,\qquad \zeta_{m,j}^Q=(x_j+m)\mathcal{J}^Q(x_j+m)\,,\qquad
\mathcal{J}^Q(t)=\hspace{-5pt}\sum_{k=-N_0}^{N_0}\hspace{-5pt}\mathcal{J}(t+kQ)\,,
$$
where $\mathcal{J}^Q$ is a finite sum approximation (up to a given integer $N_0\in\mathbb{N}$) of a $Q$-periodic version of the kernel $\mathcal{J}$.
Note that, for $m=0,...,Q-1$ and $j=0,...,N-1$, the $2NQ$ terms $J_{m,j}^Q$ and $\zeta_{m,j}^Q$
can be pre-computed.
Moreover, the integer part $E$ is approximated by a piecewise linear approximation $E_\tau$, such that each integer jump is replaced, in a
$\tau\hbox{-}$neighborhood,
by a ramp with slope $1/\tau$.
In this way the derivative $E^\prime$, needed for the Newton linearization, is replaced by $E_\tau^\prime$, which is in turn a
$\tau$-dependent piecewise constant
approximation of the  Dirac measure $\delta_\Z$.

We recall that the singularity in the derivative of the term $|Du+p|$ in \eqref{CPdisl} is handled by replacing it with zero at the points where it vanishes,
as for the $q\hbox{-}$power Hamiltonians in Section \ref{qnorms}.

Finally, we have to update the approximation of the gradient terms according to the sign of the nonlocal velocity
$c[u]:=\left(c_0(x)+L+M_p[u]\right)$ in \eqref{CPdisl}, in order to correctly select the viscosity solutions.
To be precise, if the Engquist-Osher approximation is employed, we set
$$|Du|(x_i)\sim\left\{\begin{array}{ll}
 \sqrt{{\min}^2\left\{\displaystyle\frac{U_{i+1}-U_i}{h},0\right\} +{\max}^2\left\{\displaystyle\frac{U_i-U_{i-1}}{h},0\right\}} & \mbox{if } c[U](x_i)\geq 0\,,\\\\
 \sqrt{{\max}^2\left\{\displaystyle\frac{U_{i+1}-U_i}{h},0\right\} +{\min}^2\left\{\displaystyle\frac{U_i-U_{i-1}}{h},0\right\}} & \mbox{if } c[U](x_i)< 0\,.\\
\end{array}\right.
$$

We discretize the torus $\mathbb{T}^1$ with $100$ nodes, so that the space step is $h=0.01$. Choosing $Q=10$, the allowed densities are given by $p=P/Q=0.1P$ for
$P\in\Z$. We take $(p,L)\in[-4,4]^2$ discretized with $81\times 81$ nodes and we choose $c_0(x)=2\sin(2\pi x)$.

\begin{figure}[h!]
 \begin{center}
 \begin{tabular}{cc}
\includegraphics[width=.45\textwidth]{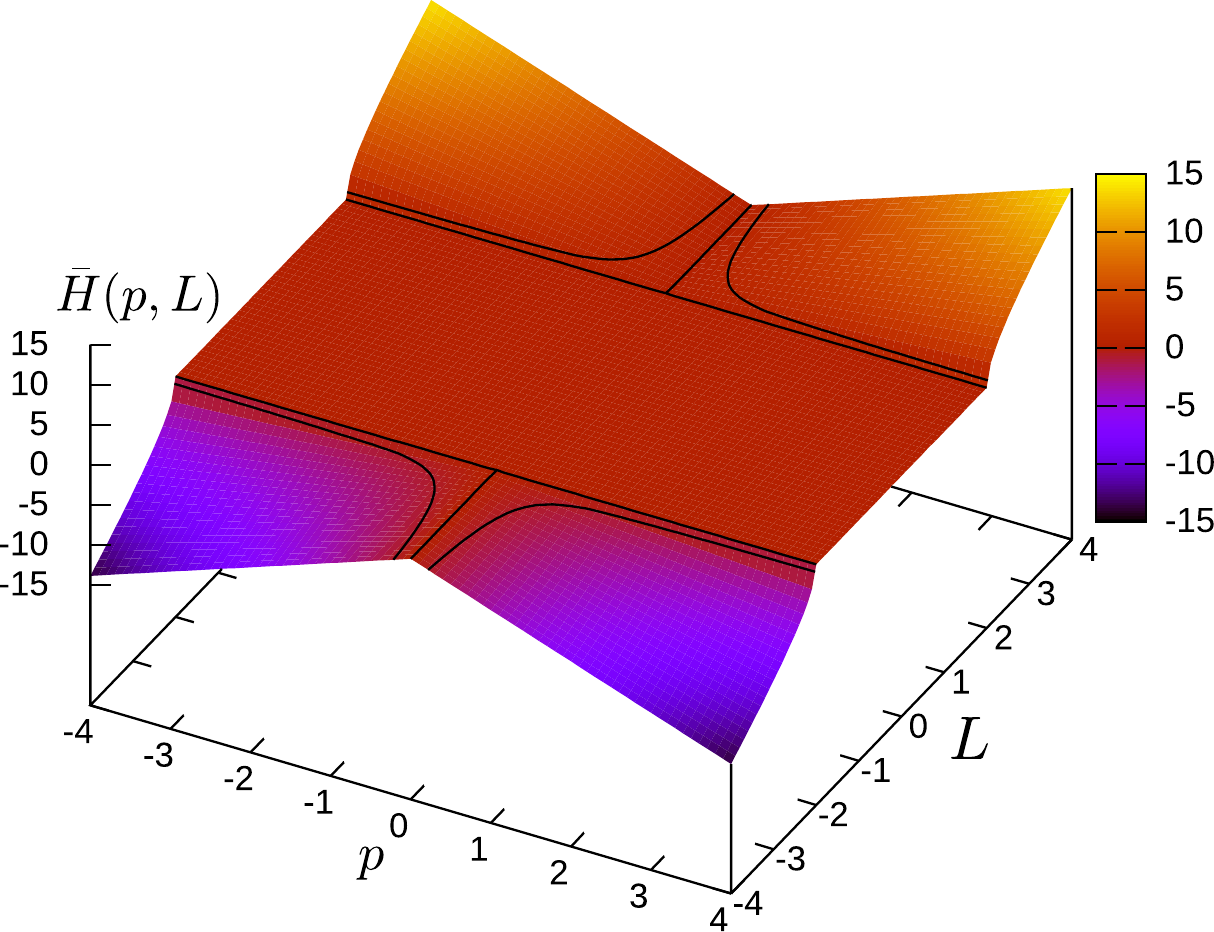} &
\includegraphics[width=.3\textwidth]{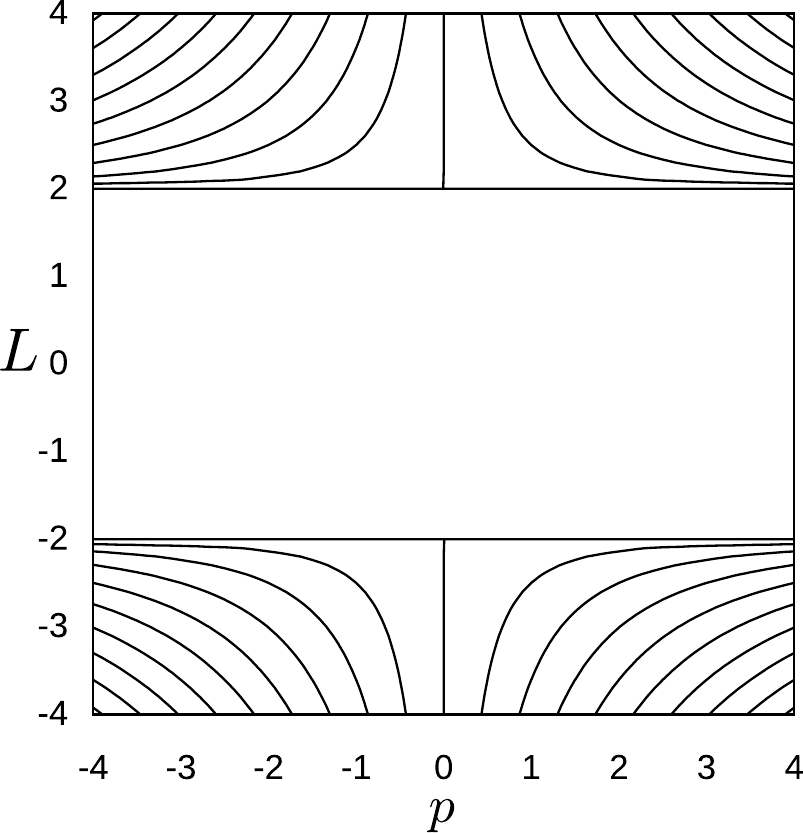} \\
(a)&(b)
\end{tabular}
\end{center}
\caption{Effective Hamiltonian for $(p,L)\in[-4,4]^2$: (a) surface and (b) level sets.}\label{hbar1d-local}
\end{figure}

We perform a preliminary test in which the dislocations do not interact, namely we completely remove the nonlocal term $M_p[u]$ in \eqref{CPdisl}
by setting $\mathcal{J}\equiv 0$. In this situation, independently of  $p$, dislocation particles are not able to get out the wells of the potential
$c_0$, unless the drift $L$ is strong enough. In the present example the critical stress is just the amplitude of the potential ($L=2$)
and we expect {\em no effective motion}, i.e., $\bar H(p,L)=0$ for every $(p,L)$ in the strip $\R\times[-2,2]$. Moreover, for $p=0$ we expect $\bar H(0,L)=0$ for all $L$ 
(i.e., no particles no motion!).
This is confirmed by Figure \ref{hbar1d-local}, in which we show the surface and the level sets of the computed effective Hamiltonian.

\begin{figure}[h!]
 \begin{center}
 \begin{tabular}{cc}
\includegraphics[width=.45\textwidth]{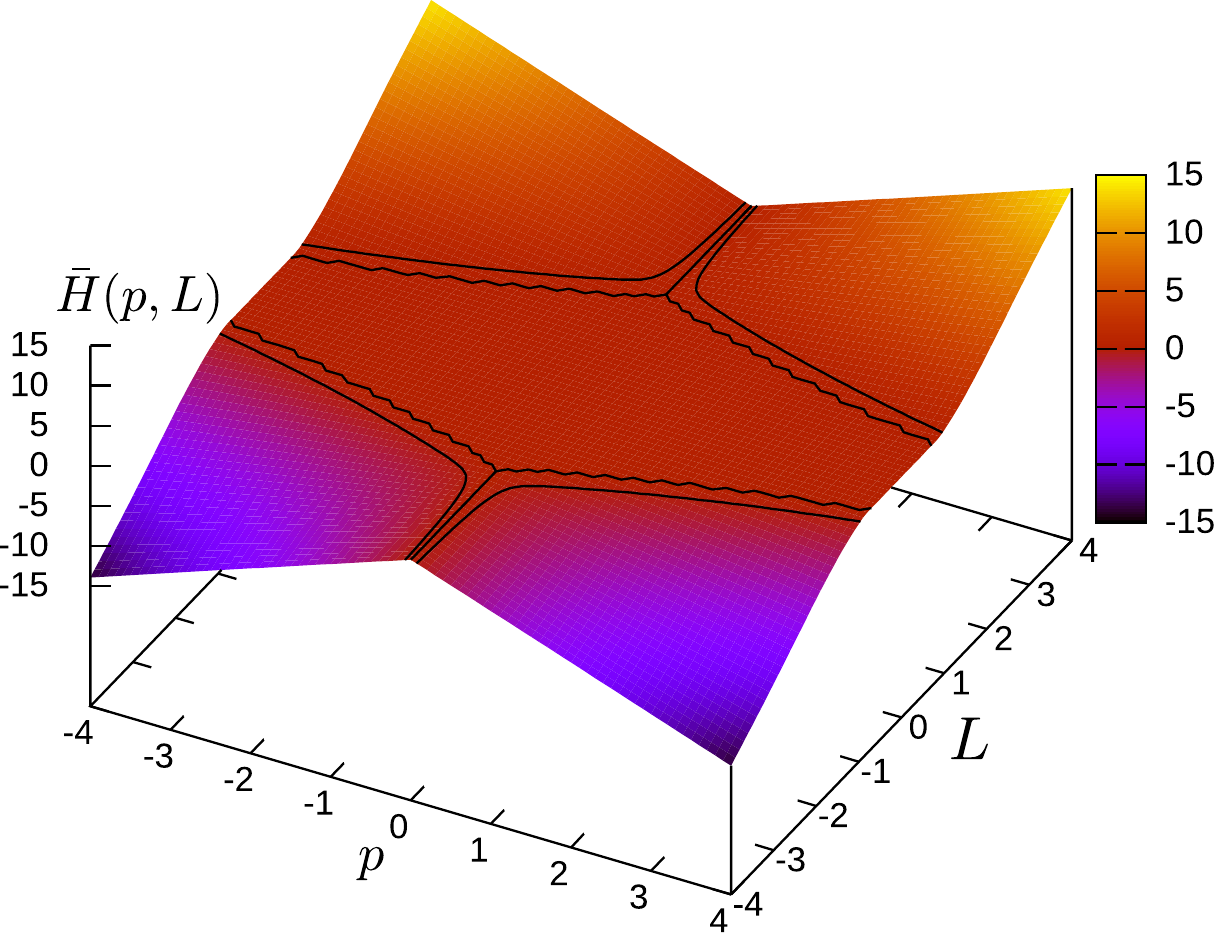} &
\includegraphics[width=.3\textwidth]{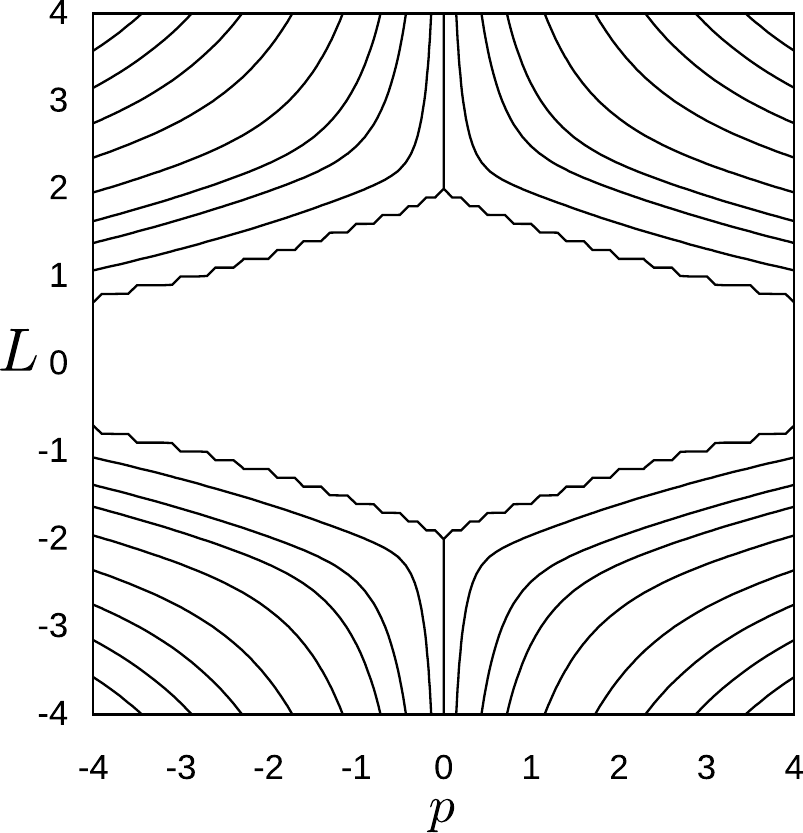} \\
(a)&(b)
\end{tabular}
\end{center}
\caption{Effective Hamiltonian for $(p,L)\in[-4,4]^2$: (a) surface and (b) level sets.}\label{hbar1d-nonlocal}
\end{figure}

We proceed with an intermediate test, studied in \cite{GHM} as a regularization of equation \eqref{CPdisl}, corresponding
in our setting to the choice
$E(\alpha)=\alpha$ (no jumps) with a kernel $\mathcal{J}$ satisfying $\int_\R \mathcal{J}(z)dz=1$ and such that its $Q\hbox{-}$periodic version is just
$\mathcal{J}^Q(t)\equiv 1/Q$. Accordingly, the nonlocal operator $M_p[u]$ reduces to a convolution of $u$ with $\mathcal{J}-\delta_0$.
In this situation, we expect the nonlocal interactions to produce an additional energy that can be sufficient to move the particles, i.e., 
$\bar H>0$ despite the external stress
is less than critical ($L<2$). Moreover, this effect should be amplified as the particle density $p$ increases (more particles more interactions). This is what is observed in Figure \ref{hbar1d-nonlocal},
in which we show the surface and the level sets of the computed effective Hamiltonian. \\
\\
Note that the boundary of the plateau of $\bar H$ is not smooth and
this reflects the fact that just adding a single particle in the wells of the potential can severely affect the collective motion.
Our result is in qualitative agreement with that found in \cite{GHM}.

\begin{figure}[h!]
 \begin{center}
 \begin{tabular}{cc}
\includegraphics[width=.45\textwidth]{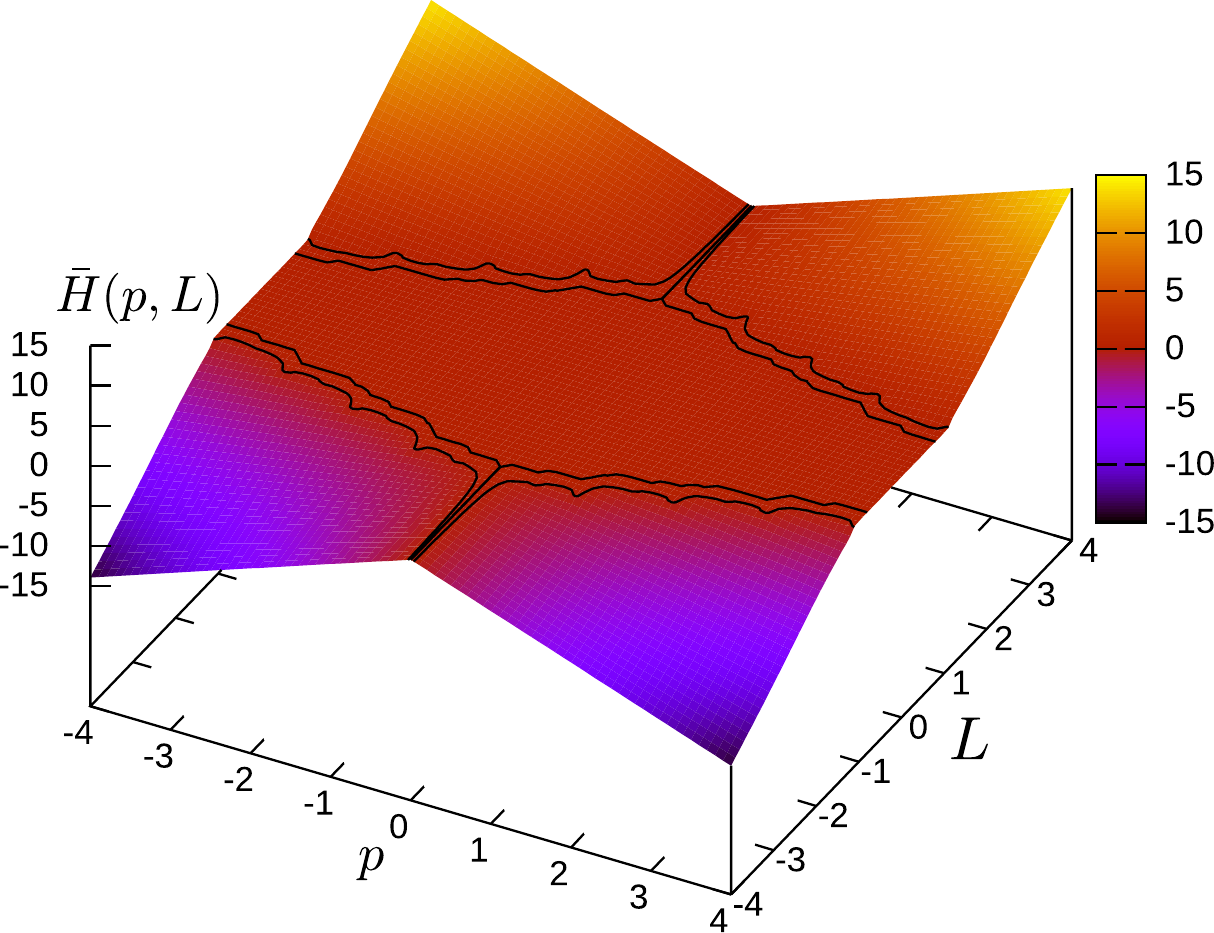} &
\includegraphics[width=.3\textwidth]{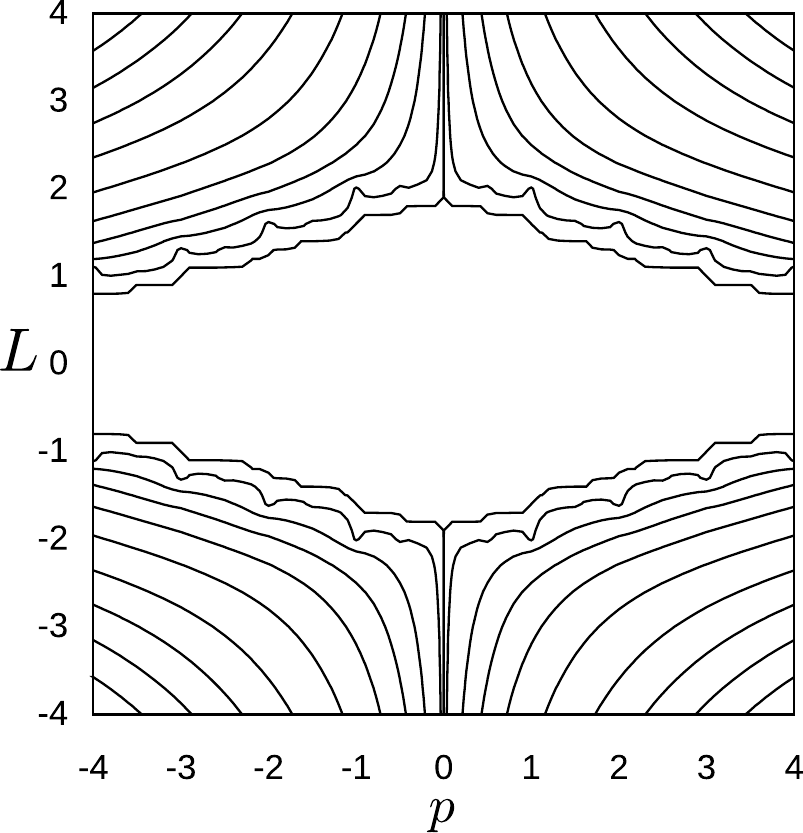} \\
(a)&(b)
\end{tabular}
\end{center}
\caption{Effective Hamiltonian for $(p,L)\in[-4,4]^2$: (a) surface and (b) level sets.}\label{hbar1d-dislocations}
\end{figure}

\begin{figure}[h!]
 \begin{center}
 \begin{tabular}{c}
\includegraphics[width=.5\textwidth]{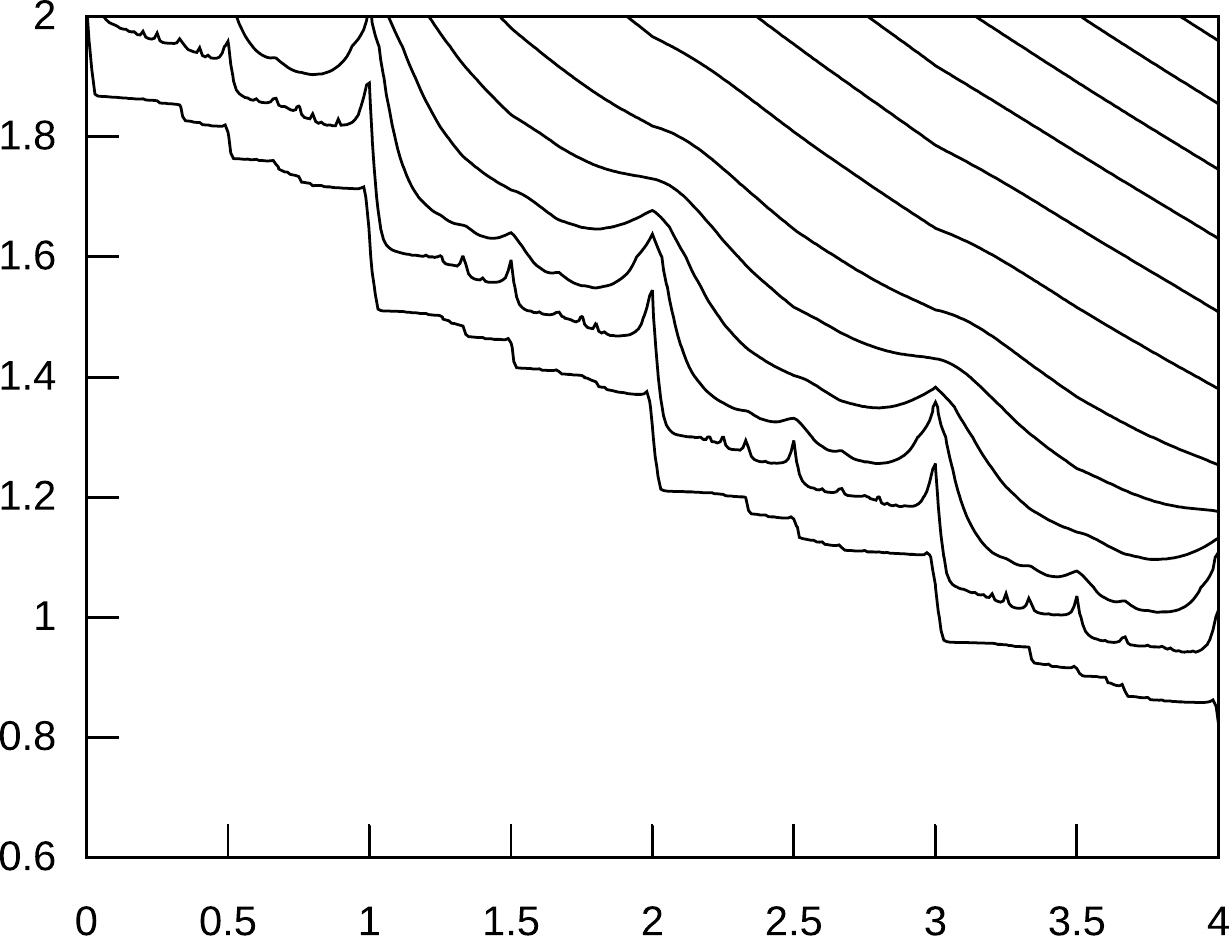}
\end{tabular}
\end{center}
\caption{Effective Hamiltonian for $(p,L)\in[0,4]\times [0.6,2]$ in the case $Q=100$.}\label{hbar1d-dislocations-zoom}
\end{figure}

Finally, let us consider the general case, studied in \cite{CCM}. We choose the kernel
$$\mathcal{J}(z)=\min\left\{C_\mathcal{J},\frac{1}{z^2}\right\}\,,$$
for a positive constant $C_\mathcal{J}$ such that $\int_\R \mathcal{J}(z)dz=1$. Moreover, we choose $N_0=100$ for the approximation of the
$Q\hbox{-}$periodic kernel $\mathcal{J}^Q$ and $\tau=5h$ for the approximation of the integer part $E$ and its derivative.
In Figure \ref{hbar1d-dislocations} we show the surface and the level sets of the computed effective Hamiltonian.
\\
The expected behavior is similar to the one described in the previous test, but here the result is quite surprising
as in the analogous simulation obtained in \cite{CCM}. We clearly see that level sets close to zero are no longer monotone,
some spikes appear at the integer densities $p\in\Z$, while smaller spikes appear at the half-integer densities $p\in\frac12\Z$.
This could sound weird from a physical point of view, since it is against the intuition that more particles reduce the critical stress needed to activate motion.
On the contrary this phenomenon reproduces what is known in the literature as {\em strain hardening}, namely the fact that a too high density of dislocations
translates in a mutual obstruction to the motion, and hence in a less effective cooperation. Nevertheless, it is still unclear to us why this hardening
occurs at so specific ratios between $P$ and $Q$. It may be related to the shape of the potential $c_0$ and it is still under investigation.
Moreover, the feeling is that the same behavior can occur at smaller scales with a kind of self-similar structure, depending on the choice of $Q$.
This is somehow supported by the simulation in Figure \ref{hbar1d-dislocations-zoom}, where we show the effective Hamiltonian
for $(p,L)\in[0,4]\times[0.6,2]$, computed using a grid of size $401\times 28$, with $N_0=1000$ for the $Q\hbox{-}$periodic
kernel approximation and $Q=100$, so that the allowed
densities are given by $p=P/Q=0.01P$ for $P\in\Z$.  We see that yet smaller spikes appear at intermediate {\em frequencies} between $\Z$ and
$\frac12\Z$.
\\
We conclude summarizing in Table \ref{table-1d-dislocations} the performance of the proposed method for the three main tests above, respectively identified by letters
($A$), ($B$) and ($C$).
\begin{table}[!ht]
 \centering
  \begin{tabular}{|c|c|c|c|}
    \hline
    Test & Av. CPU (secs) per $(p,L)$ & Av. Iterations per $(p,L)$ & Total CPU (secs)\\
    \hline
    ($A$)&0.017 &5 &115.71\\\hline
    ($B$)&0.032 & 7 &215.74\\\hline
    ($C$)&0.115 &16 &759.31\\\hline
 \end{tabular}
  \caption{Performance for the dislocation dynamics case.}\label{table-1d-dislocations}
\end{table}\\
Note the different performance between the local case ($A$) and the nonlocal cases, due to the fact that the Jacobian matrix in the
Newton linarization of the system is no longer sparse for ($B$) and ($C$). An additional and substantial computational cost appears in ($C$),
due to the computation of the nonlocal terms, depending on the interaction kernel  and the magnitude of $Q$.
Indeed, also the case $Q=100$ described above and not reported in the table is a very intensive task, it toke $10612$ seconds (about $3$ hours)
with an average CPU time per point $(p,L)$ of $1.02$ seconds, about ten times the case $Q=10$.

\section{Stationary mean field games}\label{meanfieldgames}
We consider the following class of ergodic Mean Field Games:
\begin{equation}\label{MFG}
 \left\{
  \begin{array}{ll}
   -\nu \Delta  u  +H(x, D u)+\l  =V[m] \qquad & x\in\T^n\\[4pt]
    \nu \Delta m +\mathrm{div}(m\,H_p(x, D u))=0 & x\in\T^n\\[4pt]
    \int_{\T^n} u(x)dx=0, \int_{\T^n} m(x)dx=1,\,m\ge 0\,.
  \end{array}
 \right.
\end{equation}
If $\nu>0$,  $H$ is smooth and convex,  then there exists a unique triple  $(u,m,\l)$ which is a classical solution of \eqref{MFG} (see \cite{LL}).
A finite difference scheme for \eqref{MFG} is presented in \cite{AC},
where it is proved the well-posedness of the corresponding discrete system and some convergence result. The discrete solution $(U,M,\Lambda)$ is computed by a large-$t$
approximation for both equations in \eqref{MFG}.

Here we present a direct resolution within our framework, in the simple case of the eikonal Hamiltonian in dimension two, with a cost function $f$ and a
local potential $V$, namely
$$
 \left\{
  \begin{array}{ll}
   -\nu \Delta  u  +|Du|^2+f(x)+\l  =V(m) \qquad & x\in\T^2\\[4pt]
    \nu \Delta m +2\,\mathrm{div}(m\,D u)=0 & x\in\T^2\\[4pt]
    \int_{\T^2} u(x)dx=0, \int_{\T^2} m(x)dx=1,\,m\ge 0\,.
  \end{array}
 \right.
$$
Differently from the previous sections, this problem is {\em overdetermined}. Indeed, introducing on the torus $\T^2$ a uniform grid with $N\times N$ nodes,
we end up with $2N^2+2$ nonlinear equations, corresponding to the discretization of the two PDEs and the normalization conditions, given by
$$h^2\sum_{i=0}^{N-1}U_i=0\,,\qquad\qquad h^2\sum_{i=0}^{N-1}M_i-1=0\,,$$
where $h=\frac{1}{N}$ is the space step. On the other hand, the number of unknowns is $2N^2+1$, corresponding to the degrees of freedom of $U$ and $M$
plus the additional unknown $\Lambda$. Note that we do not include the constraint $M\ge 0$ in the discretization.
Our experiments show that the normalization condition
for $M$ seems enough to force numerically its nonnegativity. This point is still under investigation.

For the sake of comparison, we present some of the tests reported in \cite{AC}.
We choose $U\equiv 0$, $M\equiv1$ and $\Lambda=0$ as initial guess and we set to $\varepsilon=10^{-6}$ the tolerance for the stopping criterion of the algorithm.
We first consider the case $V(m(x))=m^2(x)$,
$f(x)=\sin(2\pi x_1)+\cos(4\pi x_1)+\sin(2\pi x_2)$ and $\nu=1$.
We choose $N=50$ so that the size of the system is $5002\times 5001$.
\begin{figure}[h!]
 \begin{center}
 \begin{tabular}{ccc}
\includegraphics[width=.22\textwidth]{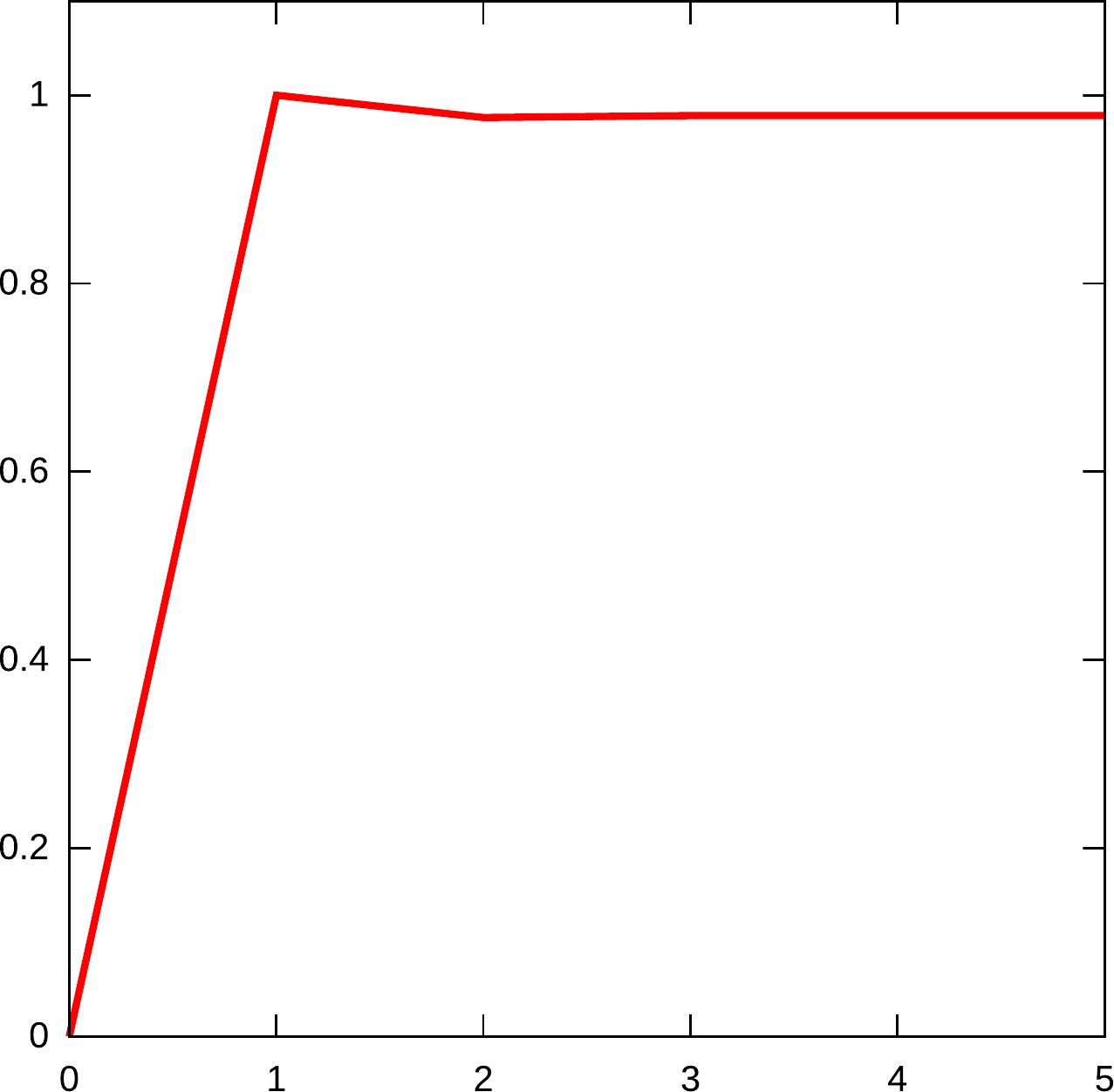} &
\includegraphics[width=.22\textwidth]{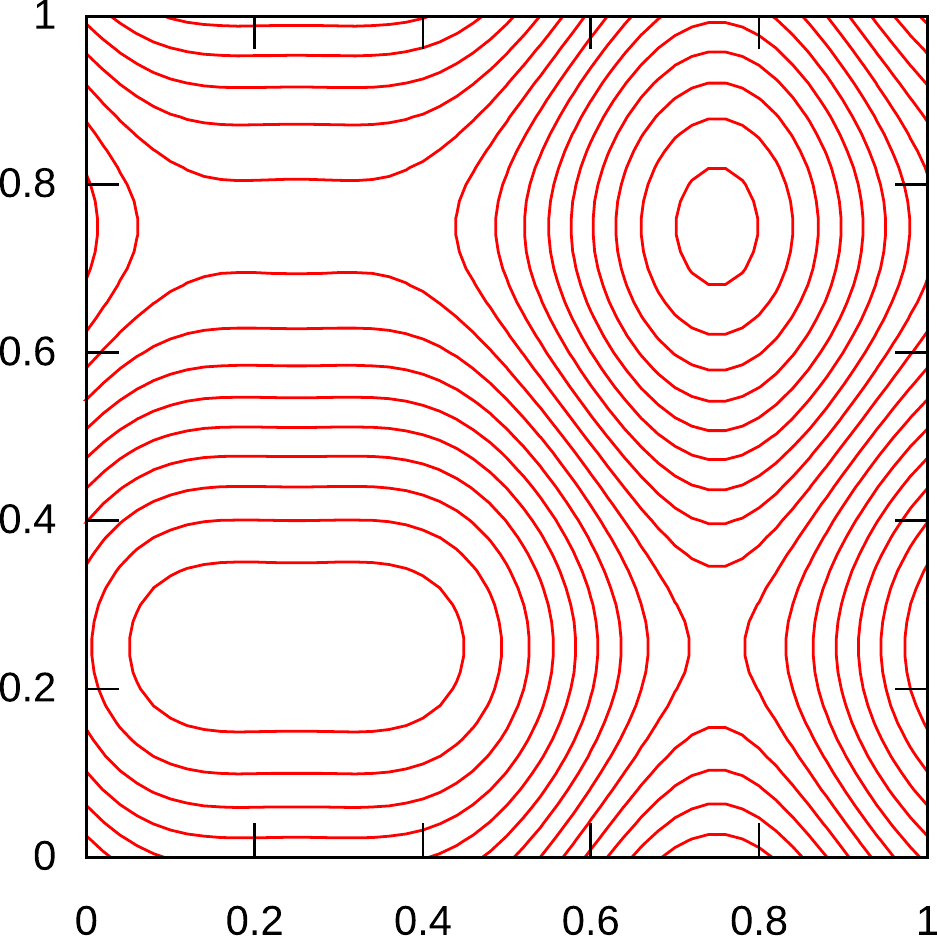} &
\includegraphics[width=.22\textwidth]{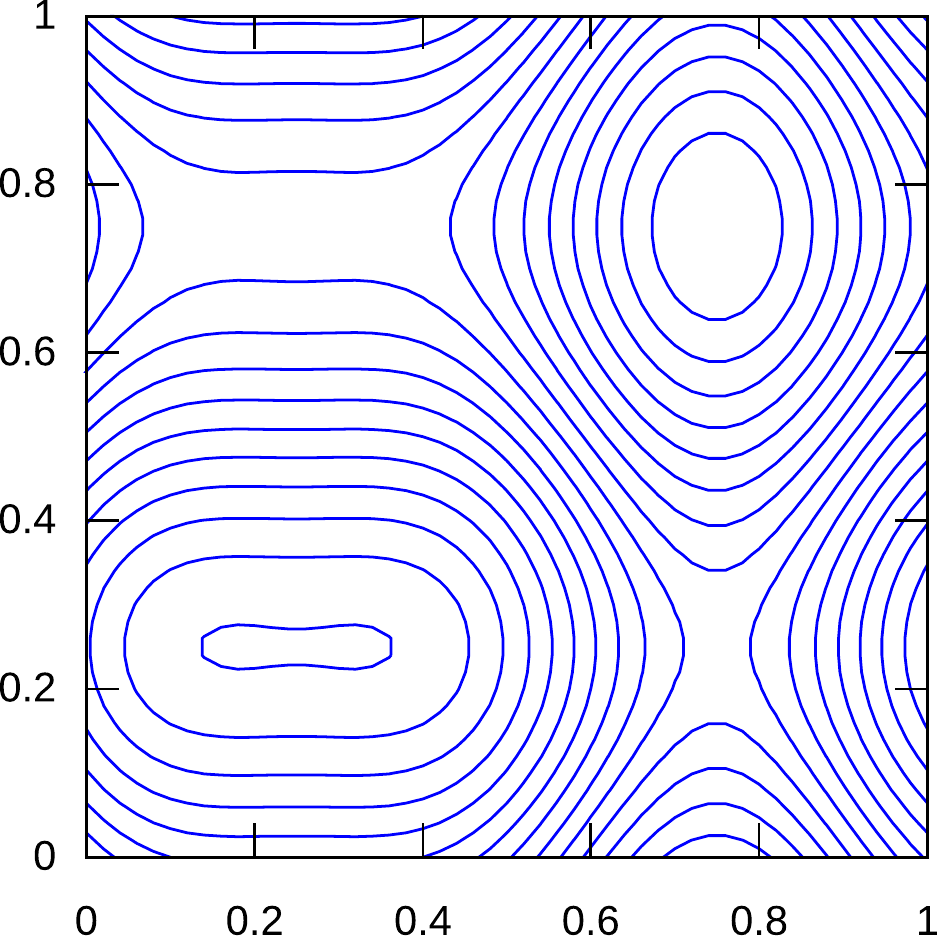} \\
(a)&(b)&(c)
\end{tabular}
\end{center}
\caption{$\Lambda$ vs number of iterations (a), level set of the solutions $U$ (b) and $M$ (c).}\label{mfg2d}
\end{figure}
In Figure \ref{mfg2d}a we show the computed $\Lambda$ as a function of the number of iterations to reach convergence.
The convergence is fast, just $5$ iterations in $8.06$ seconds. We get $\Lambda=0.9784$, which is exactly the same value reported in \cite{AC},
and the computed pair of solutions $(U,M)$, shown in Figure \ref{mfg2d}b and Figure \ref{mfg2d}c, is in qualitative agreement.

\begin{figure}[h!]
 \begin{center}
 \begin{tabular}{ccc}
\includegraphics[width=.225\textwidth]{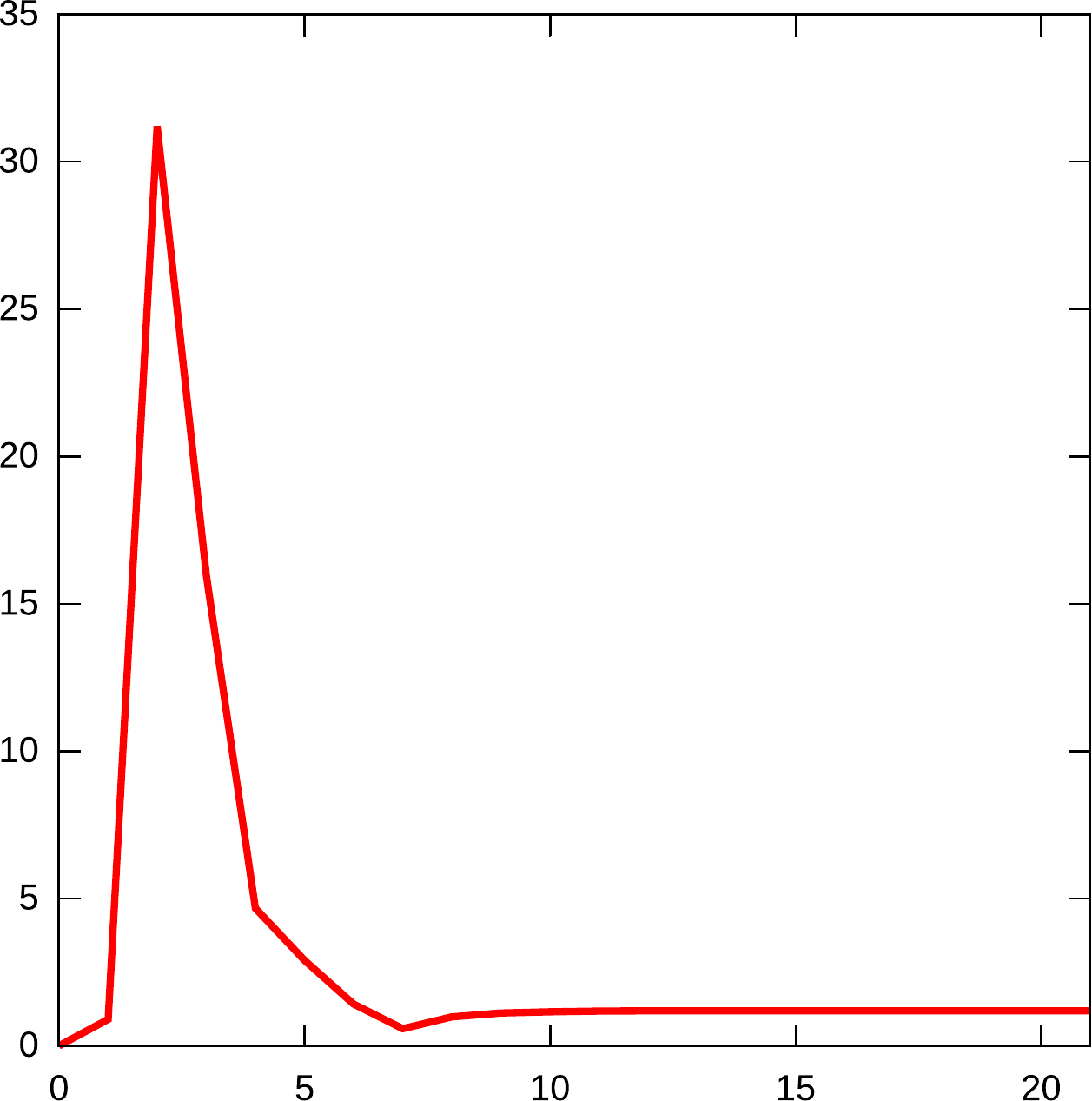} &
\includegraphics[width=.225\textwidth]{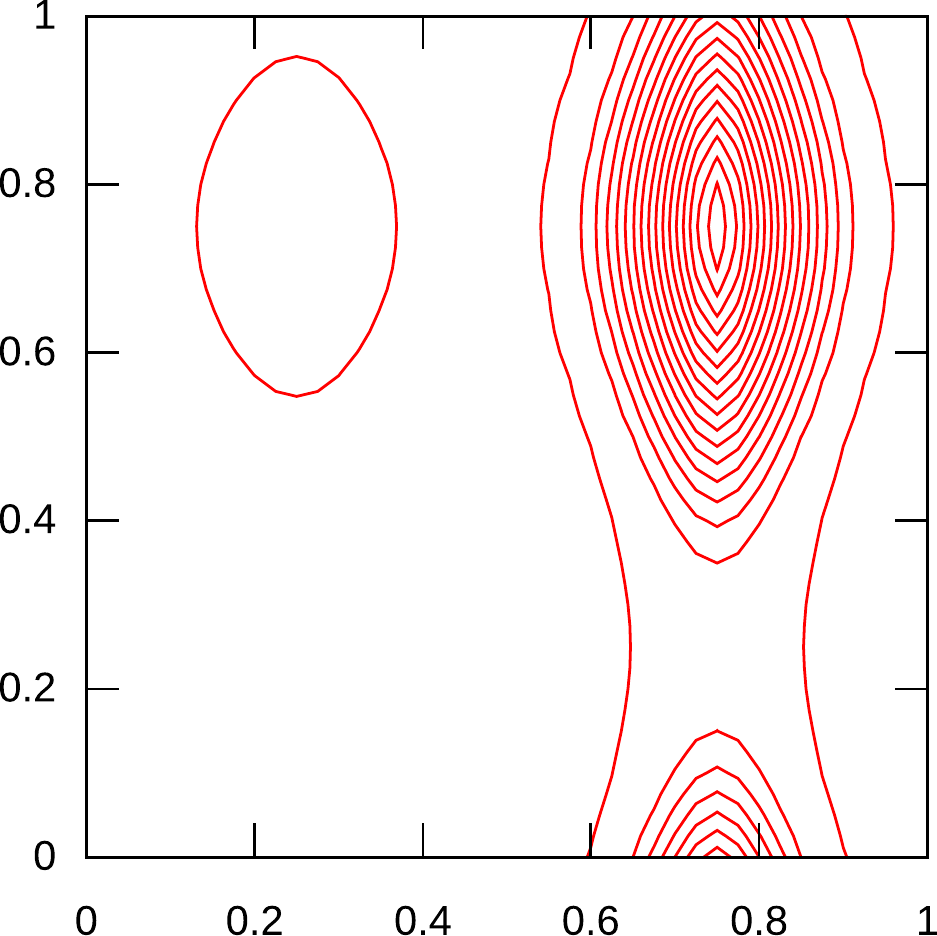} &
\includegraphics[width=.225\textwidth]{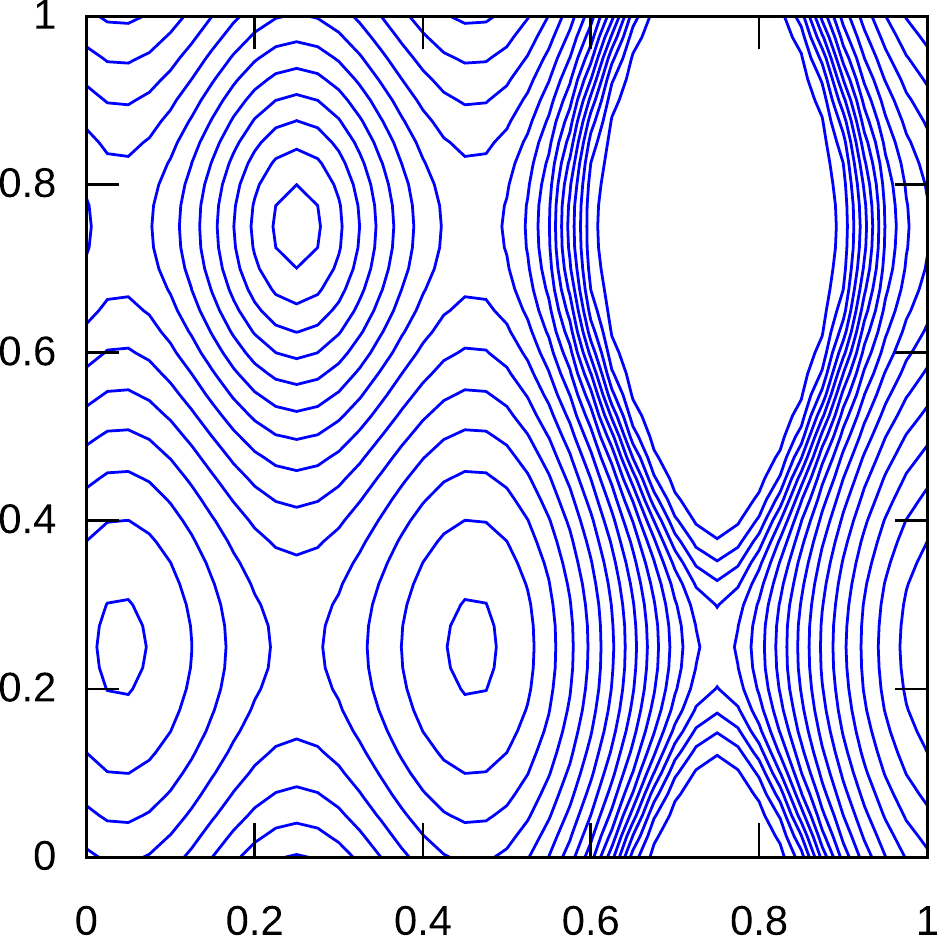} \\
(a)&(b)&(c)
\end{tabular}
\end{center}
\caption{$\Lambda$ vs number of iterations (a), level set of the solutions $U$ (b) and $M$ (c).}\label{mfg2d-2}
\end{figure}

We proceed with a case which is close to the deterministic limit, namely we repeat the previous test with $\nu=0.01$. We reach convergence in $21$ iterations
and $10.72$ seconds. In Figure \ref{mfg2d-2}a we show the computed $\Lambda$ as a function of the number of iterations, observing that the convergence is no longer
monotone. We get $\Lambda=1.1878$ which is again the same value reported in  \cite{AC}, as for the computed solutions $(U,M)$, shown in Figure \ref{mfg2d-2}b and
Figure \ref{mfg2d-2}c.

\begin{figure}[h!]
 \begin{center}
 \begin{tabular}{ccc}
\includegraphics[width=.225\textwidth]{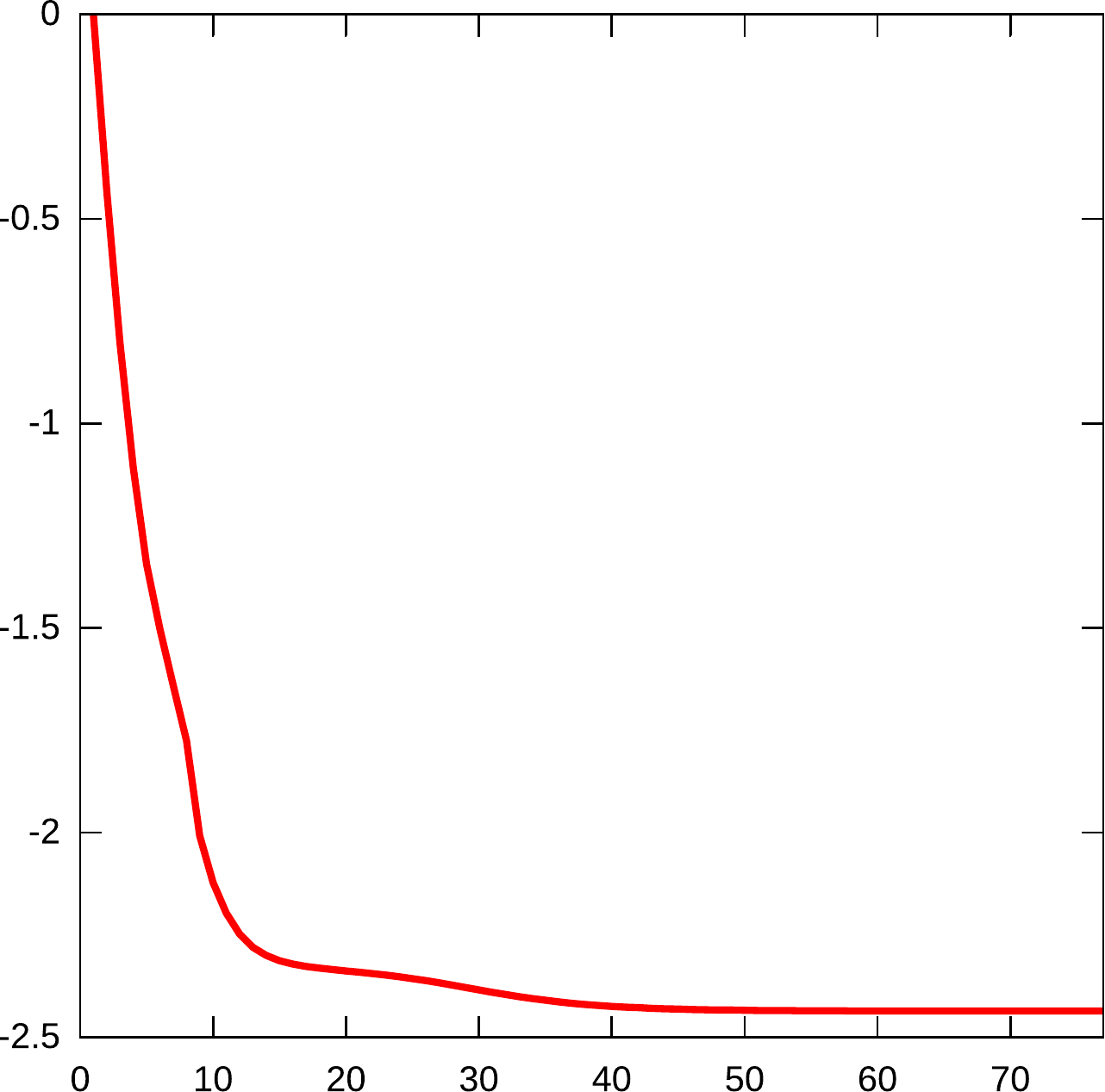} &
\includegraphics[width=.225\textwidth]{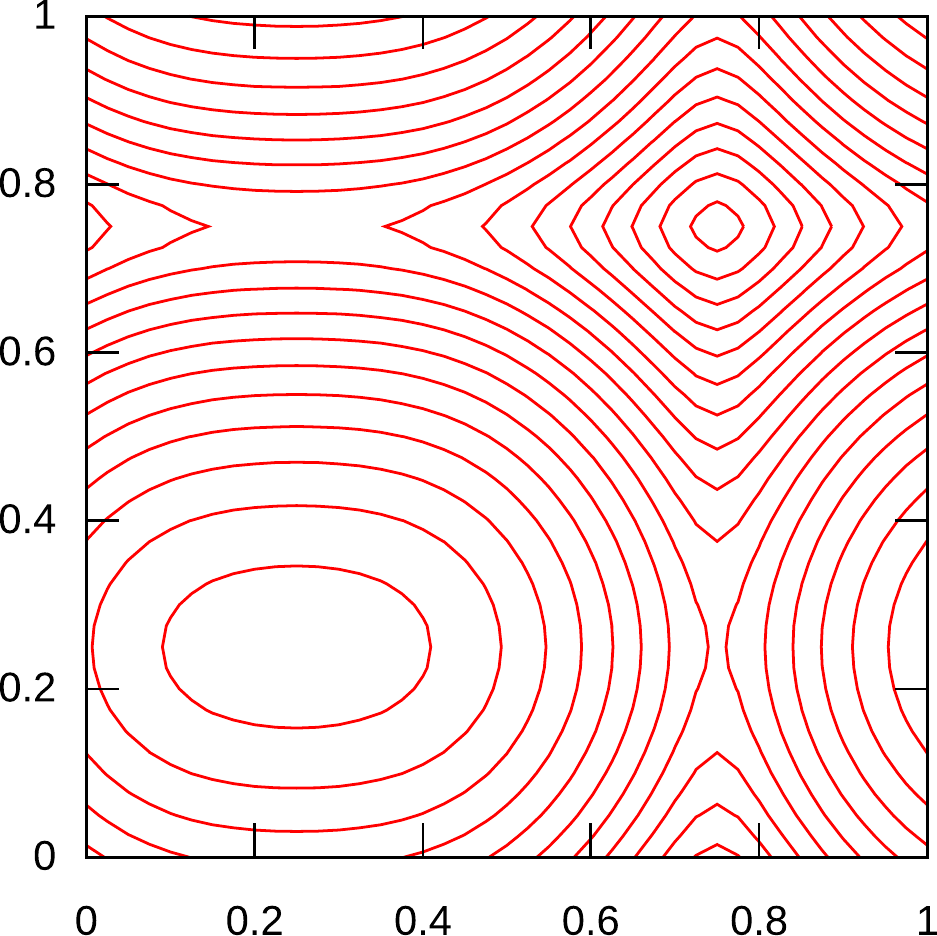} &
\includegraphics[width=.225\textwidth]{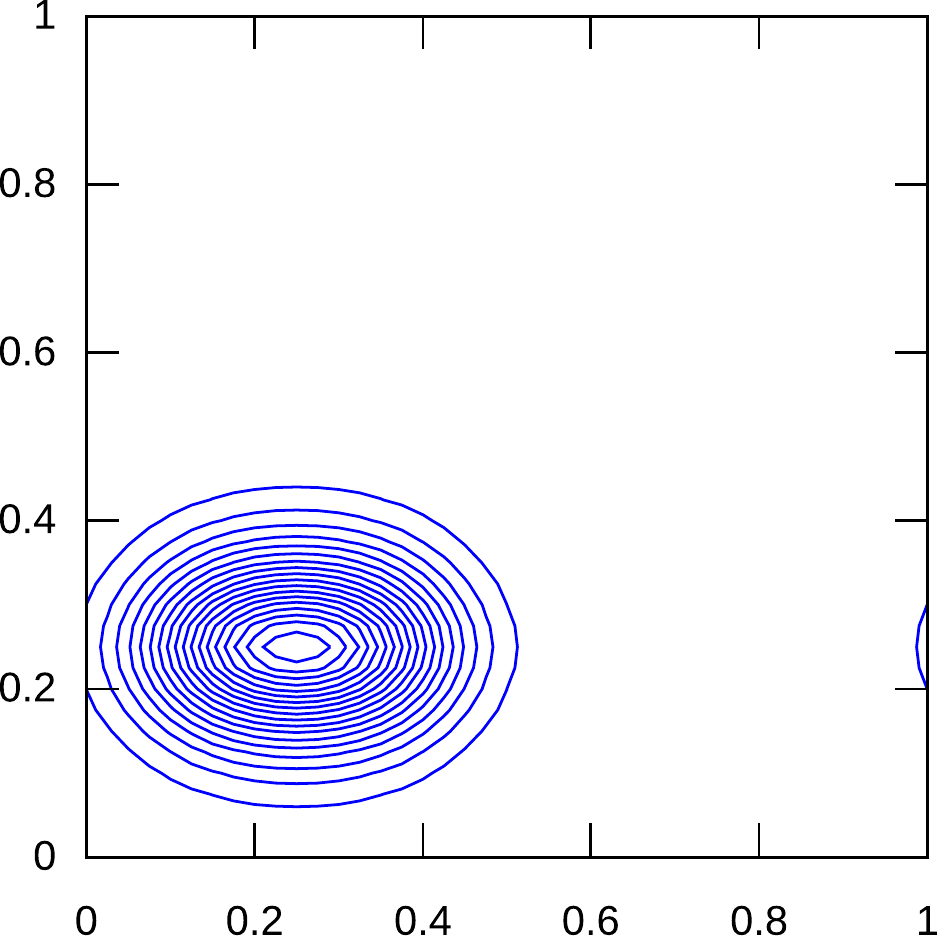} \\
(a)&(b)&(c)
\end{tabular}
\end{center}
\caption{$\l$ vs number of iterations (a), level set of the solutions $u$ (b) and $m$ (c).}\label{mfg2d-3}
\end{figure}

We conclude with the case of a nonincreasing potential $V(m(x))=-\log(m(x))$ and we choose $\nu=0.1$. The convergence is much slower than in the
previous tests, as shown in Figure \ref{mfg2d-3}a. In \cite{AC} it is not reported the computed value of $\Lambda$,
whereas we get $\Lambda=-2.4358$ in $77$ iterations and $42.33$ seconds. On the other hand,
the computed solutions $(U,M)$ are in qualitative agreement and they are shown in Figure \ref{mfg2d-3}b and Figure \ref{mfg2d-3}c.
\\
Finally, we summarize the results for reference in Table \ref{table-mfg}.
\begin{table}[!ht]
 \centering
  \begin{tabular}{|c|c|c|c|c|}
    \hline
    $V(m)$ & $\nu$ & $\Lambda$ & Iterations & Total CPU (secs)\\
    \hline
    $m^2$&1&0.9784 &5 &8.06\\\hline
    $m^2$&0.01&1.1878 &21 &10.72\\\hline
    $-log(m)$&0.1&-2.4358 &77 &42.33\\\hline
 \end{tabular}
  \caption{Performance for the 2D MFG.}\label{table-mfg}
\end{table}\\
\section{Multi-population stationary mean field games}\label{multipopMFG}
This is a generalization of \eqref{MFG} to the case of $P$ populations, each one described by a MFG-system, coupled via a potential term (see \cite{LL}).
We consider the setting recently studied in \cite{C15} for problems with Neumann boundary conditions.

Here we present the simple case in dimension one and two of an eikonal Hamiltonian with a linear local potential, namely the problem
$$
  \left\{
   \begin{array}{ll}
    -\nu \Delta  u_i  +|D u_i|^2+\l_i  =V_i(m) \qquad & \mbox{in }\Omega\,,\quad i=1,...,P\\[4pt]
     \nu \Delta m_i +2\mathrm{div}(m_i\,D u_i)=0 & \mbox{in }\Omega\,,\quad i=1,...,P\\[4pt]
     \partial_n u_i=0\,,\quad\partial_n m_i=0  & \mbox{on }\partial\Omega\,,\quad i=1,...,P\\[4pt]
   \int_{\Omega} u_i(x)dx=0\,,\quad\int_{\Omega} m_i(x)dx=1\,,\quad m_i\ge 0     & i=1,...,P\,,
    \end{array}
  \right.
$$
where $\Omega=[0,1]$ or $\Omega=[0,1]^2$, the corrector $u(x)=(u_1(x),...,u_P(x))$ and the mass density
$m(x)=(m_1(x),...,m_P(x))$ are vector functions and $\l=(\l_1,...,\l_P)\in\R^P$ is a $P$-tuple of ergodic constants. Moreover, for $i=1,...,P$,
the linear local potential $V_i$ takes the
form $$V_i(m(x))=\sum_{j=1}^P \theta_{ij} m_j(x)\,,$$ for some given weights $\theta_{ij}\in\R$, or in matrix notation
\begin{equation}\label{theta}
V=(V_1,...,V_P)\,,\qquad \Theta=(\theta_{ij})_{i,j=1,...,P}\,,\qquad V(m)=\Theta m\,.
\end{equation}
Existence and uniqueness of a solution $(u,m,\l)$ can be proved under some monotonicity assumptions on $V$ (see \cite{C15} for details).

Within our framework of inconsistent systems, the problem is again {\em overdetermined}. Indeed, discretizing $\Omega$
with a uniform grid of $N^n$ nodes ($n=1,2$), we end up with $P(2N^n+2)$ equations in the $P(2N^n+1)$ unknowns $(U,M,\Lambda)$.
Again, we do not include the constraint $M\geq 0$ in the discretization, as before the normalization condition for $M$ seems enough
to force numerically its nonnegativity.

In the special case \eqref{theta} uniqueness is guaranteed assuming that $\Theta$ is positive semi-definite and the solution is explicitly given,
for $i=1,...,P$, by $u_i\equiv 0$, $m_i\equiv 1$ and $\l_i=\sum_{j=1}^P\theta_{ij}$.
By dropping this condition, the trivial solution is still found, but we expect to observe other more interesting solutions.

We start with some experiments in dimension $n=1$, in the case of $P=2$ populations. We choose the coupling matrix (not positive semi-definite)
$$\Theta=\left(\begin{array}{cc}
                0 & 1 \\ 1 & 0
               \end{array}
\right)$$
so that the potential for each population only depends on the other population. Moreover, we discretize the interval $\Omega=[0,1]$ with $N=100$
uniformly distributed nodes,
we set to $\nu=0.05$ the diffusion coefficient  and to $\varepsilon=10^{-6}$ the tolerance for the stopping criterion of the algorithm.
To avoid the trivial solution, we choose non constant initial guesses, such as
piecewise constant pairs with zero mean for $U$ and piecewise constant pairs with mass one for $M$.

\begin{figure}[h!]
 \begin{center}
 \begin{tabular}{cccc}
\includegraphics[width=.2\textwidth]{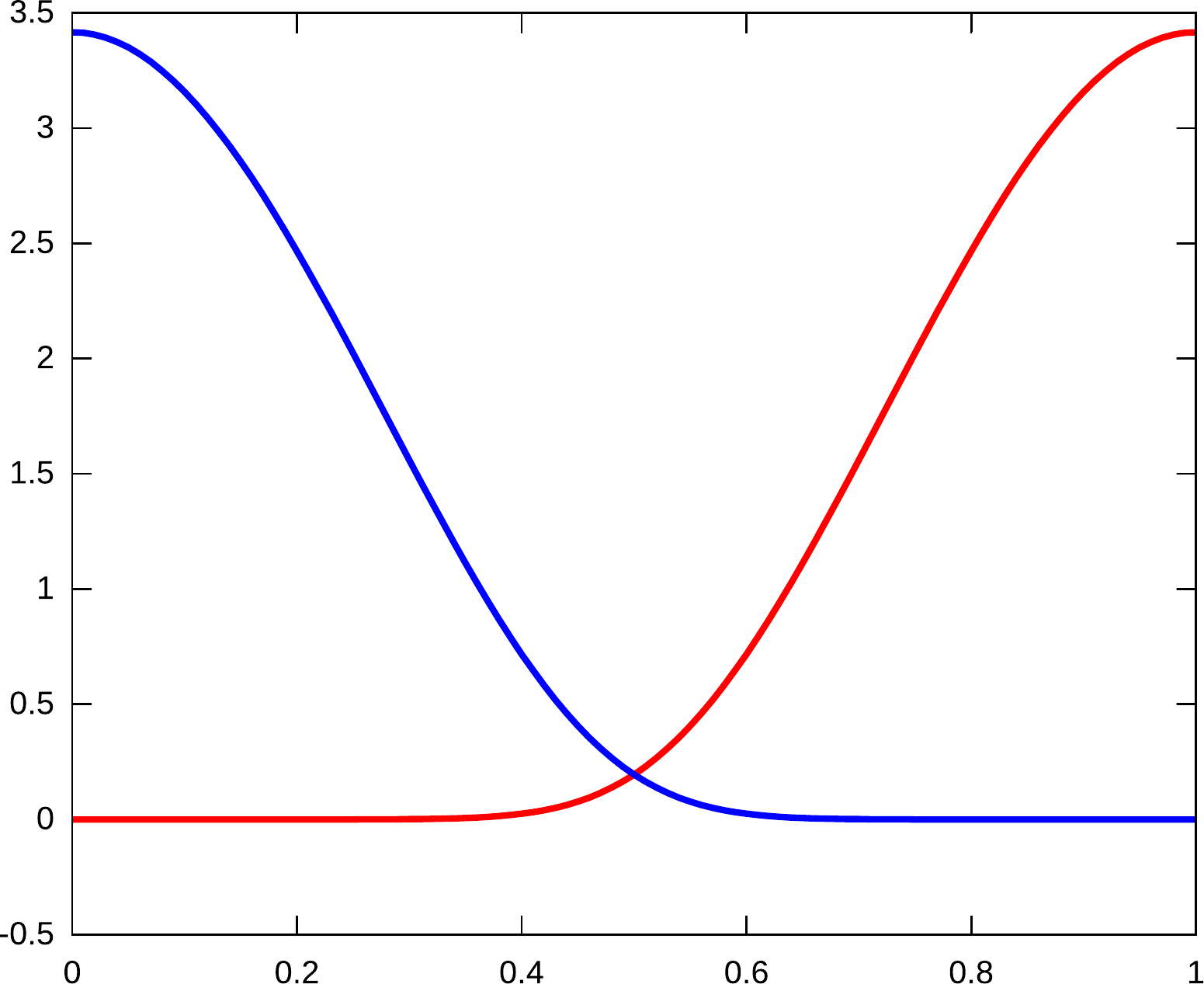} &
\includegraphics[width=.2\textwidth]{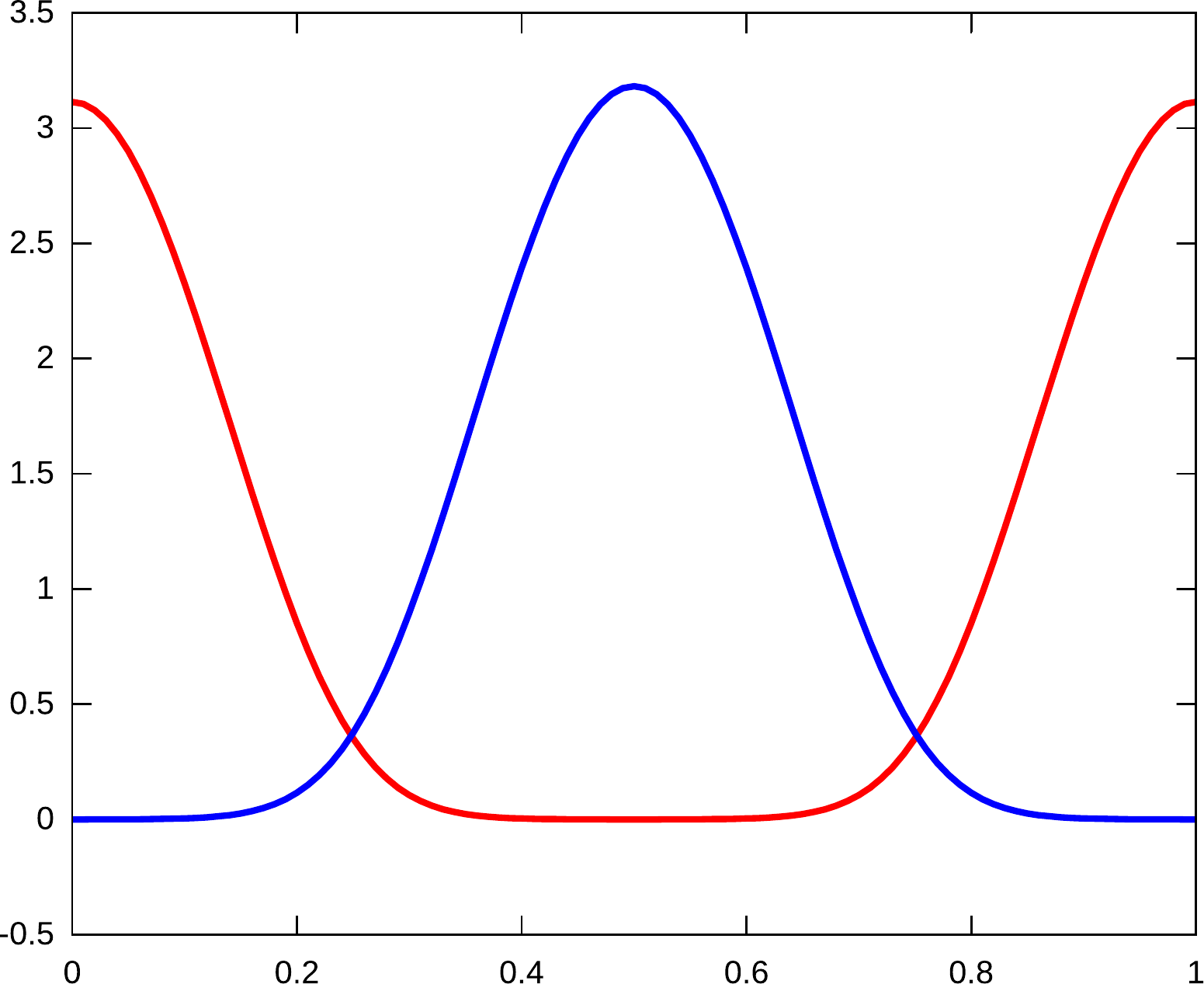} &
\includegraphics[width=.2\textwidth]{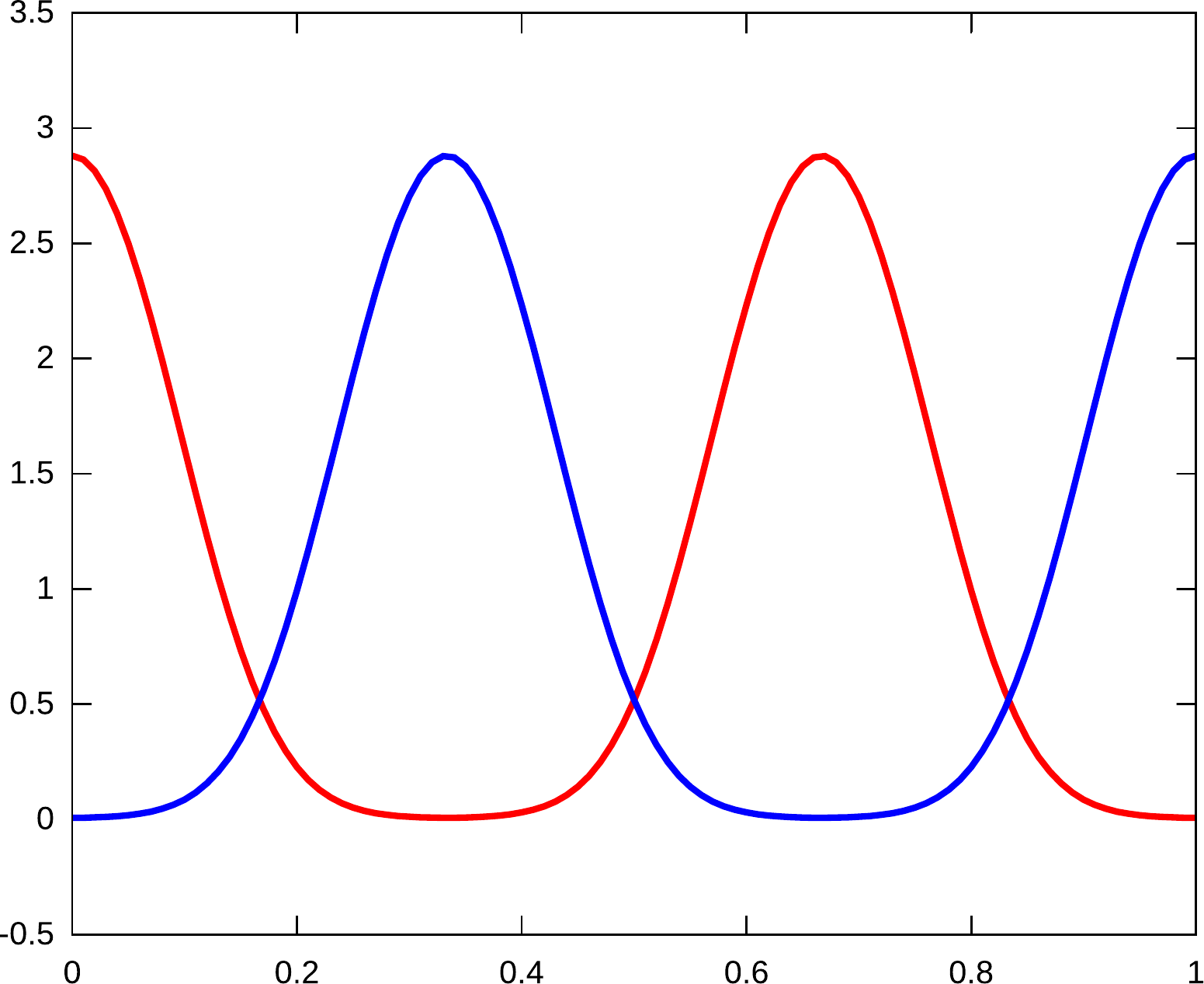} &
\includegraphics[width=.2\textwidth]{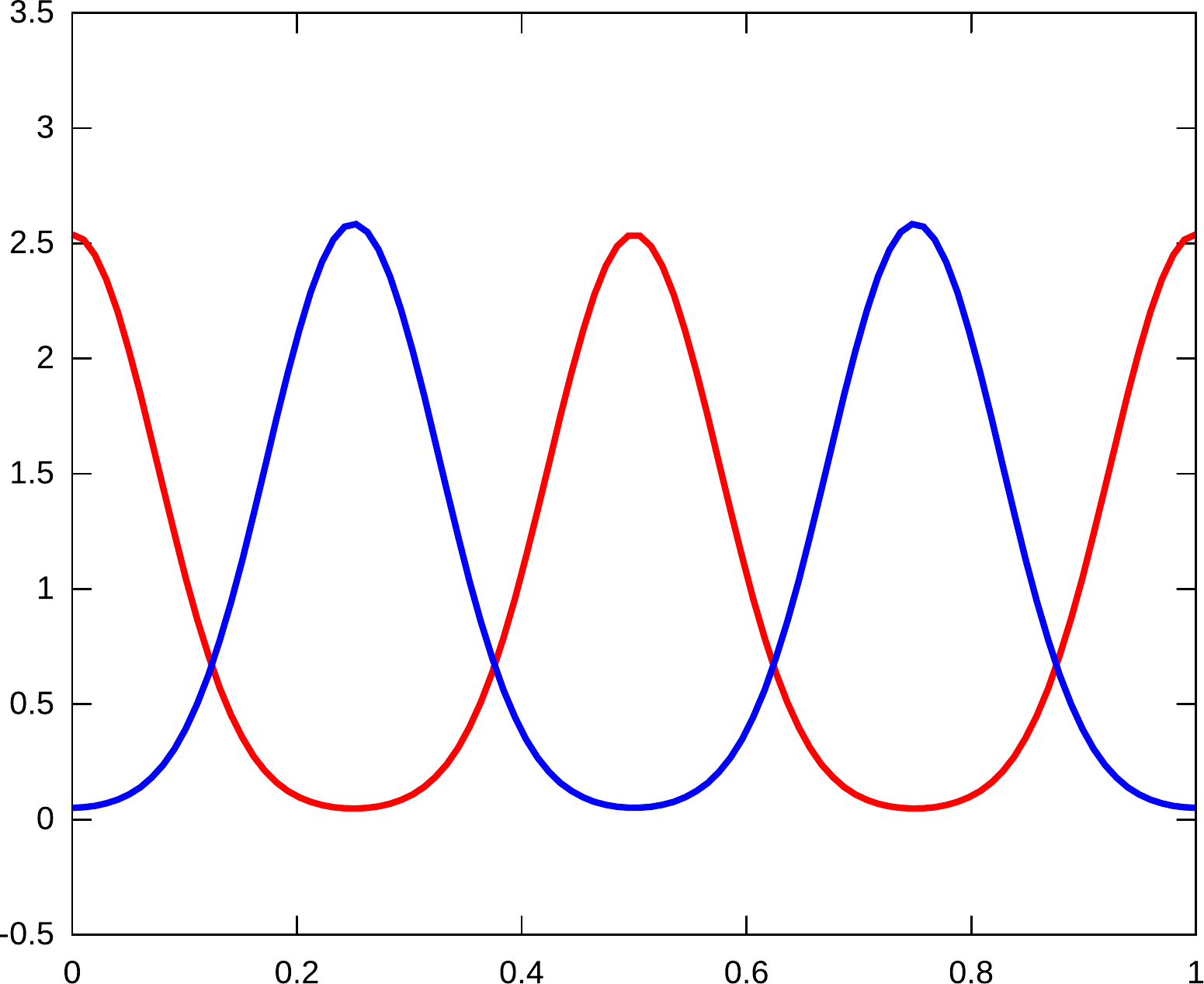} \\
\includegraphics[width=.2\textwidth]{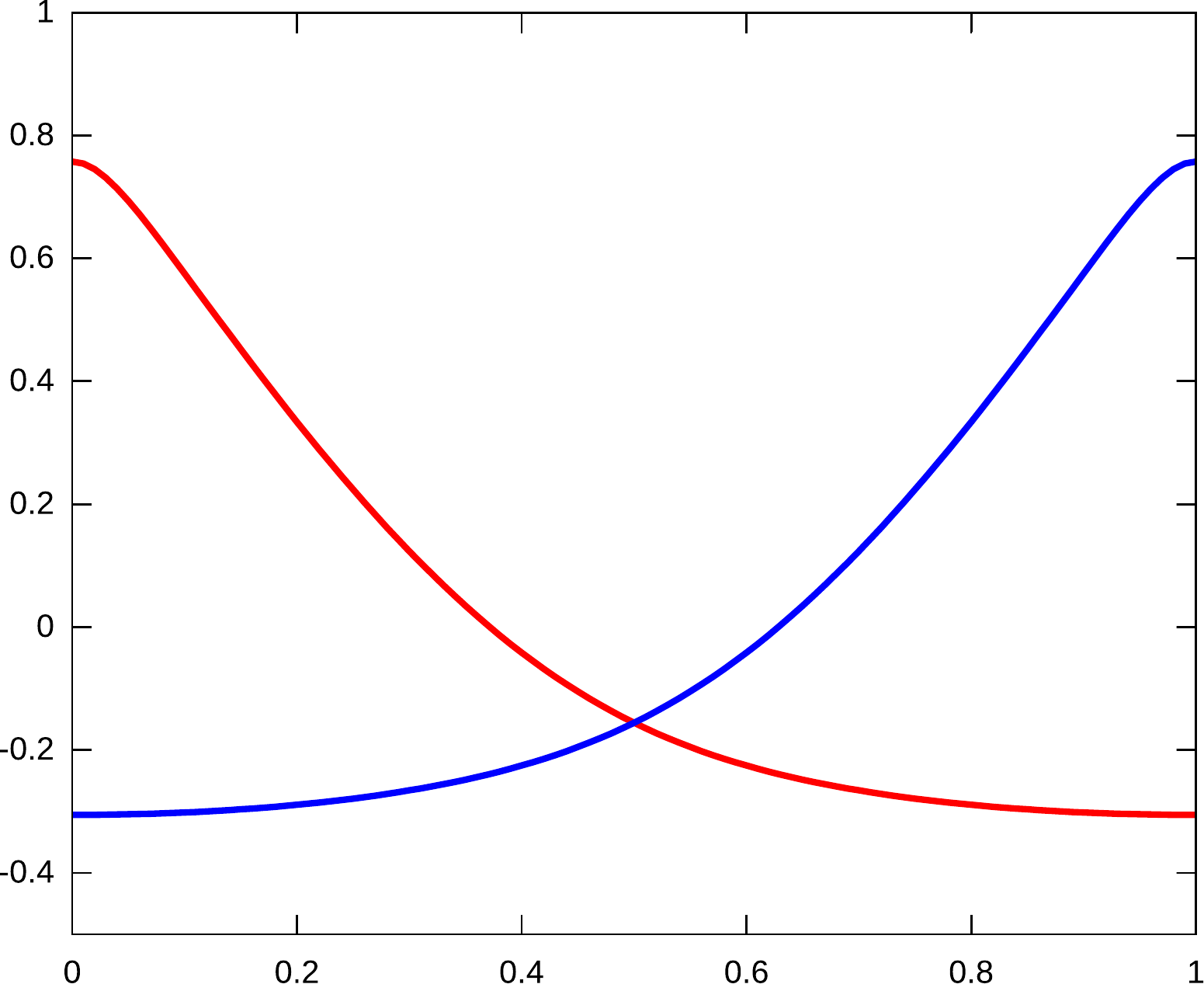} &
\includegraphics[width=.2\textwidth]{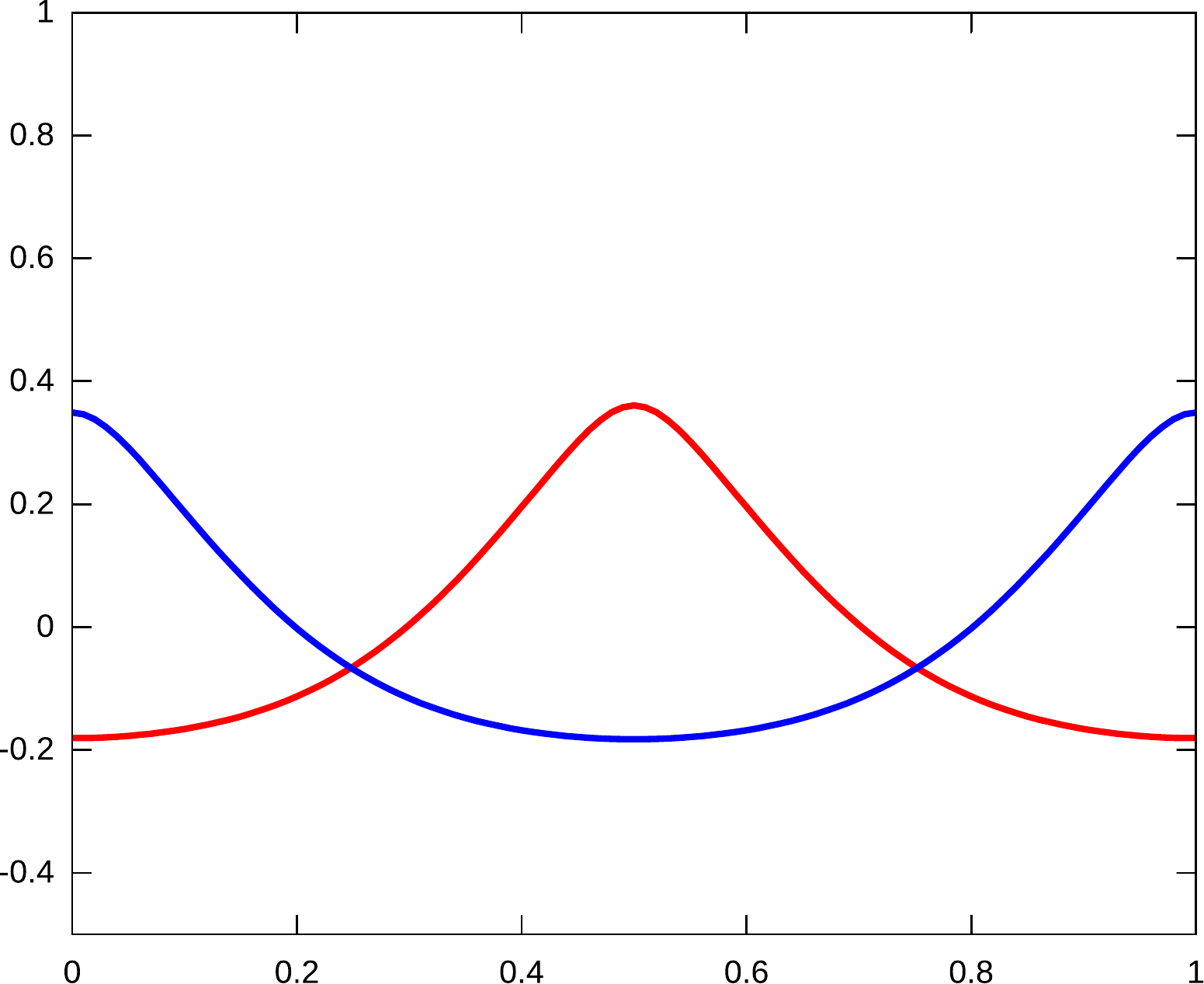} &
\includegraphics[width=.2\textwidth]{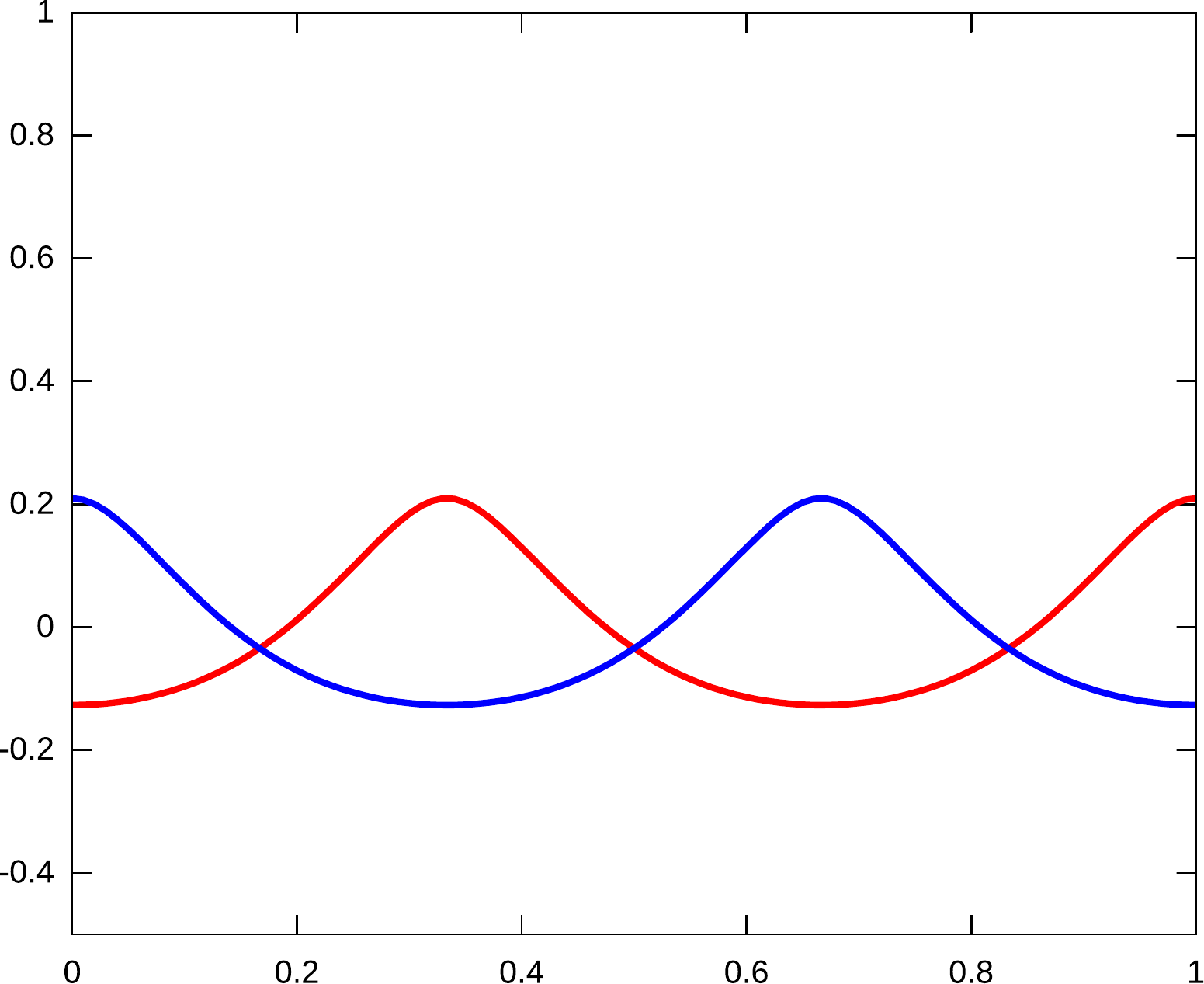} &
\includegraphics[width=.2\textwidth]{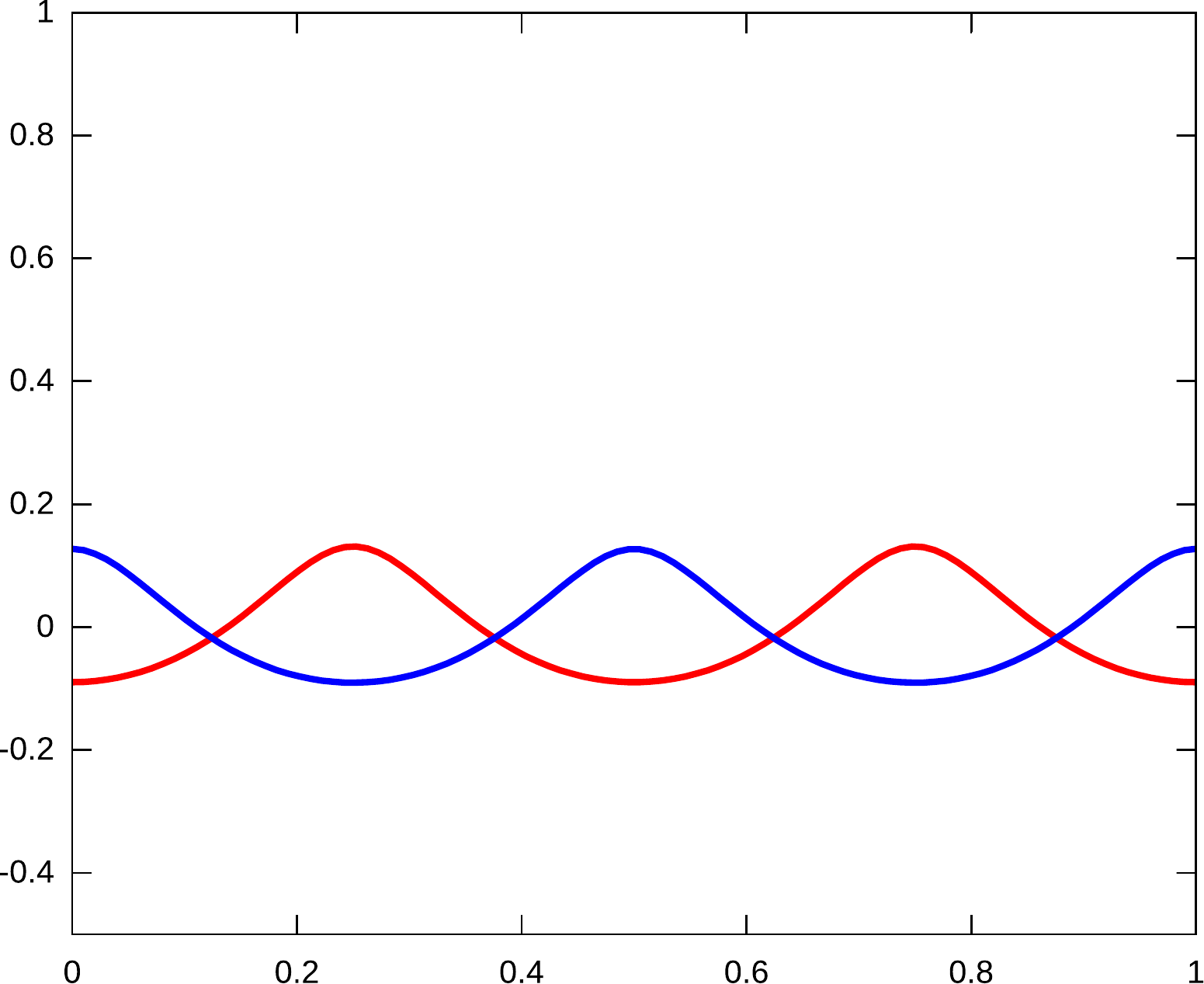} \\
(a)&(b)&(c)&(d)
\end{tabular}
\end{center}
\caption{Two-population MFG solutions ($\nu=0.05$): mass densities $M=(M_1,M_2)$ (top panels) and corresponding correctors $U=(U_1,U_2)$ (bottom panels).}\label{mfgm1d-1}
\end{figure}

We did not find a general criterion to select some specific
solution, but we experienced that the algorithm is very sensitive to the grid size.
Figure \ref{mfgm1d-1} shows four computed solutions. In the top panels we show the densities $M=(M_1,M_2)$, while in the bottom panels  the corresponding
correctors $U=(U_1,U_2)$. Moreover, in Table \ref{table-mfgm1d-1} we report the performance of the method for each case, including the computed value of
the ergodic pair $\Lambda=(\Lambda_1,\Lambda_2)$. \\
\begin{table}[!ht]
 \centering
  \begin{tabular}{|c|c|c|c|}
    \hline
    Test & $\Lambda=(\Lambda_1,\Lambda_2)$ & Iterations & Total CPU (secs)\\
    \hline
    ($a$)&$(0.03921,0.03921)$ &19 &0.14\\\hline
    ($b$)&$(0.14362,0.14048)$ &48 &0.33\\\hline
    ($c$)&$(0.29498,0.29498)$ &28 &0.19\\\hline
    ($d$)&$(0.49574,0.48481)$ &31 &0.21\\\hline
 \end{tabular}
  \caption{Performance for the 1D two-population MFG.}\label{table-mfgm1d-1}
  \vspace{-15pt}
\end{table}

\begin{figure}[h!]
 \begin{center}
 \begin{tabular}{cccc}
\includegraphics[width=.2\textwidth]{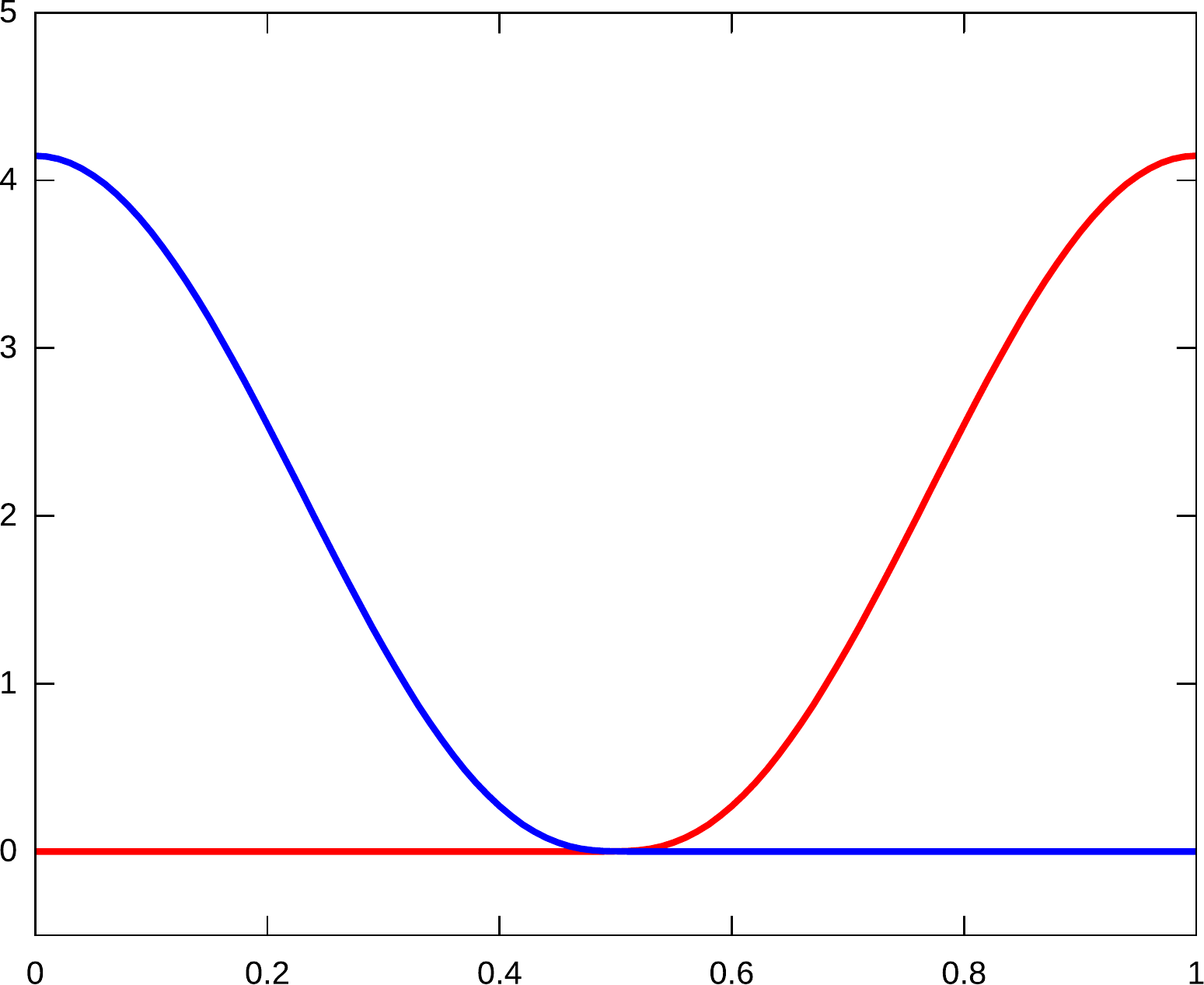} &
\includegraphics[width=.2\textwidth]{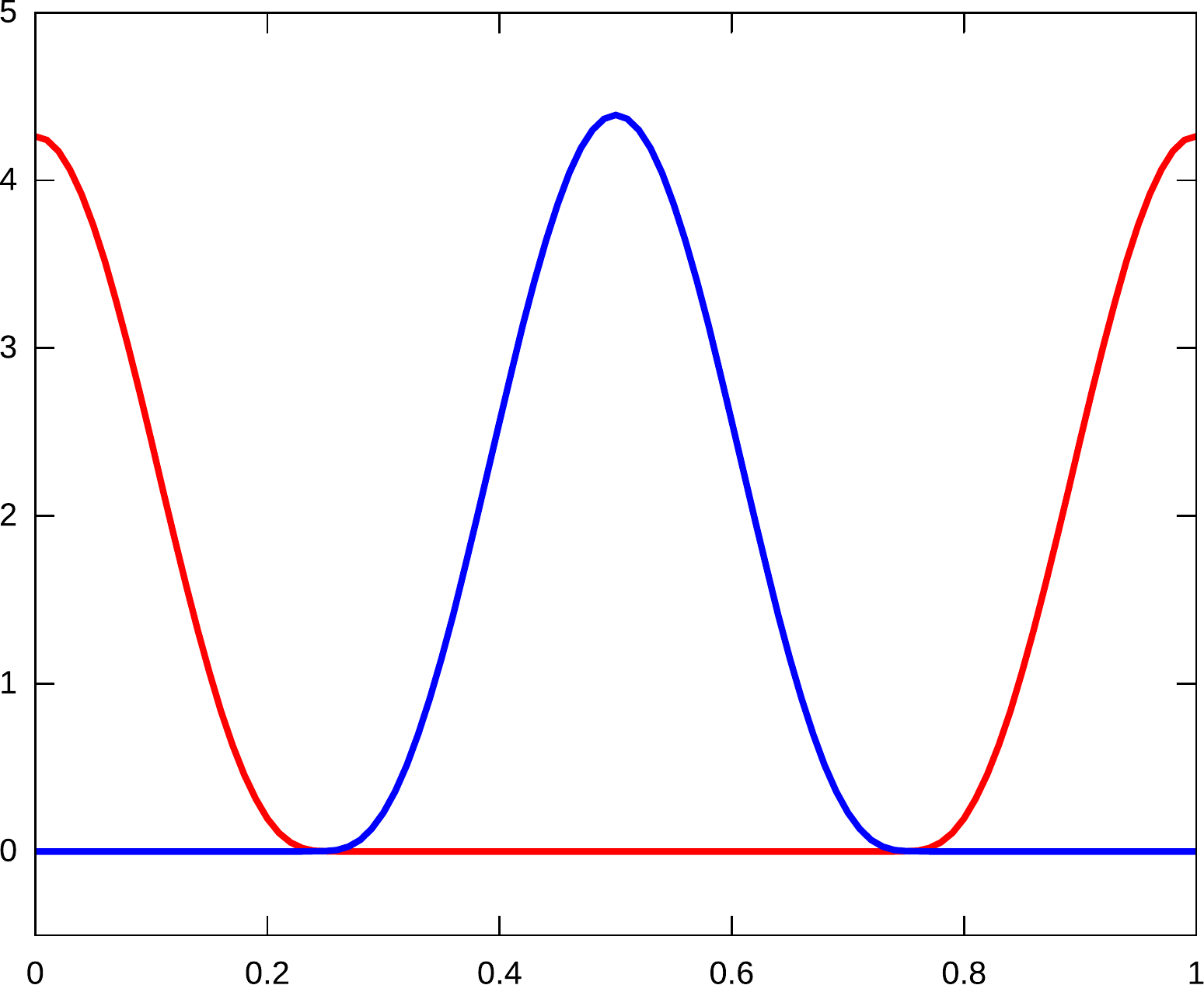} &
\includegraphics[width=.2\textwidth]{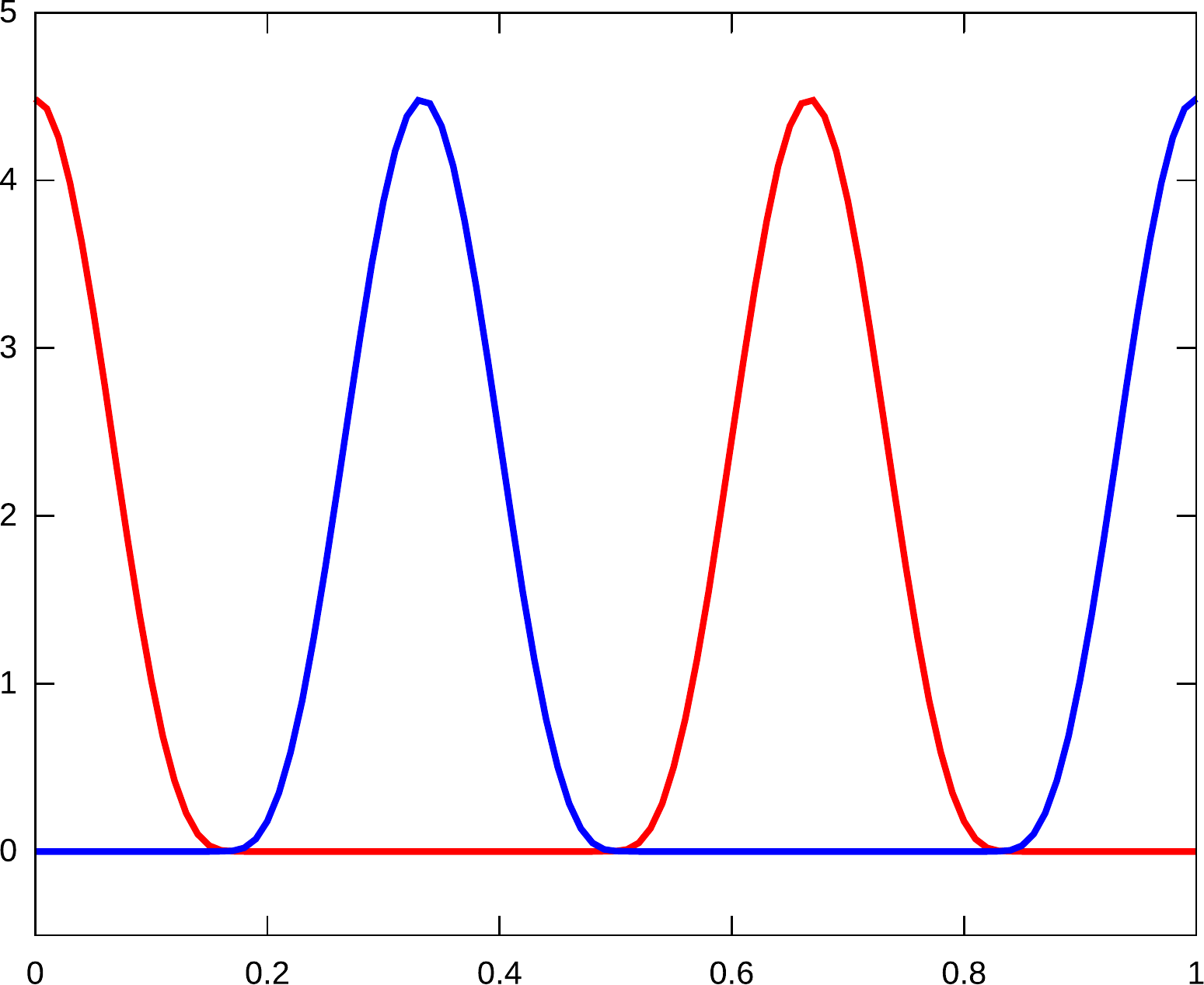} &
\includegraphics[width=.2\textwidth]{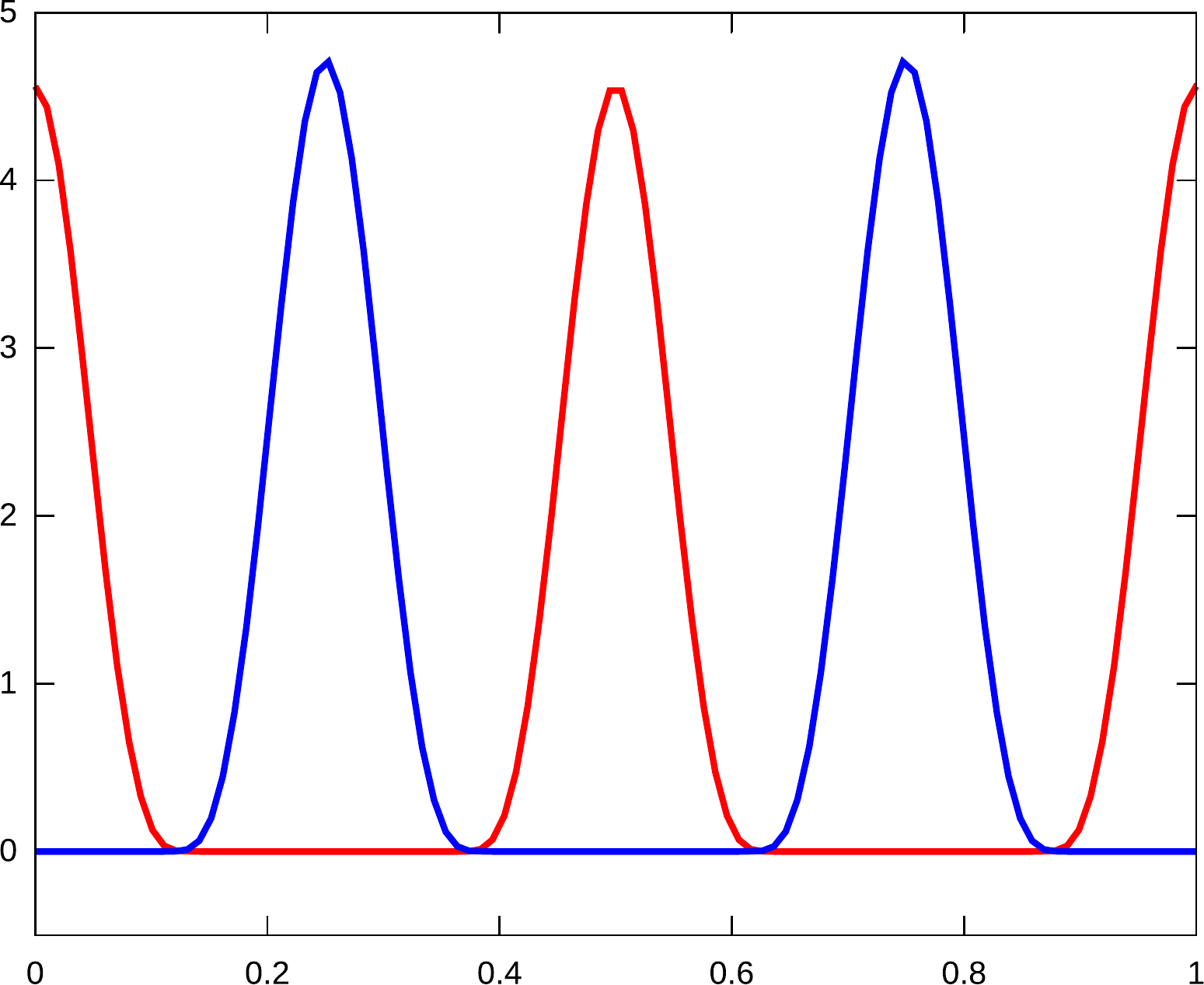} \\
\includegraphics[width=.2\textwidth]{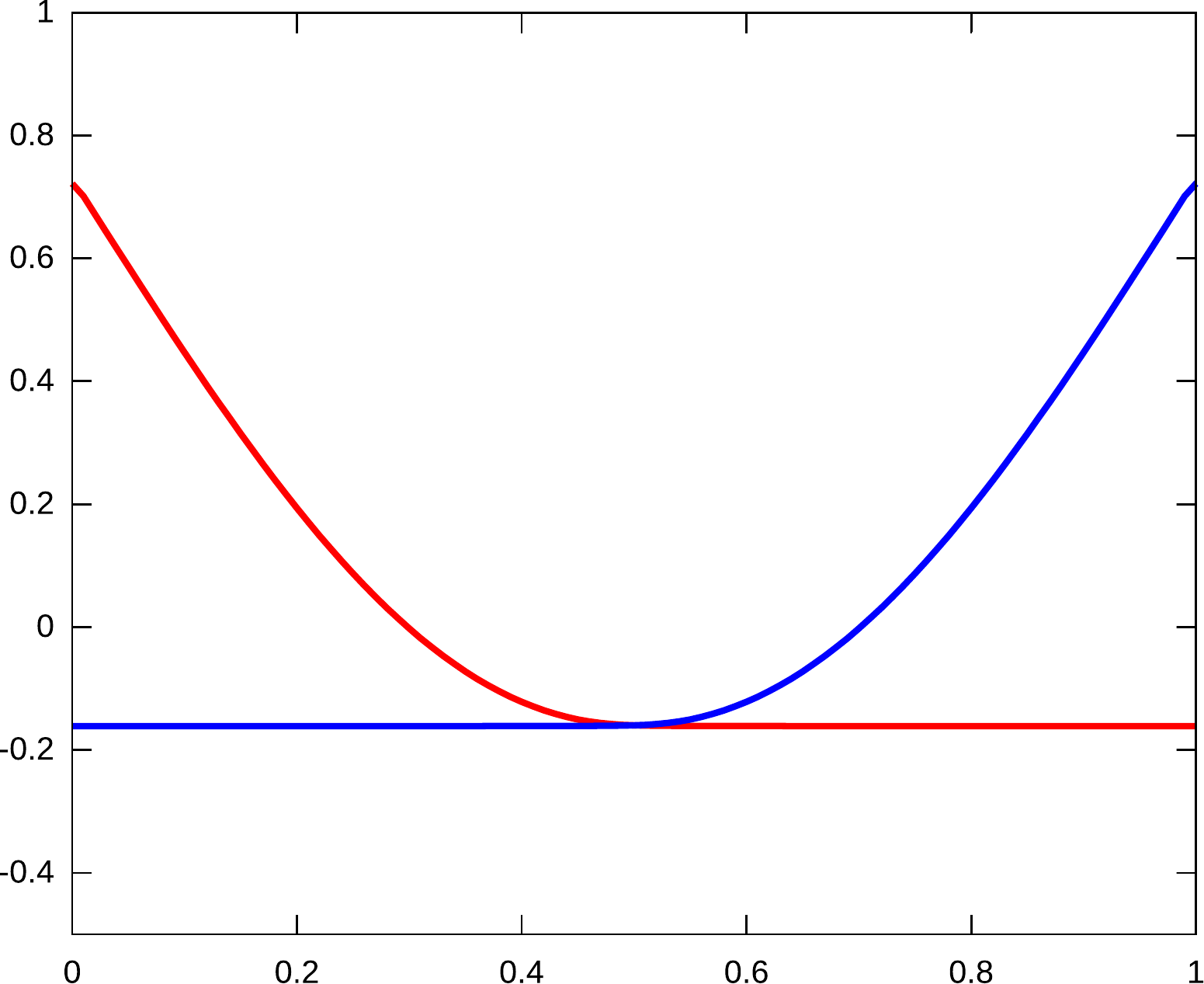} &
\includegraphics[width=.2\textwidth]{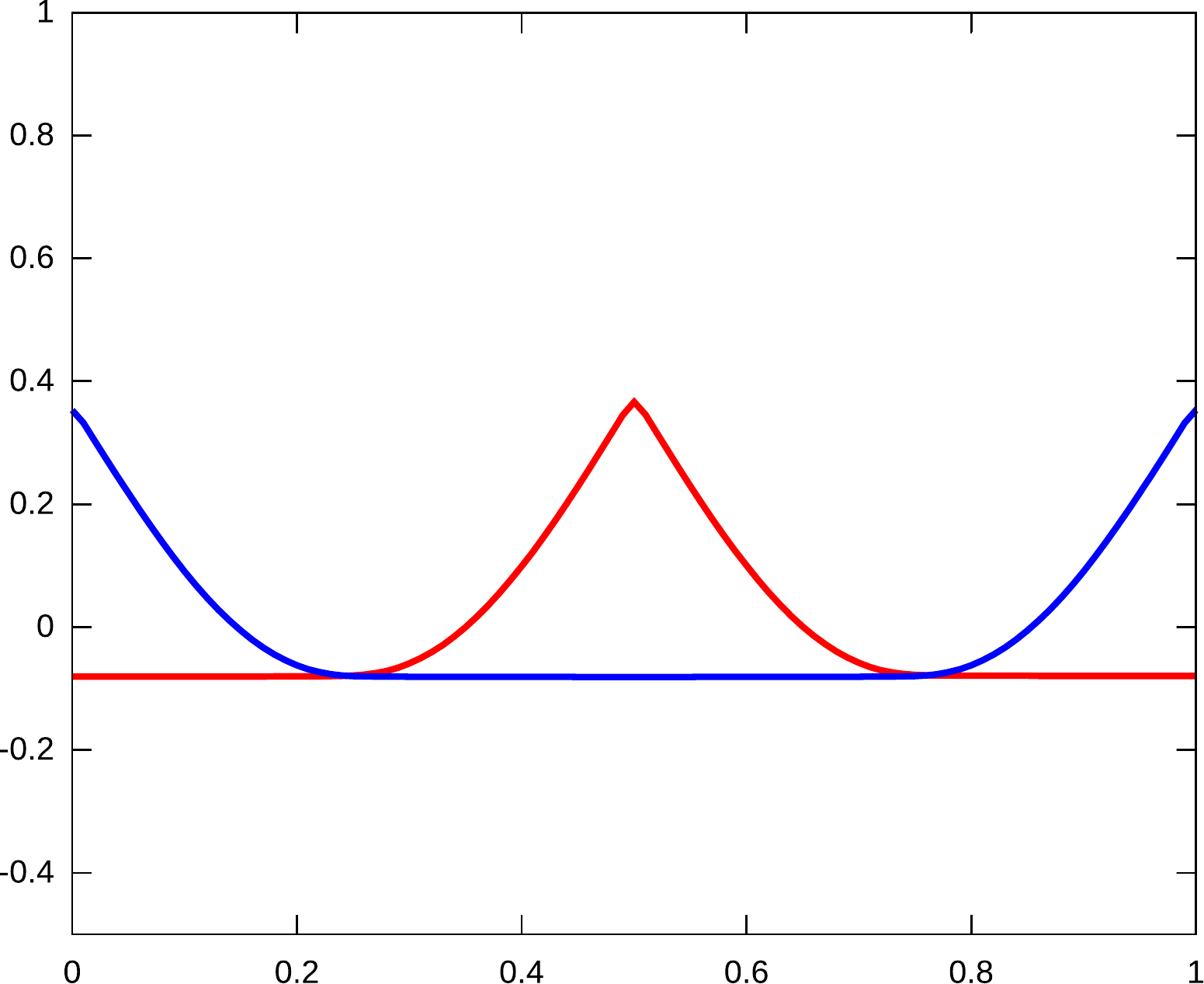} &
\includegraphics[width=.2\textwidth]{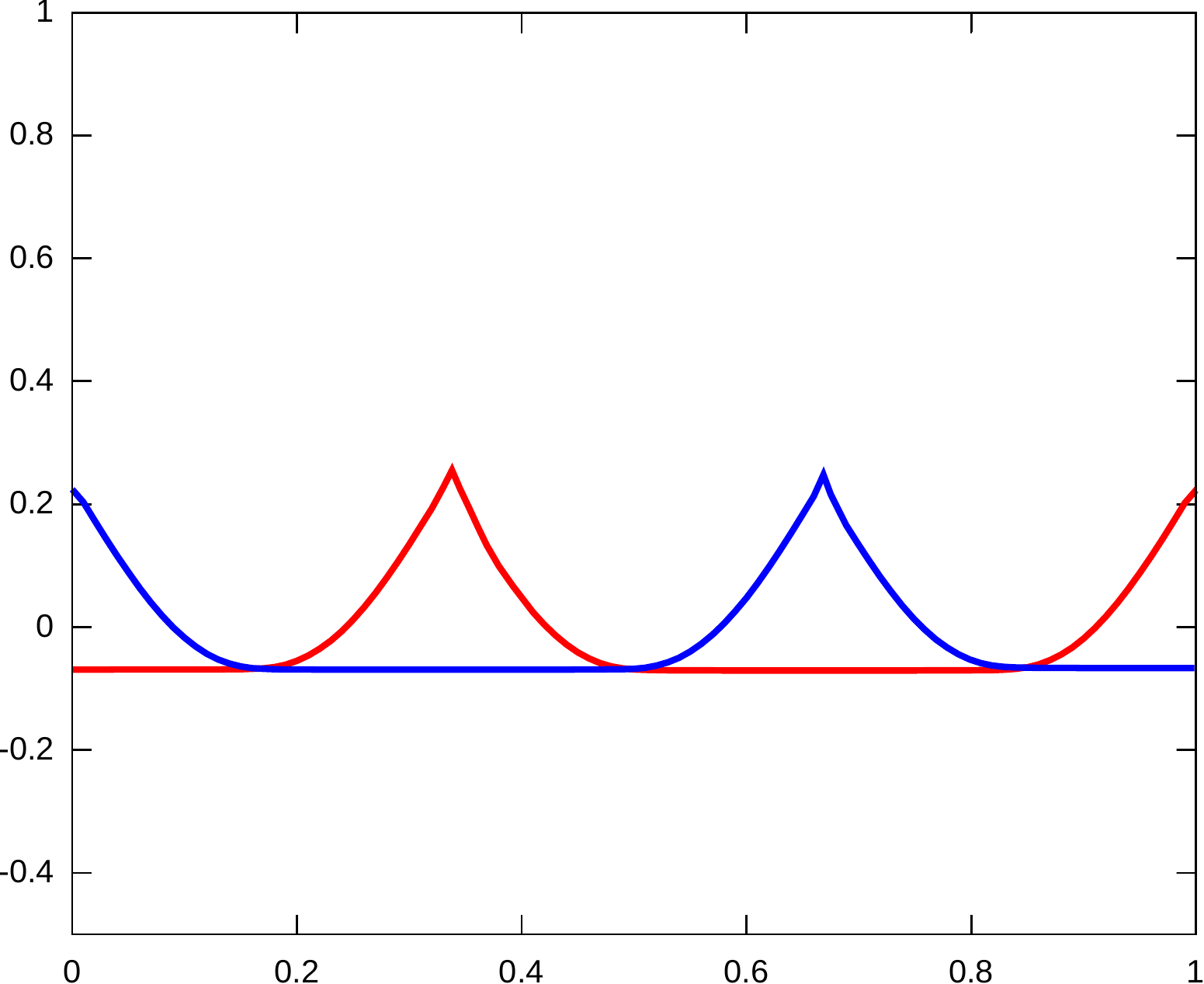} &
\includegraphics[width=.2\textwidth]{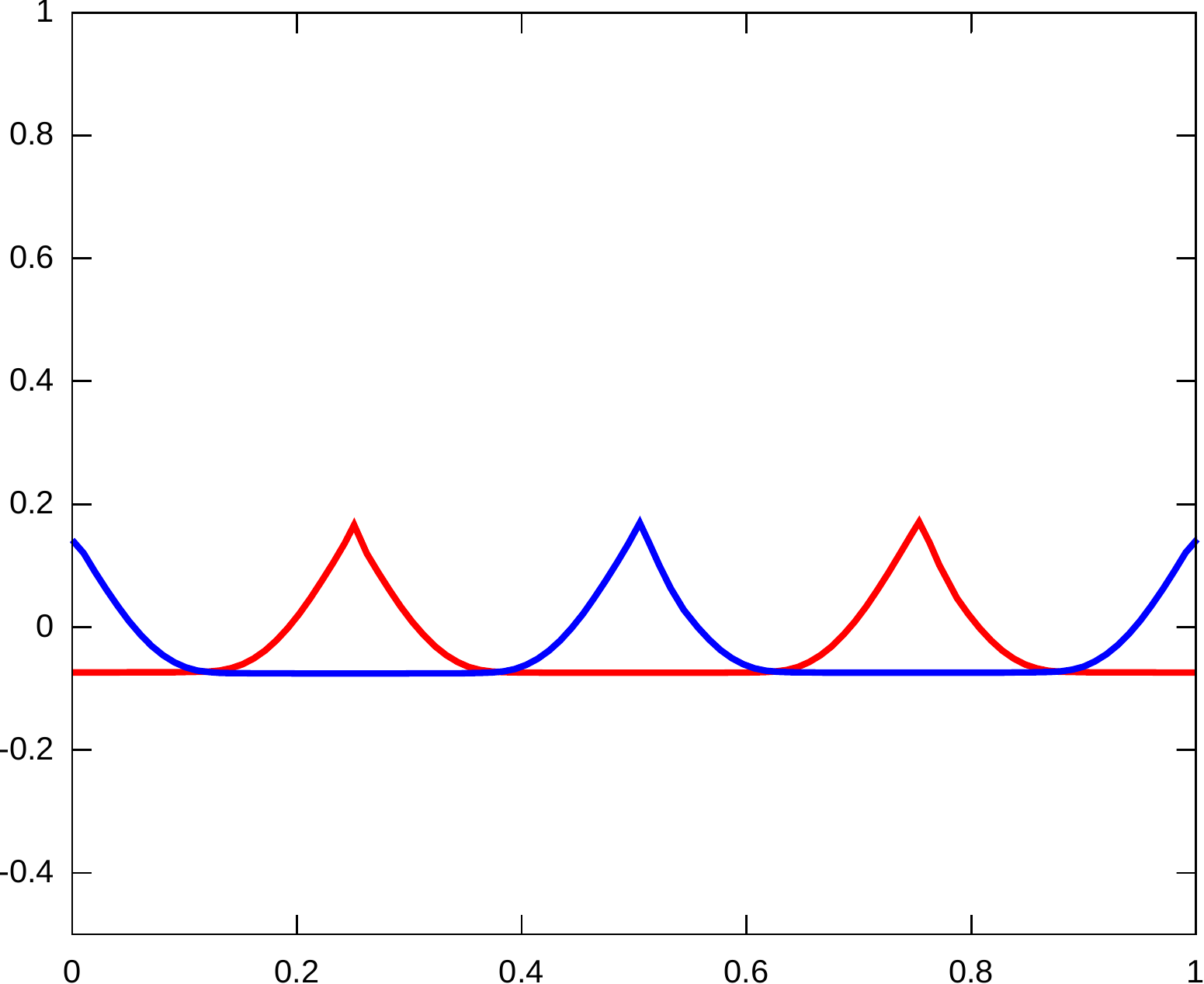} \\
(a)&(b)&(c)&(d)
\end{tabular}
\end{center}
\caption{Two-population MFG solutions ($\nu=10^{-4}$): mass densities $M=(M_1,M_2)$ (top panels) and corresponding correctors $U=(U_1,U_2)$ (bottom panels).}
\label{mfgm1d-1e0}
\end{figure}

Segregation of the two populations is expected (see \cite{C15}) and clearly visible.
This phenomenon can be enhanced by reducing the diffusion coefficient, as shown in Figure \ref{mfgm1d-1e0}, where $\nu=10^{-4}$, close to the deterministic limit.
\\

We finally consider the more complex and suggestive two dimensional case. We discretize the square $\Omega=[0,1]^2$ with $25\times 25$
uniformly distributed nodes and we push the diffusion $\nu$ up to $10^{-6}$, in order to observe segregation among the populations.
Moreover, we choose the interaction matrix as before, with all the entries equal to $1$ except for the diagonal, which is set to $0$.

\begin{figure}[h!]
  \begin{center}
 \begin{tabular}{ccccc}
\includegraphics[width=.16\textwidth]{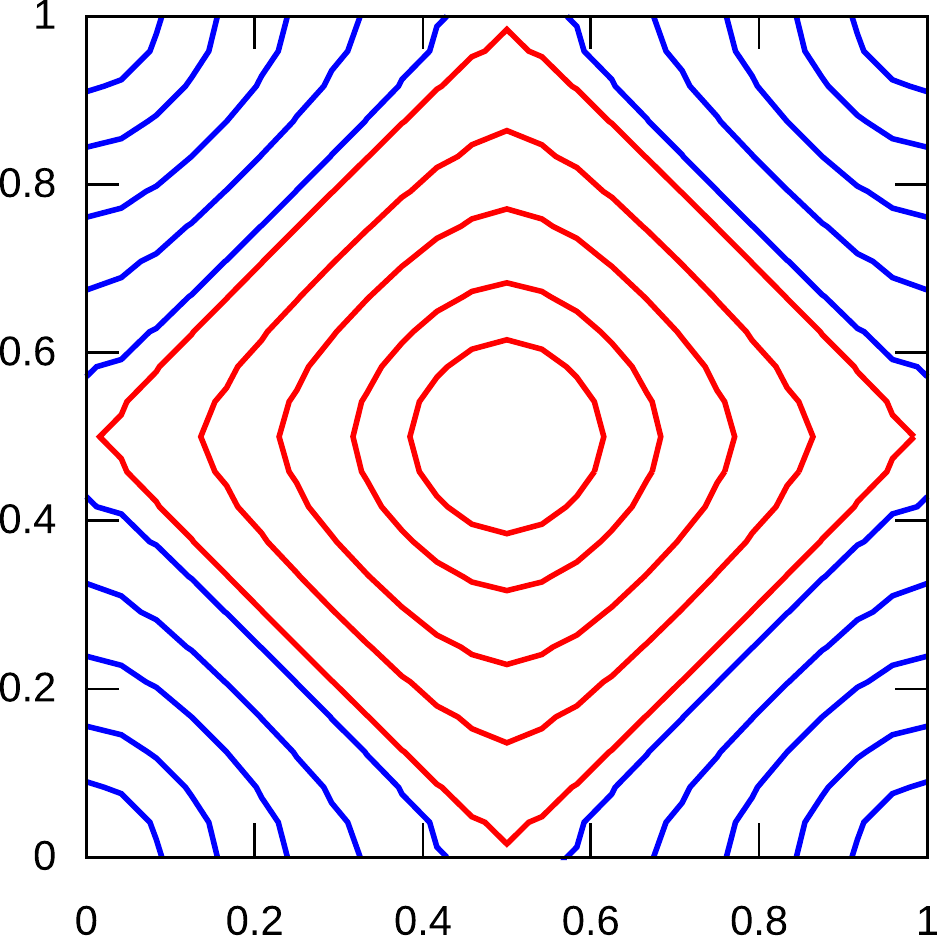} &
\includegraphics[width=.16\textwidth]{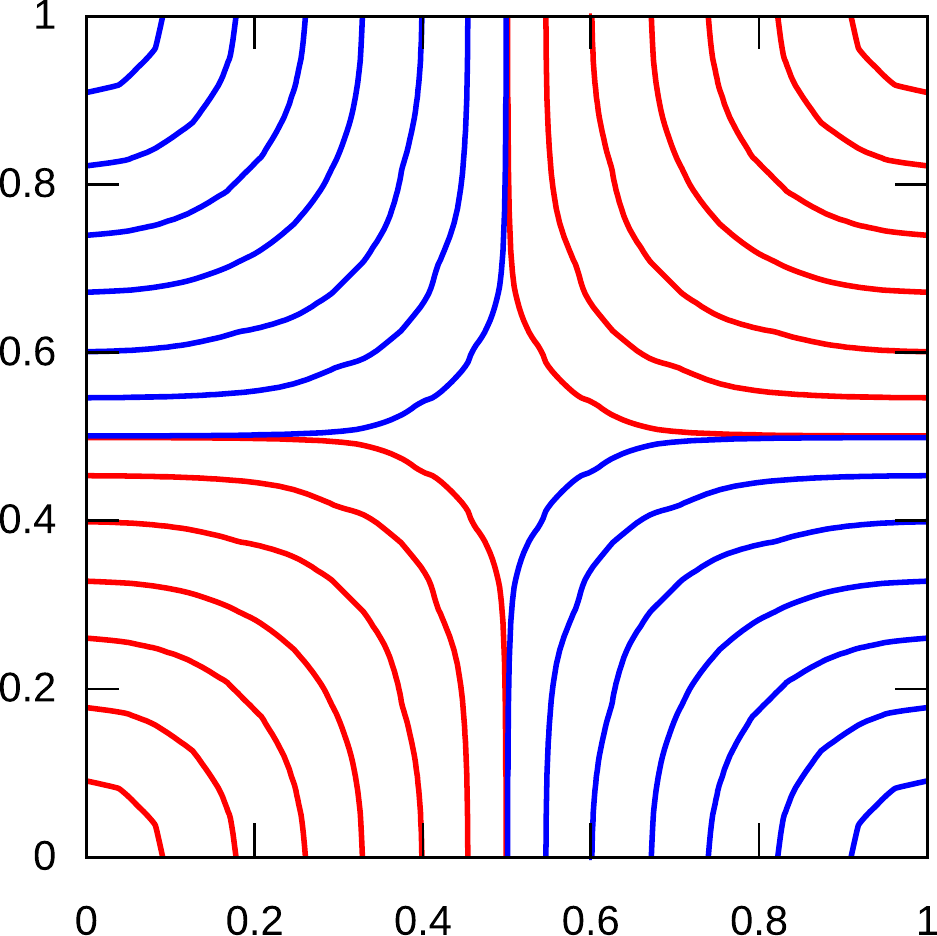} &
\includegraphics[width=.16\textwidth]{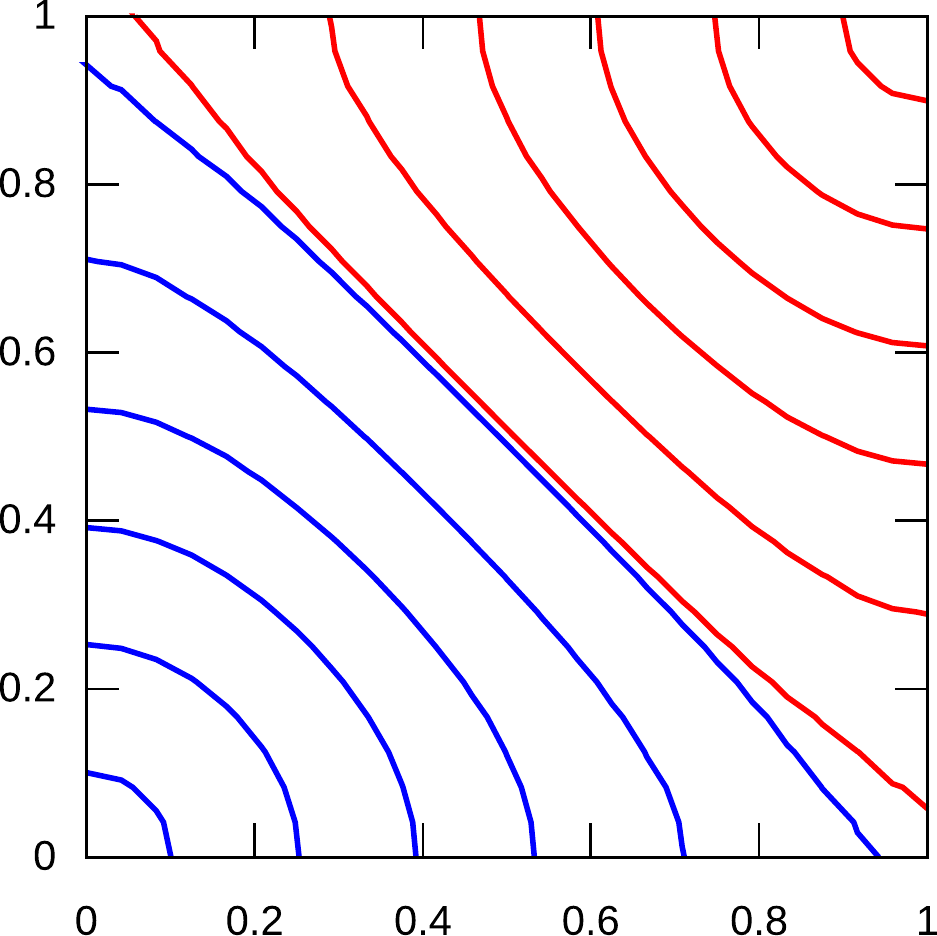} &
\includegraphics[width=.16\textwidth]{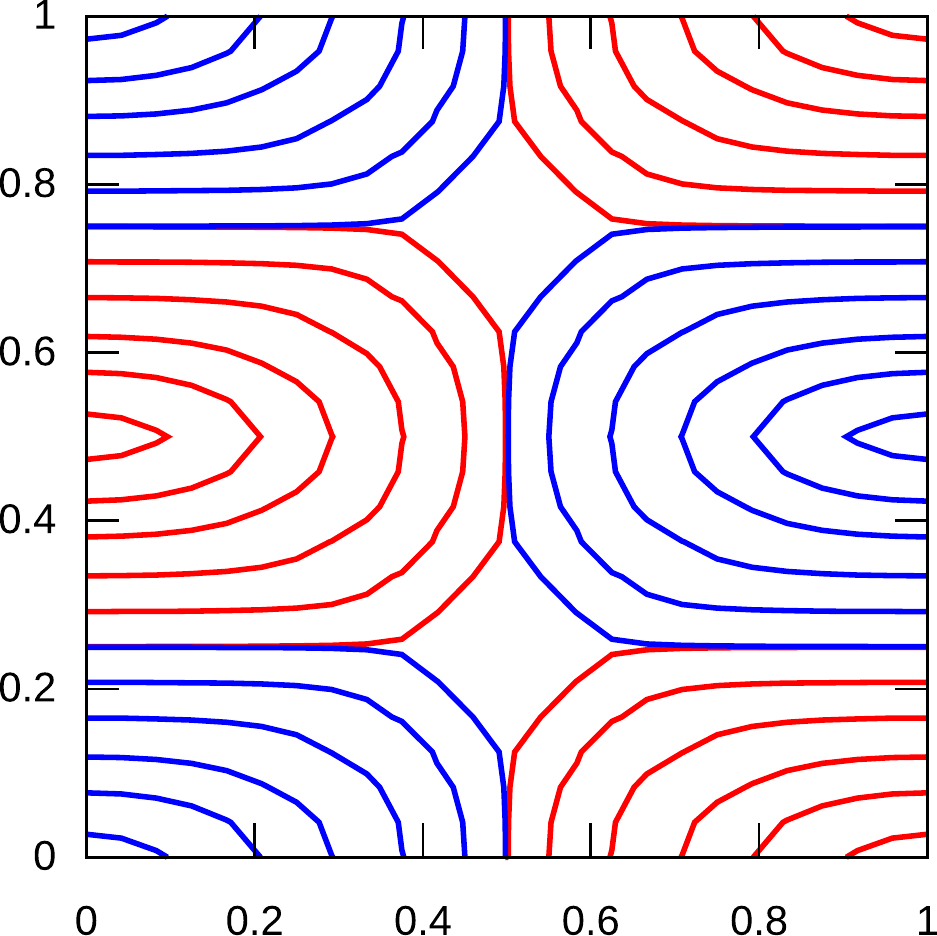} &
\includegraphics[width=.16\textwidth]{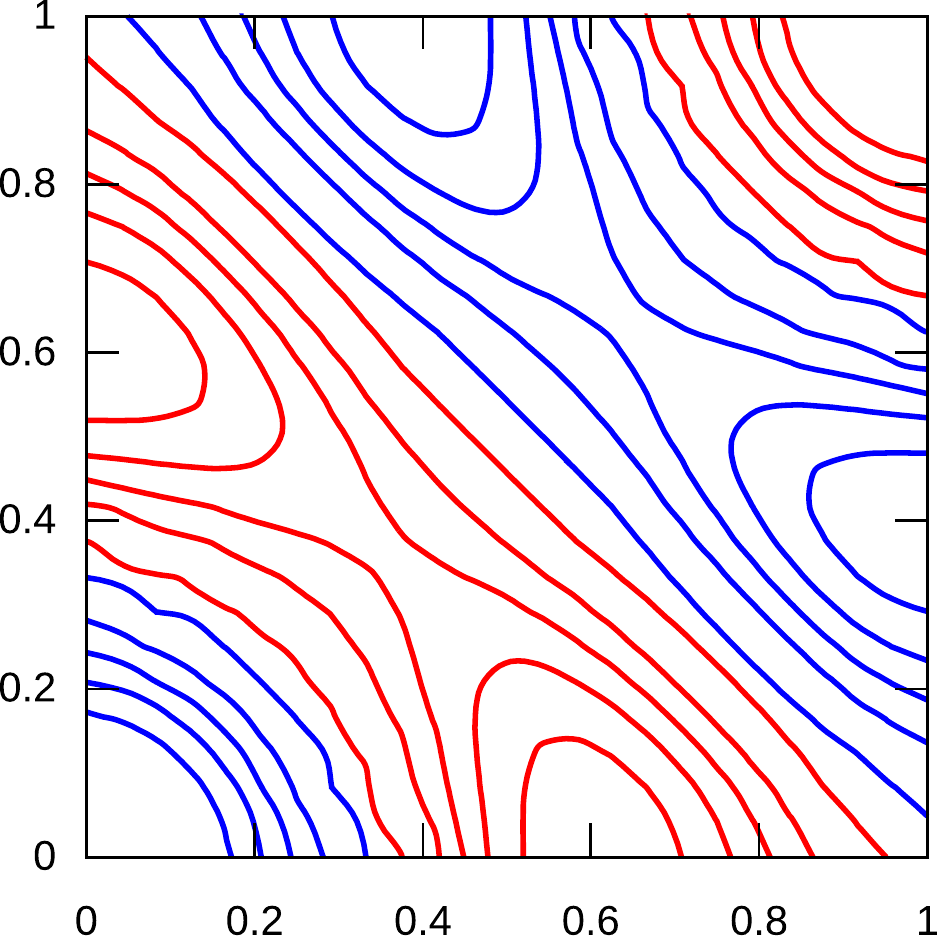} \\
\includegraphics[width=.16\textwidth]{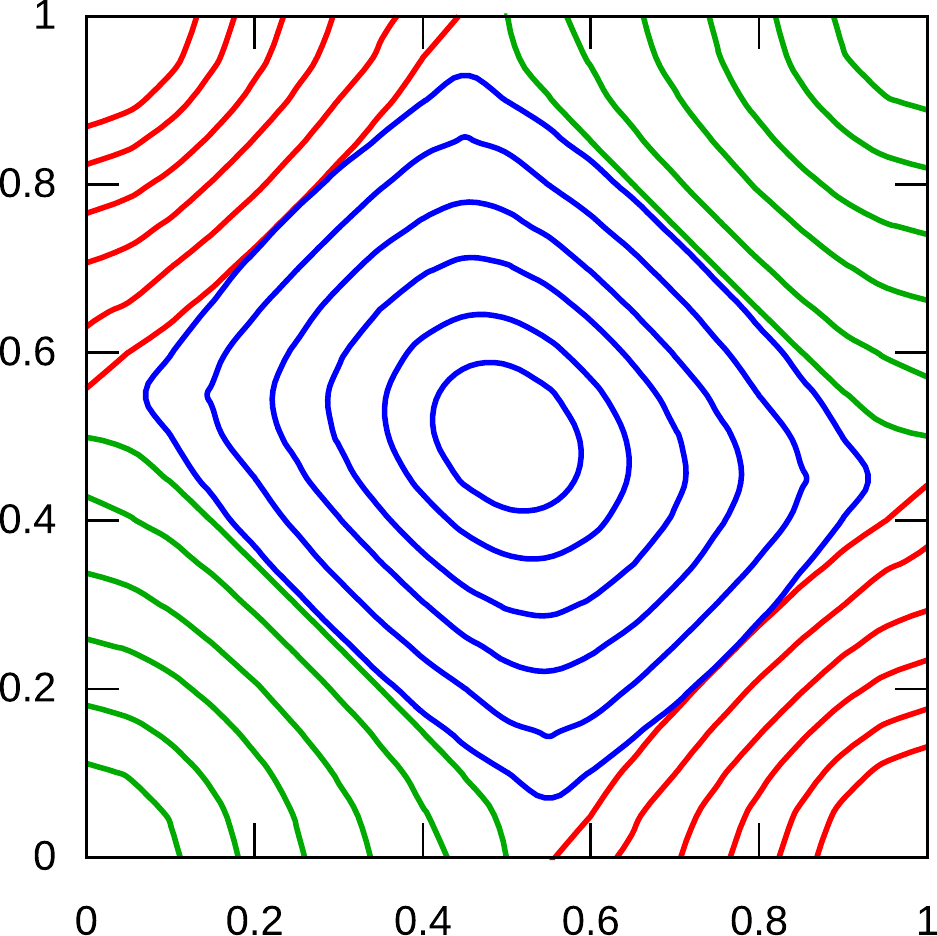} &
\includegraphics[width=.16\textwidth]{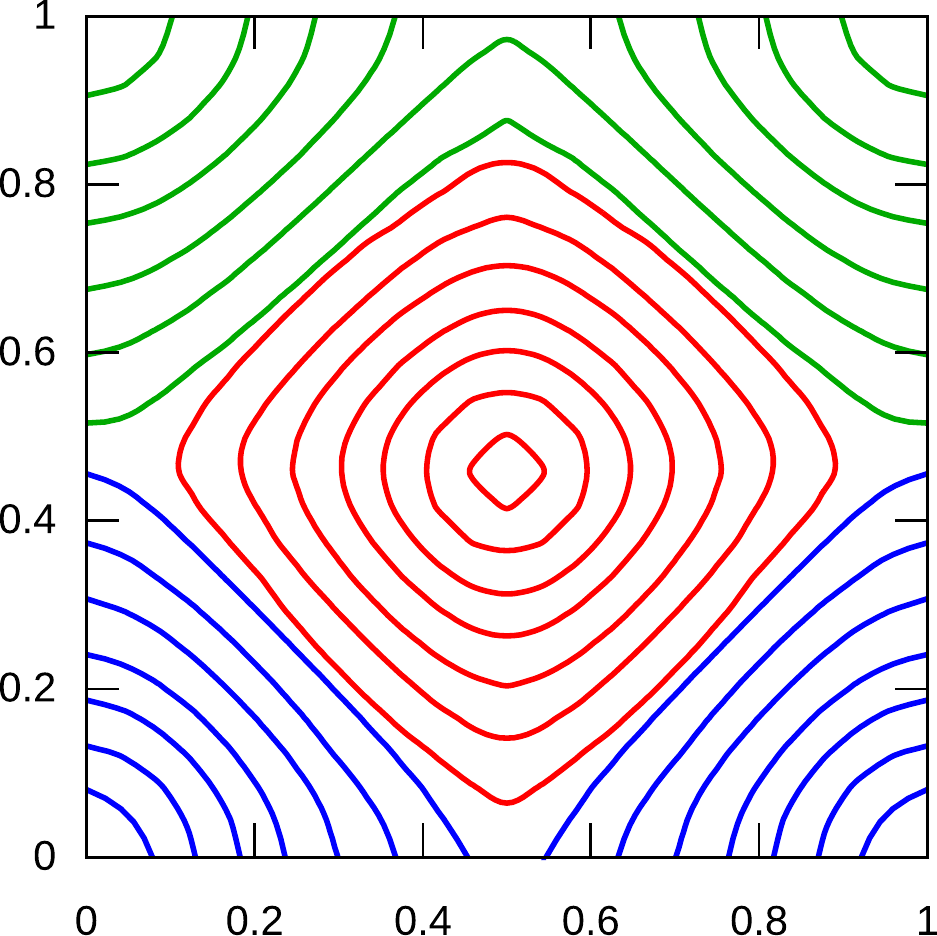} &
\includegraphics[width=.16\textwidth]{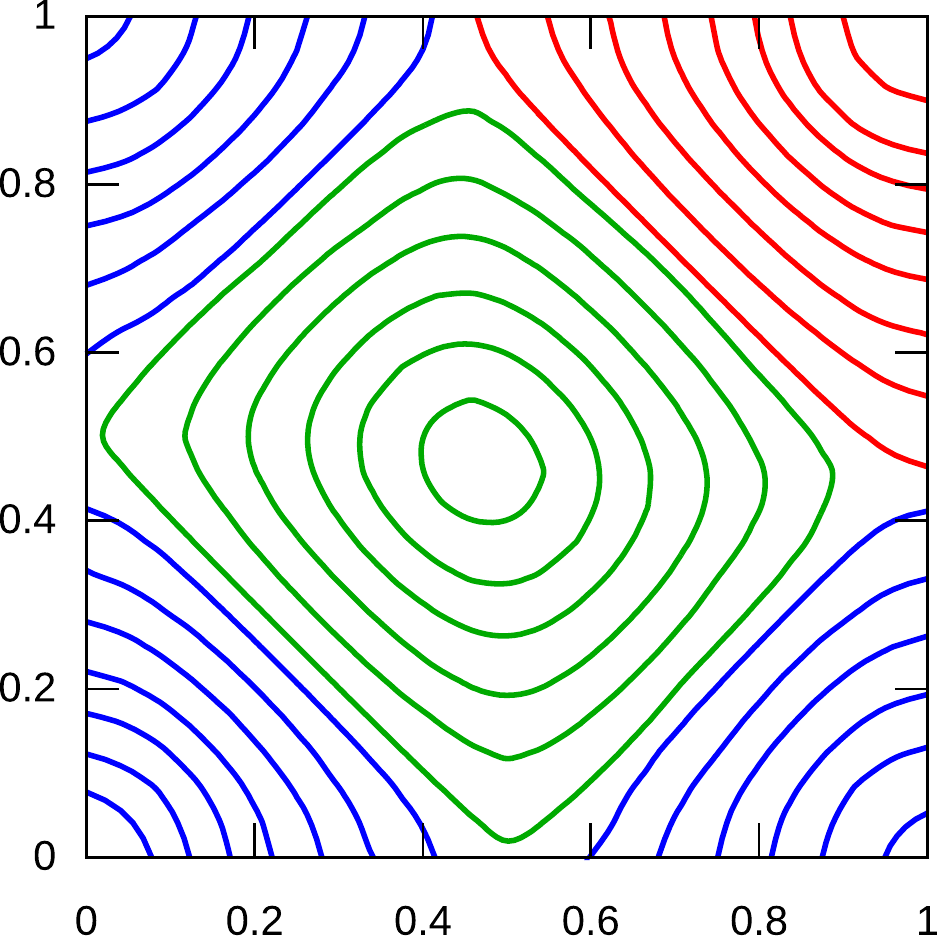} &
\includegraphics[width=.16\textwidth]{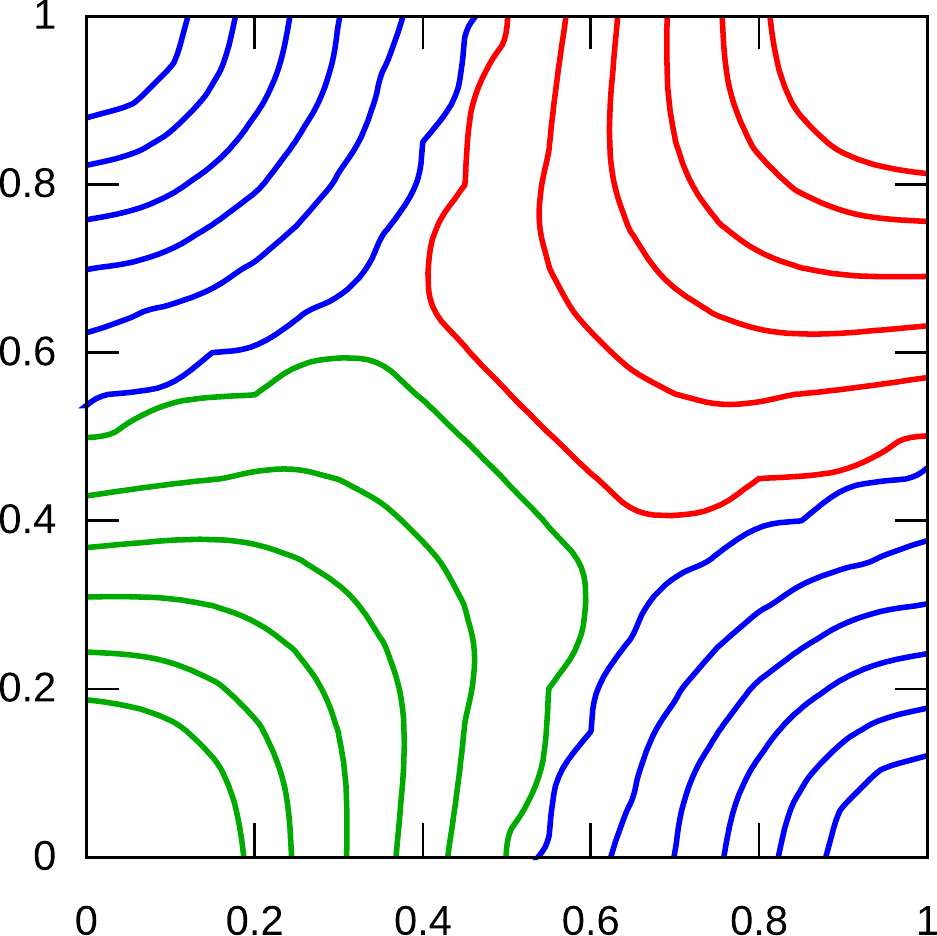} &
\includegraphics[width=.16\textwidth]{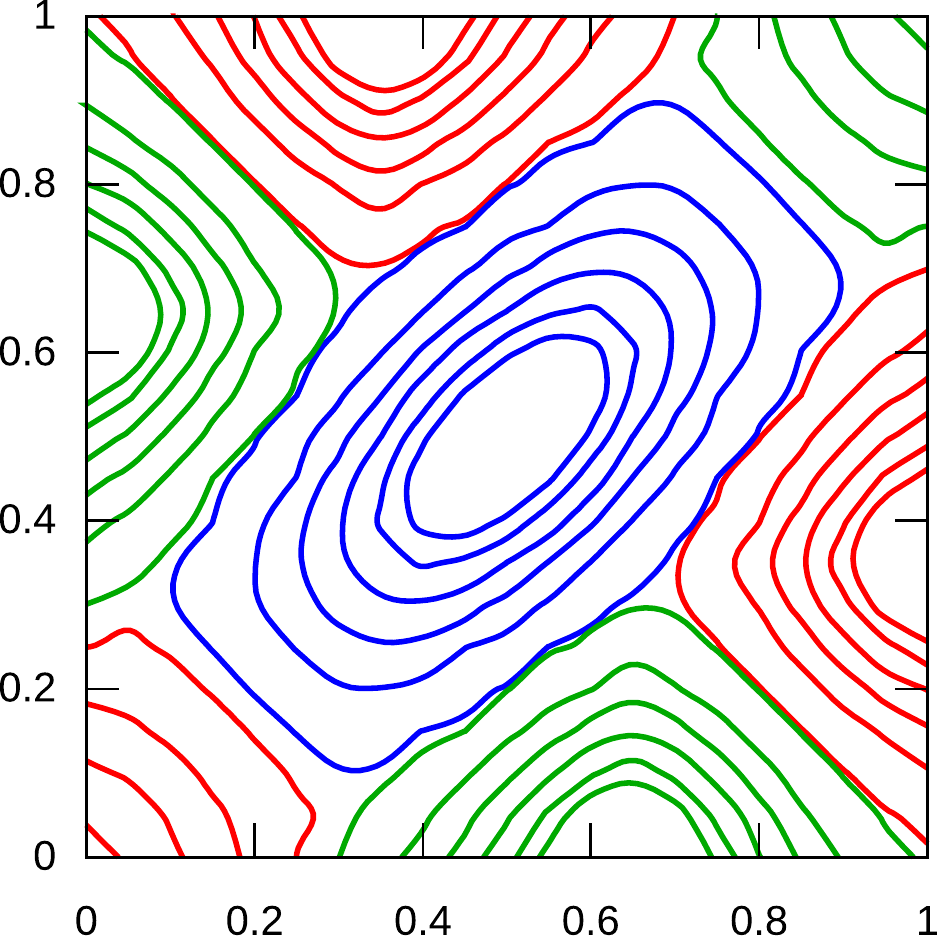} \\
\includegraphics[width=.16\textwidth]{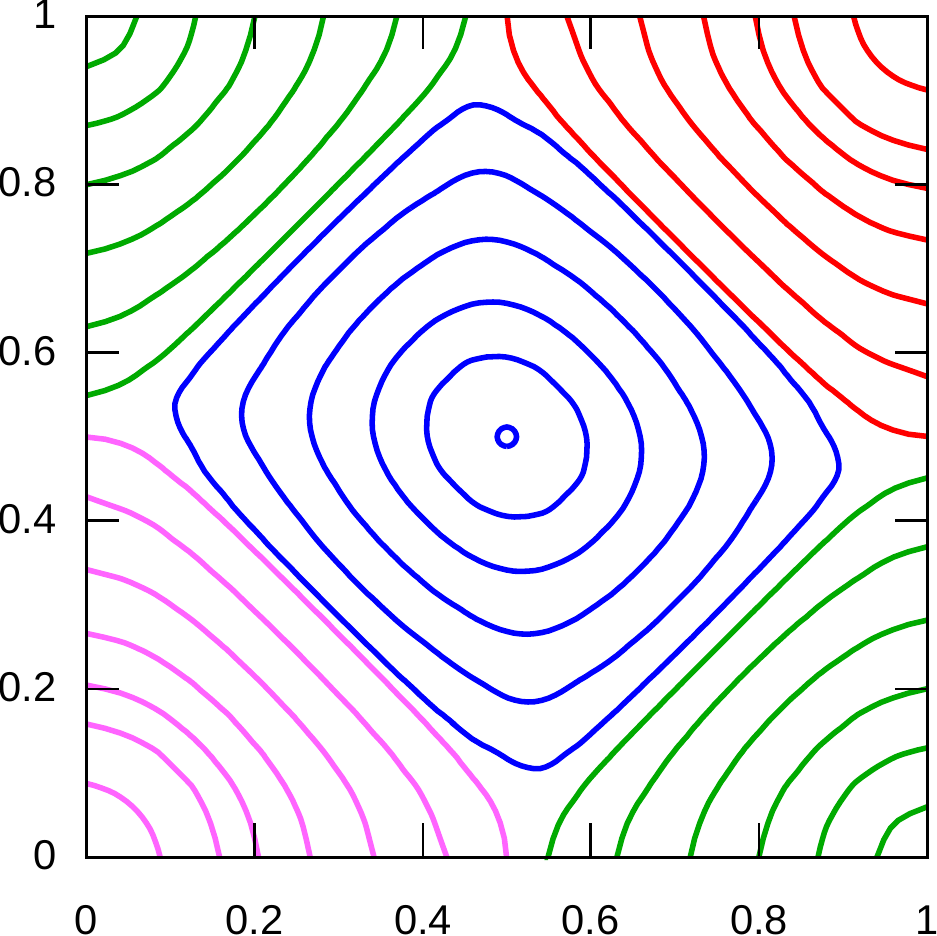} &
\includegraphics[width=.16\textwidth]{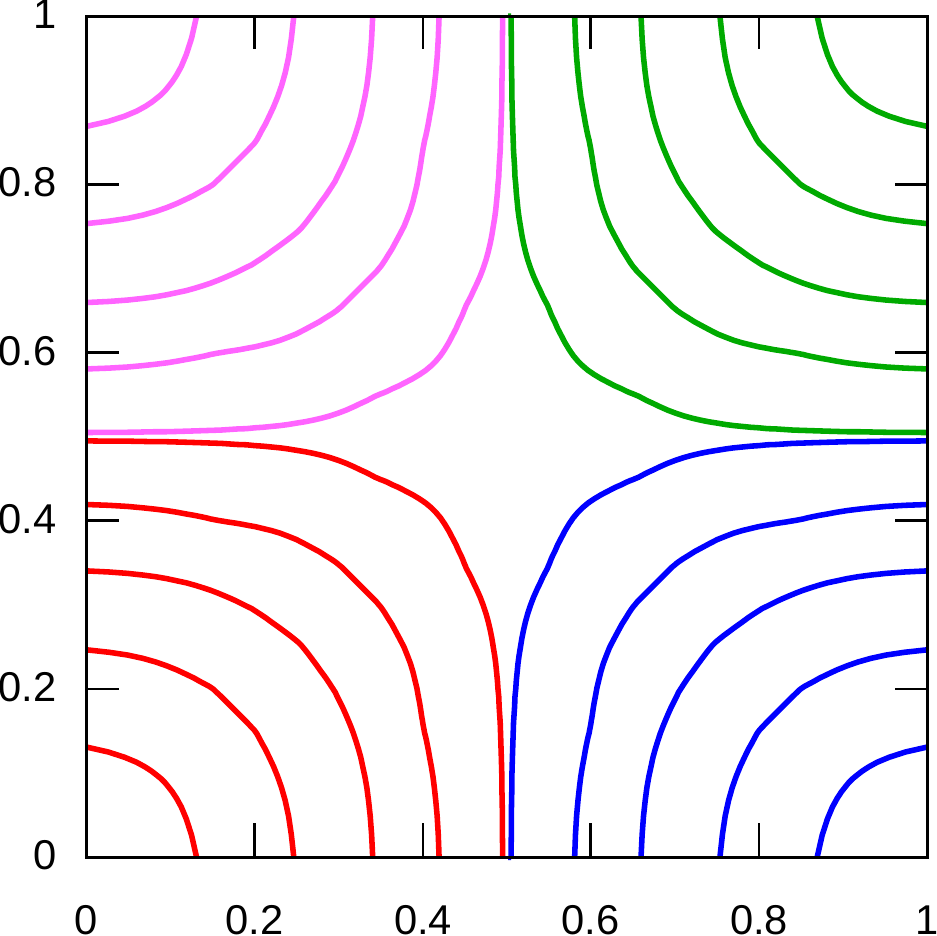} &
\includegraphics[width=.16\textwidth]{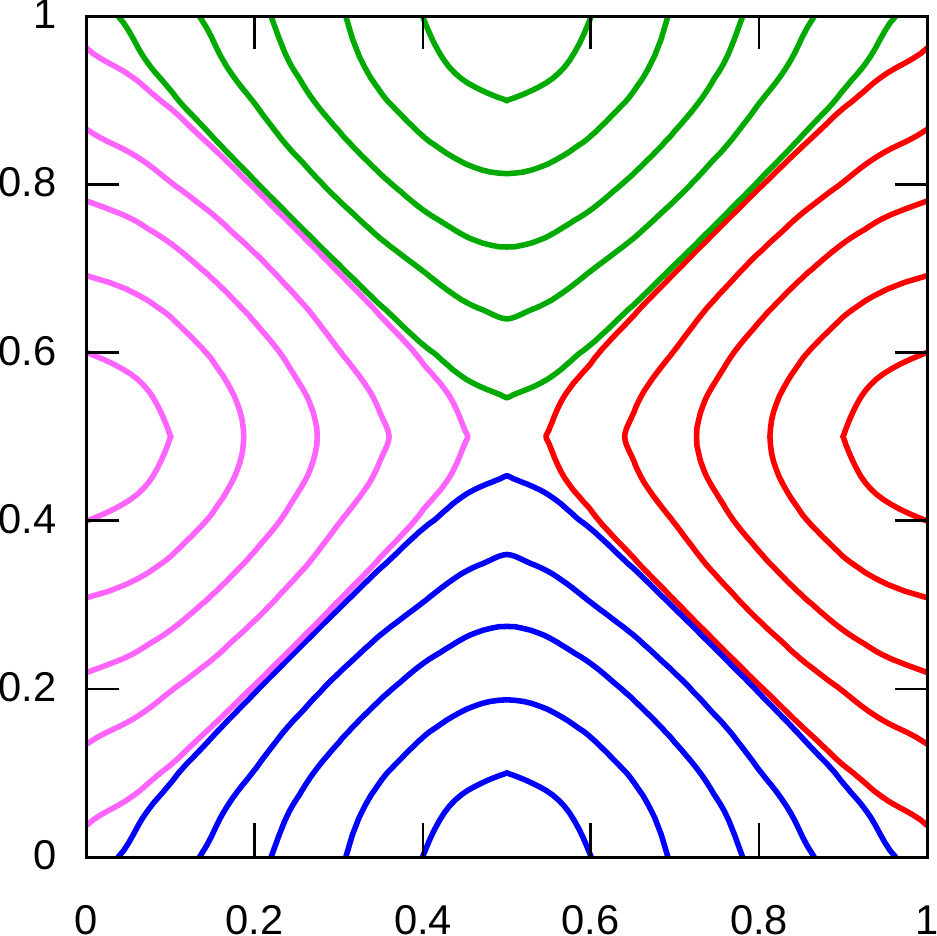} &
\includegraphics[width=.16\textwidth]{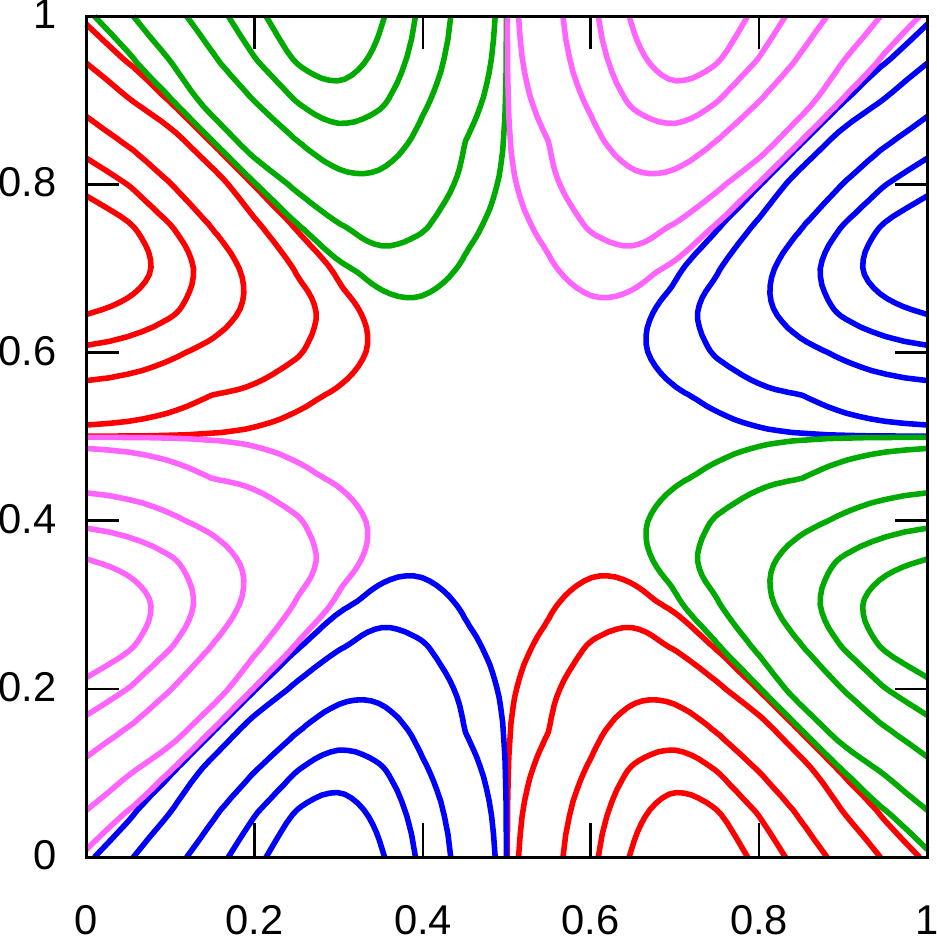} &
\includegraphics[width=.16\textwidth]{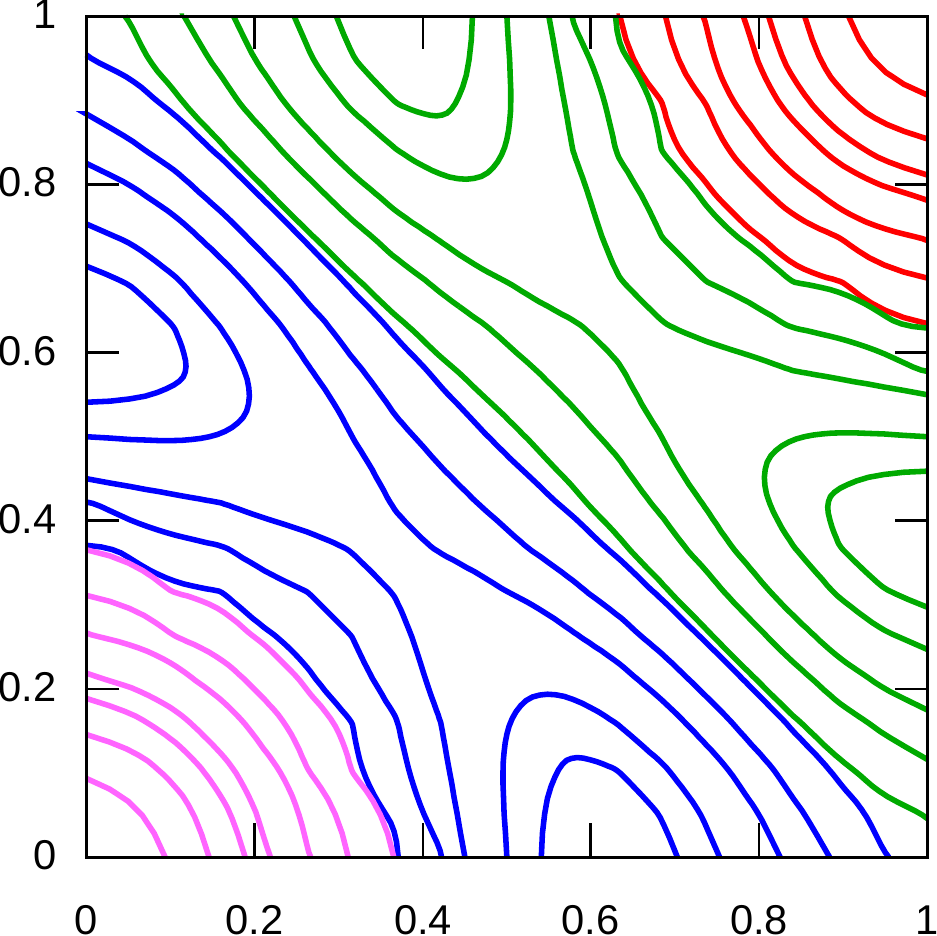} \\
\end{tabular}
\end{center}
\caption{Multi-population MFG mass densities.}
\label{mfgm2d-1}
\end{figure}

Here it is worth noting that we have almost no control on the outcome of the computation. Despite we tried to initialize
the mass densities and the correctors by means of Gaussian distributions supported in small balls at given points, the result is unpredictable.
Figure \ref{mfgm2d-1} shows a rich collection of solutions, corresponding to $P=2$ (top panels), $P=3$ (middle panels)
and $P=4$ (bottom panels) populations. We clearly see how the populations compete to share out all the domain.
\section*{Conclusions}\label{conclusions}
We presented a new approach for the numerical solution of ergodic problems involving Hamilton-Jacobi equations. 
The proposed Newton-like method for inconsistent nonlinear systems 
is able to solve first and second order nonlinear cell problems arising in very different contexts,
e.g., for scalar convex and nonconvex Hamiltonians, weakly coupled systems, dislocation dynamics and mean field games, also in the case of
more competing populations. A very large collection of numerical simulations shows the performance of the algorithm, including some experimental convergence
analysis. We reported both numerical results and computational times, in order to allow future comparison.

\section*{Acknowledgements}
The authors would like to thank F.J. Silva who, talking about mean fields games, pronounced the magic words ``Gauss-Newton method'', 
and also M. Cirant for fruitful discussions on the multi-population extension.


\begin{thebibliography}{1}

\bibitem{AP}
{\sc Y. Achdou and S. Patrizi},
{\em Homogenization of first order equations with $u/\epsilon$-periodic Hamiltonian: rate of convergence as $\epsilon \to 0$ and numerical approximation of the effective Hamiltonian},
Math. Models Methods Appl. Sci., 21 (2011), pp. 1317--1353.

\bibitem{AC}
{\sc Y. Achdou and I. Capuzzo Dolcetta},
{\em Mean field games: numerical methods},
SIAM J. Numer. Anal., 48 (2010), no. 3, pp. 1136--1162.

\bibitem{ACC}
{\sc Y. Achdou, F. Camilli and I. Capuzzo Dolcetta},
{\em Homogenization of Hamilton-Jacobi equations: Numerical Methods},
Math. Models Methods Appl. Sci., 18 (2008), pp. 1115--1143.

\bibitem{Ba}
{\sc N. Baca\"er},
{\em Convergence of numerical methods and parameter dependence of min-plus eigenvalue problems,
Frenkel-Kontorova models and homogenization of Hamilton-Jacobi equations},
Math. Model. Numer. Anal., 35 (2001), pp. 1185--1195.

\bibitem{BI}
{\sc A. Ben-Israel},
{\em A Newton-Raphson method for the solution of systems of equations},
J. Math. Anal. Appl., 15 (1966), pp. 243--252.

\bibitem{CCM16}
{\sc S. Cacace, F. Camilli and C. Marchi},
{\em A numerical method for Mean Field Games on networks},
preprint \texttt{http://arxiv.org/abs/1507.03731}

\bibitem{CCM}
{\sc S. Cacace, A. Chambolle and R. Monneau},
{\em A posteriori error estimates for the effective Hamiltonian of dislocation dynamics},
Numerische Mathematik, 121 (2012), no. 2, pp. 281--335.

\bibitem{CM}
{\sc F. Camilli and C. Marchi},
{\em Rates of convergence in periodic homogenization of fully nonlinear uniformly elliptic PDEs},
Nonlinearity 22 (2009), no. 6, pp. 1481--1498.

\bibitem{C15}
{\sc M. Cirant}, {\em Multi-population mean field games systems with Neumann boundary conditions},
J. Math. Pures Appl., 103 (2015), no. 5, pp. 1294--1315.

\bibitem{DZ}
{\sc A. Davini and M. Zavidovique},
{\em Aubry sets for weakly coupled systems of Hamilton-Jacobi equations},
SIAM J. Math. Anal. 46 (2014), pp. 3361--3389.

\bibitem{SPQR}
{\sc T. Davis},
{\em SuiteSparse},
\texttt{http://faculty.cse.tamu.edu/davis/suitesparse.html}

\bibitem{E}
{\sc L.C. Evans}, Some new PDE methods for weak KAM theory, Calc. Var. Partial Differ.
Equ. 17 (2003), no. 2, 159--177.

\bibitem{FIM}
{\sc N. Forcadel, C. Imbert and R.Monneau},
{\em Homogenization of some particle systems with two-body interactions and of the dislocation dynamics},
Discrete Contin. Dyn. Syst. 23 (2009), no. 3, pp. 785--826.

\bibitem{FR}
{\sc M. Falcone and M. Rorro},
{\em On a Variational Approximation of the Effective Hamiltonian}, in
Numerical Mathematics and Advanced Applications, Springer Berlin Heidelberg, 2008, pp. 719--726.

\bibitem{GHM}
{\sc M.-A. Ghorbel, P. Hoch and R. Monneau},
{\em A numerical study for the homogenization of one-dimensional models describing the motion of discolations},
Int. J. Comput. Sci. Math. 2 (2008), no. 1-2, pp. 28--52.

\bibitem{GWP}
{\sc M.J. Glencross, E.T. Wong and M. Planitz},
{\em Three slants on the generalised inverse},
Math. Gaz., 63 (1979), pp. 173--185.

\bibitem{G}
{\sc D. Gomes},
{\em A stochastic analogue of Aubry-Mather theory},
Nonlinearity, 15 (2002), pp. 581--603.

\bibitem{GO}
{\sc D. Gomes and A. Oberman},
{\em Computing the effective Hamiltonian using a variational approach},
SIAM J. Control Optim., 43 (2004), pp. 792--812.

\bibitem{LL}
{\sc J-M. Lasry and P-L. Lions},
{\em Mean field games},
Jpn. J. Math., 2 (2007), no. 1, pp. 229--260.

\bibitem{LPV}
{\sc P-L. Lions, G. Papanicolaou and S.R.S. Varadhan},
{\em Homogenization of Hamilton-Jacobi equations}, unpublished.

\bibitem{Q}
{\sc J. Qian},
{\em Two approximations for effective Hamiltonians arising from homogenization of Hamilton-Jacobi equations},
UCLA, Department of Mathematics, preprint, 2003.

\bibitem{R06}
{\sc M. Rorro},
{\em An approximation scheme for the effective Hamiltonian and applications},
Appl. Numer. Math., 56 (2006), no. 9, pp. 1238--1254.

\end{thebibliography}
\end{document}